\documentclass[10pt]{elsarticle-ppn}
\usepackage{fancyhdr}
\usepackage{marginnote}
\usepackage{makecell} 
\usepackage{graphicx} 
\usepackage[shortlabels]{enumitem}
\usepackage{amsmath}
\usepackage{tikz}
\usepackage{frcursive}
\usepackage{mathtools}

\usepackage{amsfonts}
\usepackage{mathrsfs}
\usepackage{amsthm}
\usepackage{amssymb}
\usepackage{color}
\usepackage{dsfont}
\usepackage{geometry}
\newgeometry{tmargin=2cm, bmargin=3cm, lmargin=1.7cm, rmargin=1.7cm}
\usepackage{enumitem, hyperref}
\usepackage{tcolorbox}

\DeclareMathAlphabet{\mathpzc}{OT1}{pzc}{m}{it}

\newtheorem{theo}{\bf Theorem}

\newtheorem{coro}{\bf Corollary}[section]
\newtheorem{lem}[coro]{\bf Lemma}
\newtheorem{rem}[coro]{\bf Remark}
\newtheorem{defi}[coro]{\bf Definition}

\newtheorem{prop}[coro]{\bf Proposition}

\def\avenorm#1{\mathchoice%
          {\mathop{\kern 0.2em\vrule width 0.6em height 0.69678ex depth -0.58065ex
                  \kern -0.545em \|{#1}\|}}%
          {\mathop{\kern 0.1em\vrule width 0.5em height 0.69678ex depth -0.60387ex
                  \kern -0.495em \|{#1}\|}}%
          {\mathop{\kern 0.1em\vrule width 0.5em height 0.69678ex depth -0.60387ex
                  \kern -0.495em \|{#1}\|}}%
          {\mathop{\kern 0.1em\vrule width 0.5em height 0.69678ex depth -0.60387ex
                  \kern -0.495em \|{#1}\|}}}

\newcommand{\barint}{
         \rule[.036in]{.12in}{.009in}\kern-.16in
          \displaystyle\int  }

\def\aa{{\mathfrak{a}}}
\def\fb{{\mathfrak{b}}}
\def\R{{\mathbb{R}}}

\def\ba{{\mathbf{a}}}
\def\bb{{\mathbf{b}}}
\def\rn{{\mathbb{R}^{N}}}

\def\w{\widetilde}

\def\rp{{[0,\infty )}}

\def\ve{{\varepsilon}}
\def\vr{{\varrho}}
\def\vt{{\vartheta}}

\newcommand{\wt}{\widetilde}

\newcommand{\vp}{\varphi}

\newcommand{\dv}{\mathrm{div}}
\newcommand{\V}{V}
\newcommand{\VD}{{V_D}}
\newcommand{\vo}{{{v_0}}}
\newcommand{\uo}{{{u_0}}}

\def\d{\,{\rm d}}
\newcommand{\ds}{\,{\rm d}s}
\newcommand{\dx}{\,{\rm d}x}
\newcommand{\dy}{\,{\rm d}y}

\newcommand{\dt}{\,{\rm d}t}
\newcommand{\dr}{\,{\rm d}r}
\newcommand{\dD}{\,{\rm d}D}
\newcommand{\dtau}{\,{\rm d}\tau}
\newcommand{\rd}{{\rm d}}

\def\LL{\mathrm{L}}

\makeatletter
\def\namedlabel#1#2{\begingroup
    #2%
    \def\@currentlabel{#2}%
    \phantomsection\label{#1}\endgroup
}
\makeatother

\newcommand{\RR}{\mathbb{R}}

\newcommand{\B}{\mathcal{B}}
\newcommand{\Mo}{M}

\newcommand{\linE}{{\mathsf{E}}}
\newcommand{\linI}{{\mathsf{I}}}
\newcommand{\cE}{{\mathcal{E}}}
\newcommand{\cI}{{\mathcal{I}}}

\newcommand{\X}{\mathcal{X}_p}
\newcommand\RN{\mathbb{R}^N}

\numberwithin{equation}{section}

\definecolor{darkblue}{rgb}{0.05, .05, .65}
\definecolor{darkgreen}{rgb}{0.05, .70, .05}
\definecolor{darkred}{rgb}{0.8,0,0}

\begin{document}

\begin{frontmatter}

\title{\textbf{Refined asymptotics for the Cauchy problem\\ for the fast $p$-Laplace evolution equation}}

\author[1]{Matteo Bonforte}\ead{matteo.bonforte@uam.es}
\author[2]{Iwona Chlebicka\corref{mycorrespondingauthor}}
\cortext[mycorrespondingauthor]{Corresponding author}
\ead{i.chlebicka@mimuw.edu.pl}
\author[3]{Nikita Simonov}\ead{nikita.simonov@sorbonne-universite.fr}

\address[1]{Departamento de Matem\'a{}ticas, Universidad Aut\'o{}noma de Madrid, ICMAT -- Instituto de Ciencias
Matem\'a{}ticas, CSIC-UAM-UC3M-UCM, Campus de Cantoblanco, 28049 Madrid, Spain}
\address[2]{Institute of Applied Mathematics and Mechanics,
Faculty of Mathematics, Informatics and Mechanics,
University of Warsaw, ul.~Banacha 2, 02-097 Warsaw, Poland
}
\address[3]{Laboratoire Jacques Louis Lions (CNRS UMR n$^\circ$~7598), Sorbonne Universit\'e, 4 place Jussieu, 75005 Paris, France}

\begin{abstract} 

Our focus is on the fast diffusion equation driven by the $p$-Laplacian operator, that is $\partial_t u=\Delta_p u$ with $1<p<2$, posed in the whole space $\mathbb{R}^N$, $N\geq 2$. The nonnegative solutions, in rescaled variables, are expected to converge in time toward a stationary profile. While such convergence had been previously established for $p$ close to $2$, no quantitative rates were known, and the asymptotic behaviour remained poorly understood across the full fast diffusion range. In fact, the long time behaviour of solutions to the $p$-Laplace Cauchy problem drastically change in different subranges of the $p$. Some of them are analysed here for the first time.

In this work, we provide the convergence rates for nonnegative, integrable solutions in the so-called good fast diffusion range, $p_c=\tfrac{2N}{N+1} <p<2$, where mass is conserved.  We prove that solutions converge to a self-similar profile with matching mass, with explicit rates measured in relative error. Our constructive proof is based on a new entropy method that remains effective even when the entropy is not displacement convex, where optimal transport techniques fail.

In the very fast diffusion range $1<p<p_c$, we give the first asymptotic analysis near the extinction time. This regime poses additional challenges: mass is not conserved, solutions vanish in finite time, and no fundamental solutions exist. We found new critical exponents -- especially in high dimensions -- that give rise to markedly different qualitative behaviour depending on the value of $p$.

We also establish convergence rates for the gradients of radial solutions in the good fast diffusion range, again measured in relative error. Finally, we analyze the structural properties required for the entropy method to apply, thereby opening a broader investigation into the basin of attraction of Barenblatt-type profiles, particularly in the singular case of $p$ close to $1$.
\end{abstract}


\end{frontmatter}

\newpage

\setcounter{tocdepth}{1}

\tableofcontents

\section{Introduction}
In this paper we are interested in the long-time behaviour of solutions to the following Cauchy problem
\begin{equation}\tag{CPLE}\label{CPLE}
\begin{cases}
\partial_t u=\dv\left(|\nabla u|^{p-2}\,\nabla u\right)\quad  x\in\rn,\quad t>0\,,\\
u(0, x)=\uo(x)\ge0\,,\quad x\in\rn\,,
\end{cases}
\end{equation}
where the exponent $1<p<2$. The above equation is called the \emph{p-Laplace evolution} equation, it is a nonlinear, gradient-driven diffusion equation which is singular in the considered regime $1<p<2$, linear when $p=2$, and degenerate for $p>2$. Nonlinear evolution equations involving the $p$-Laplace operator $\Delta_p u=\dv\left(|\nabla u|^{p-2}\,\nabla u\right)$ have attracted the attention of the mathematical community for over a half of century. The interest come from both the preeminent role of the power-growth operators in the modelling of fluid dynamics~\cite{Lad} and intrinsic mathematical challenges, see e.g.~\cite{DiBenBook,KaminVazquez1988,Carrillo-Jungel-Markowich-Toscani-Unterreiter,KuMi-guide} and references therein, and the very recent probabilistic interpretation in terms of a game theoretical approach~\cite{delTeso2024,Lewicka2017,Lewicka2021,Lewicka2022, Lewicka2020book}.

The existence and regularity of solutions has been widely studied, see for instance~\cite{JLVSmoothing,Vazquez2007,DiBenBook}. Indeed, it is known that, under the assumption $0\leq u_0\in\mathrm{L}^1_{\mathrm{loc}}(\mathbb{R}^N)$, the (unique) solution to \eqref{CPLE} is a continuous curve in $\LL^1_{\mathrm{loc}}(\rn)$, that is $u\in C\left(\left[0, \infty\right), \LL^1_{\mathrm{loc}}(\rn)\right)$, it remains nonnegative, and  $u(t,\cdot)\in C^{1,\alpha}(\rn)$ for some $\alpha\in(0,1)$ (when $p$ is close to one, one also need to ask higher integrability of $u_0$), see Section~\ref{sec:prelim} for more details.

The long time behaviour of~\eqref{CPLE} depends both on the spatial dimension $N$ and the exponent $p$. It is natural, however, to divide analysis into three distinct cases: $1<p<2$, $p=2$, and $p>2$.  When $p=2$,~\eqref{CPLE} reduces to the classical \emph{heat equation}. In the present setting (nonnegative, integrable initial data), the asymptotic behaviour of solutions is well described by the Gaussian profile. We refer to~\cite{Vazquez2017,J_ngel_2016} for detailed results on convergence rates and the refined asymptotic analysis in various contexts. When $p>2$, the dynamics differs significantly from the linear case. Notably, solutions exhibit \emph{finite speed of propagation}: compactly supported initial data generate solutions that remain compactly supported for all times. In this case, the long-time behaviour of solutions is described by the fundamental solution, commonly known as Barenblatt, or Barenblatt-Pattle solution, which is itself compactly supported for all times. Compared to the linear case, the available results are more limited. For a discussion of convergence rates and fine asymptotic behaviour, we refer to~\cite{KaminVazquez1988,Vespri2016, Agueh2003, Agueh2008, DelPino2002, DelPino2003}.  As final comment, we also remark that, in the cases $p=2$ and $p>2$, \emph{conservation of mass} holds for non-negative $\mathrm{L}^1(\mathbb{R}^N)$ solutions.

The aim of this paper is to present new results and give a broad perspective on the long-time behaviour for~\eqref{CPLE} in the parameter range $1<p<2$. In this regime, as in the linear case, solutions exhibit \emph{infinite speed of propagation}. However, in contrast with the linear case, they develop \emph{fat tails}: compactly supported initial data lead to solutions with algebraic spatial decay (power law), rather than a Gaussian profile typical of the linear scenario. The second main distinction from the case $p\ge2$ is the failure \emph{conservation of mass} when $p$ is close to $1$. Indeed,  for every $t>0$ it holds
\begin{equation}\label{p_c}
\int_{\rn}u(t,x)\dx = \int_{\rn}\uo(x)\dx\quad\mbox{only if}\quad p\ge  p_c:=\frac{2N}{N+1}\,,
\end{equation}
see for instance~\cite{Vazquez:2023aa}.
On the one hand, when mass is preserved, solutions starting from non-negative initial data become positive for all times, and we analyze the asymptotic behaviour as $t\to\infty$, analogously to the case $p\ge 2$. On the other hand, the failure of mass conservation for $1<p<p_c$ significantly complicates the asymptotic analysis of~\eqref{CPLE}, especially when compared to the better-understood regime $p_c<p<2$. Even less is known in the critical case $p=p_c$.  The loss of mass in the very fast diffusion range allows a wide class of solutions to extinguish in finite time: the asymptotic behaviour in this case needs to be analyzed close to the extinction time, which is not explicit and depends on the initial datum.  For this reasons, among others discussed below, it is natural to distinguish between the \emph{good fast diffusion range}, defined by $p_c < p < 2$, and the \emph{very fast diffusion range}, where $1 < p \leq p_c$. In the very fast diffusion range the spatial dimension $N$ also plays a significant role in the asymptotic analysis. We identify new critical exponents at which previously unrecognized qualitative changes in the behavior of solutions emerge. Remarkably, they appear in high dimensions only, namely when $N\ge6$. Let us begin by describing the good fast diffusion range, where the behaviour of solutions is better understood.

\subsection{The good fast diffusion range: \texorpdfstring{$p_c<p<2$}{pc< p< 2} and \texorpdfstring{$N\geq2$}{N>=2}}
In this range the behaviour of solutions to~\eqref{CPLE} for large times is described by the means of the \emph{fundamental solution} (also called the \emph{Barenblatt} or \emph{Barenblatt-Pattle} solution) defined as
\begin{equation}\label{beta}
\B_M(t,x):= t^\frac{1}{2-p}\,\left[b_1\, t^\frac{\beta\,p}{p-1}\,M^\frac{\beta\,p\,(p-2)}{p-1} + b_2\, |x|^\frac{p}{p-1}\right]^\frac{1-p}{2-p}\quad\mbox{where}\quad\beta=\frac{1}{p-N(2-p)}
\end{equation}
and where $M$ represents the (conserved) mass of $\B_M$, $b_1$ and $b_2$ are given numerical constants (see
 \eqref{b1} and \eqref{b2} for their definitions).  The Barenblatt solution $\B_M$ takes a Dirac delta $M\,\delta_0$ as its initial datum (in the sense of distribution). While it is probably not so evident from the expression in~\eqref{beta}, $\B_M$ has a \emph{self-similar} form which can be understood through the profile $\VD:\rn\to\rp$ and the following formula
\begin{equation}\label{vD}
\VD(y):=\left(D+\tfrac{2-p}{p}|y|^\frac{p}{p-1}\right)^{-\frac{p-1}{2-p}}\,.
\end{equation}
Having $R_T:[0,\infty)\to[0,\infty)$ defined as
\begin{equation}\label{R(t)}
    R_T(t):=\left(\tfrac{t+T}{\beta}\right)^{\beta}\,,
\end{equation}
we have
\begin{equation}\label{Barenblatt-via-vD}
\B_M(t+T, x)=R^{-N}_T(t)\VD\left(\frac{x}{R_T(t)}\right)\qquad\mbox{and}\qquad D:=\beta^{N\beta\frac{2-p}{p-1}}\,\frac{b_1}{M^{p\beta\frac{2-p}{p-1}}}
\,.
\end{equation}
For non-negative initial data $\uo\in \LL^1(\rn)$, solutions to~\eqref{CPLE} \emph{relax to self-similarity} (see~\cite{KaminVazquez1988})  and the precise result can be stated as follows
\begin{equation}\label{mild.convergece.intro}
\|u(t,\cdot) - \B_M(t,\cdot)\|_{\LL^1(\rn)}\rightarrow0\quad\mbox{and}\quad t^{N\beta}\|u(t,\cdot) - \B_M(t,\cdot)\|_{\LL^\infty(\rn)}\rightarrow0\quad\mbox{as}\quad t\rightarrow\infty\,,
\end{equation}
where $M=\int_{\rn}\uo(x)\dx$, and the factor $t^{N\beta}$ in front of the $\LL^\infty$-norm is necessary to get a meaningful result, since the $\LL^\infty$ norm of solutions decays in time as $t^{-N\beta}$. By interpolation, similar results can be obtained for $\LL^q$-norms, for $1< q<\infty$. We notice that, in order to speak of results as~\eqref{mild.convergece.intro}, we shall often use the term \textit{convergence}, even if the Barenblatt function $\B_M$ is not a stationary profile.
However, it is possible, and sometimes very useful, to perform a change of variable in order to transform $\B_M(t+T, x)$ into the stationary profile $\VD(y)$ and, henceforth, rescale equation~\eqref{CPLE} into a nonlinear \emph{Fokker-Planck}  type  equation, see~\cite{J_ngel_2016}. Indeed, if
\begin{equation}\label{eq:v}\tag{R-CPLE}
v(\tau, y)={R_T^N(t)}u(t,x)\,,\quad\mbox{then}\quad
\left\{\begin{array}{ccl}
\partial_\tau v(\tau,y)& =&\dv_y\Big(|\nabla v(\tau,y)|^{p-2}\nabla v(\tau,y)+y \,v(\tau,y)\Big)\,,\\
v(x,0)& =&\vo(x)\,.
\end{array}\right.
\end{equation}
We notice that the initial datum $u_0$ is transformed accordingly to the formula $v_0(y)={R_T^N(0)}u_0(x)$, whereas
\begin{equation}
\label{change:rho-u}
\tau:=\ln \frac{R_T(t)}{R_T(0)}\,\qquad \text{and}\qquad y:=\frac{x}{R_T(t)}\,.
\end{equation}
The main advantage is that now the problem~\eqref{eq:v} has a family of stationary solutions $\VD$ given by \eqref{vD}\,. Indeed a simple computation shows that $|\nabla \VD(y)|^{p-2} \nabla \VD(y) = - y\, \VD(y)$, for any $y\in\rn$. We notice that at least when $p_c<p<2$, the parameter $D$ is completely determined by the mass of the function $\VD$.  In the present regime, most of the time we shall choose $T=\beta$ in \eqref{R(t)}, which sets the initial datum $\vo=\uo$ (since $R_\beta(0)=1$) and  makes most of the computations simpler. Namely, we pick
\begin{equation*}
R_\beta(t)=\big(1+\tfrac t\beta\big)^\beta\,,\quad\mbox{such that}\quad \VD(y)= R_\beta^N(t)\, \B_M(t+\beta, x)\,.
\end{equation*}
We shall call the functions $\VD$ the \emph{stationary Barenblatt profiles} or simply \emph{Barenblatt profiles} when no confusion arises. We also notice that, among all the non-stationary Barenblatt solutions defined in~\eqref{beta}, only $\B_M(t+\beta, x)$ is transformed to the stationary one through the change of variables defined above, the rest of functions from the family $\B_M(t+T, \cdot)$ do not become stationary after this change of variables. Since the equation (in the original variables) is translation invariant, there is no loss of generality in assuming the profiles to be centered at the origin throughout the paper.

Convergence results are much better understood in these new variables, for instance~\eqref{mild.convergece.intro} becomes
\begin{equation}\label{mild-convergence-R}
  \|v(\tau,\cdot) - \VD(\cdot) \|_{\LL^1(\rn)}\rightarrow0\quad\mbox{and}\quad \|v(\tau,\cdot) - \VD(\cdot)\|_{\LL^\infty(\rn)}\rightarrow0\quad\mbox{as}\quad \tau\rightarrow\infty\,.
\end{equation}
As a consequence of the change of variables, there is no factor $t^{N\beta}$ in front of the $\LL^\infty$ norm and now $v(\tau)$ converges to a stationary profile as $\tau\to\infty$. It is well known that, without any additional assumptions, such results are sharp with respect to the strength of the norm of convergence. At the same time, in the good fast diffusion range, no convergence rates can be established without imposing further hypotheses.
Indeed, counterexamples can be constructed using similar techniques to those in~\cite{Bonforte2020}. However, results as in~\eqref{mild-convergence-R} do not take into account neither the tail behaviour of the Barenblatt nor of the solution $v$ itself, and one may ask whether we can obtain a finer description of the tail behaviour for solutions to~\eqref{eq:v}. This was done in~\cite{Bonforte2020b} where solution with the same polynomial tail behaviour of the Barenblatt profile have been completely characterized.  Indeed, in~\cite[Theorem 1.1]{Bonforte2020b}, a stronger convergence result was proven with a sharp description of the tail behaviour of the initial datum for solutions to~\eqref{eq:v}.  The main result of that paper is the characterization of the \emph{uniform convergence in relative error} (UCRE), that is, for $N\ge1$ and $p_c<p<2$, we have
\begin{equation}\label{UCRE/Xp}
\left\|\frac{v(\tau,\cdot) - \VD(\cdot)}{\VD(\cdot)}\right\|_{\LL^\infty(\rn)}\xrightarrow[{\tau\to +\infty}]{}0\quad\mbox{if and only if}\quad \|\vo\|_{\X}:=\sup_{R>0}\,R^{\frac{p}{2-p}-N}\,\int_{|x|\ge R} \vo(x)\dx <\infty\,.
\end{equation}
A natural question arises in the view of the above result:
\begin{equation}\label{Q1}\tag{Q-1}
\begin{array}{cc}
    \textit{Is it possible to prescribe an explicit rate for the uniform converge in relative error?}
\end{array}
\end{equation}
 This is the main issue that we want to address in this paper.  For initial data in the class $\X$ in the original variables~\eqref{CPLE} the convergence rate must be at most polynomial (exponential for the problem~\eqref{eq:v}). This phenomenon has been carefully shown in~\cite{Kim} for a different -- yet related --  equation, but the same reasoning applies to~\eqref{CPLE}. It can be easily seen by analysing the case  ``of shifted-in-time Barenblatt''  $u(t,x)=\B_M(t+T,x)$ or ``of shifted-in-space Barenblatt'' $u(t,x)=\B_M(t, x+x_0)$ versus their asymptotic profile. Indeed, the relative error $|\B_M(t+T,x)/\B_M(t,x)-1|$  is of order $ t^{-1}$ (of order $e^{-\frac{\tau}{\beta}}$ in the self-similar variables, i.e.  for~\eqref{eq:v}), while $|\B_M(t,x+x_0)/\B_M(t,x)-1|$ is of order $t^{-\beta}$  (of order $e^{-\tau}$ for~\eqref{eq:v}). In what follows, we give a positive answer to question~\eqref{Q1}.

We base our analysis on the \emph{entropy method}, which we briefly outline below (see also~\cite{Arnold2004}): the entropy is finite when solutions have finite $|y|^{\frac{p}{p-1}}$ moments and this gives rise to the new critical exponent
\begin{equation}\label{p_M}
p_M:=\frac{3(N+1)+\sqrt{(N+1)^2+8}}{2(N+2)}\in(p_c,2).
\end{equation}
Indeed, the $|y|^{\frac{p}{p-1}}$ moments are ``preserved'' along the flow only when $p>p_M$, meaning that $\int_{\rn}|y|^{\frac{p}{p-1}}\vo(y)\dx<\infty$ implies that $\int_{\rn}|y|^{\frac{p}{p-1}}v(\tau,y)\dx<\infty$ for all $\tau>0$.

Our main result provides an \emph{explicit and uniform} rate of convergence for the uniform relative error as long as $\vo\in \X$ and $p_M<p<2$, but in view of the above discussion, for lower values of $p$ it is necessary to impose extra assumptions, which in our case read:
\begin{equation}\label{stronger.GHP.intro}\tag{H}
  \mbox{There exists $D_1>D_2>0$ such that}\quad \V_{D_1}(y) \le \vo(y) \le \V_{D_2}(y)\qquad \forall\,y\in\rn\,.
\end{equation}
Our main result in the whole good fast diffusion range  reads as follows.
\begin{theo} \label{theo:RECR}
Let $N\ge3$, $p_c<p<2$, $0\le \vo \in \LL^1(\rn)\cap\X$,  $M:=\|\vo\|_{\LL^1(\rn)}>0$ and $D=D(M)$ as in~\eqref{Barenblatt-via-vD}. Assume $v$ is a weak solution to~\eqref{eq:v} with initial datum $\vo$ and when $p_c<p\le p_M$ assume moreover that $\vo$ satisfies~\eqref{stronger.GHP.intro}.
Then there exist $\tau_\star=\tau_\star(p,N,M, \|\vo\|_{\X})>0$, $K_\star=K_\star(p,N,M, \|\vo\|_{\X})>0$ and $\sigma=\sigma(p,N)>0$ such that
\begin{equation}\label{convergence.relative.error}
  \left\|\frac{v(\tau,\cdot)-\VD(\cdot)}{\VD(\cdot)}\right\|_{\LL^\infty(\rn)} \le K_\star\,e^{-\sigma\,\tau}\qquad\forall\, \tau\ge \tau_\star\,.
\end{equation}
\end{theo}

\begin{rem} \rm
(i) We notice that in the case $N=2$, Theorem~\ref{theo:RECR} still holds as long as $p\neq \frac{3}{2}$.

\noindent (ii) We notice that in the case $p_M<p<2$, we only need to assume that the initial datum  $\vo$ belongs to $\X$. This is the minimal (and hence the optimal) assumption for uniform convergence in relative error. This was shown in~\cite[Theorem 1.1]{Bonforte2020b}, together with explicit counterexamples by means of initial data in $\LL^1(\rn)\setminus\X$ for which such property simply fails, see~\cite[Proposition 5.1]{Bonforte2020b}.  At the same time, Theorem~\ref{theo:RECR} gives (at least in the range $p_M<p<2$) a uniform convergence rate for the whole class $\X$ which is \emph{independent of the initial datum}.

\noindent (iii) In the case $p_c<p\le p_M$, in order to ensure that the relative entropy is finite we ask for the assumption~\eqref{stronger.GHP.intro}. On one hand, this is essential to ensure that $v(\tau)$ is in the ``stability regime'', i.e. it is close uniformly in relative error to the stationary state: this is where our delicate localized analysis of (possibly sharp) convergence rates takes place. On the other hand, while this assumption may appear very strong, the convergence result~\eqref{UCRE/Xp} shows that a condition similar to \eqref{stronger.GHP.intro} holds true for all solutions with $\vo\in\X$ for ``large times'', namely that $(1-\varepsilon(\tau)) \VD(y) \leq v(\tau, y) \leq (1+\varepsilon(\tau)) \VD(y)$, with $\varepsilon(\tau)\to 0$ as $\tau\to\infty$.  Indeed, in the proofs we only need condition~\eqref{stronger.GHP.intro} to be true for large times, meaning for $\tau\geq \tau_0\geq 0$. We have chosen to set $\tau_0=0$, i.e. we assume~\eqref{stronger.GHP.intro} for small times as well, to focus our attention in understanding the convergence rates and in order to simplify the already technical proofs that contain many parameters.
\end{rem}
As explained earlier, the proof of Theorem~\ref{theo:RECR} uses the entropy method. We now present the main components involved, namely the entropy functional and the relative Fisher information. The {\em relative entropy} is given by:
\begin{flalign}\label{cal-E}
\cE[v(\tau,\cdot)|\VD]&:=\frac{1}{\gamma(\gamma-1)}\int_{\rn}\left\{v^\gamma(\tau,y)-\V_D^\gamma(y)-\gamma \V_D^{\gamma-1}(y)\big[v(\tau,y)-\VD(y)\big]\right\} \dy\,,\quad\mbox{where}\quad \gamma:=1-\frac{2-p}{p-1}\,.
\end{flalign}
The relative entropy production, i.e. minus the time derivative of the entropy along the \eqref{eq:v}~flow, is often called the {\em relative Fisher information} or simply \emph{Fisher information} and has the expression
\begin{multline}\label{cal-I-def}
\cI[v(\tau,\cdot)|\VD]:=\frac{1}{|\gamma-1|^p}\int_{\rn}v(\tau,y)\big(\nabla v^{\gamma-1}(\tau,y)-\nabla  \V_D^{\gamma-1}(y) \big) \cdot\Big(\bb[ v^{\gamma-1}(\tau,y)]-\bb[\V_D^{\gamma-1}(y)]\Big) \dy\,, \\
\quad\mbox{where}\quad\bb[\phi]:=|\nabla\phi|^{p-2}\nabla\phi\,.
\end{multline}
When no confusion arises, we shall write $\cE[v(\tau)]$ and $\cI[v(\tau)]$ instead of $\cE[v(\tau,\cdot)|\VD]$ and $\cI[v(\tau,\cdot)|\VD]$.

The \emph{entropy method} consists in proving that the entropy functional converges to zero exponentially fast. The entropy may be considered a sort of ``nonlinear distance'', adapted to the flow in order to get (possibly sharp) results. It has to be noticed that it also controls from above more standard distances: thanks to the Csisz\'ar--Kullback inequality (Lemma~\ref{CK.inq}), we can easily infer the $\LL^1$-convergence of $v(\tau)$ towards $\VD$.  Then, by a Gagliardo type interpolation and uniform regularity estimates, it is possible to extend the convergence of $v(\tau)$ towards $\VD$ in different $L^q$ norms up to $L^\infty$, or even up to some $C^\alpha$ norms, with (almost) the same convergence rates.

It is important to stress that $\cE$ or $\cI$ \textit{need not to be finite}, even for smooth bounded integrable solutions. In the range of parameters where $\cE[v(\tau)]$ and $\cI[v(\tau)]$ may be a priori unbounded, it is a delicate task to show that $\cE[v_0]$ and $\cI[v_0]$ bounded implies $\cE[v(\tau)]$ and $\cI[v(\tau)]$ \emph{bounded for all positive times (or at least integrable in time)}. This remark applies to the whole range of parameters $p\in (1,2)$.

In order to prove exponential convergence decay of the entropy functional, the first step consists in proving that, for a~solution $v$ to problem~\eqref{eq:v}, it holds (in a suitably weak sense)
\begin{equation} \label{dE=-I}
\frac{\d}{\d\tau} \cE[v(\tau,\cdot)|\VD]  = - \cI[v(\tau,\cdot)|\VD]\,.
\end{equation}
We shall clarify all these details in Subsection~\ref{ssub:entropy-fisher}. Having this relation, it remains to prove that the Fisher information controls the entropy functional, at least along the flow. This is a key step in the entropy method: we need the so-called \emph{entropy -- entropy production inequality}, namely that for a positive constant $c>0$ the following inequality holds
\begin{equation}\label{entropy-entropy-production-inq.intro}
\cI[v(\tau)|\VD]\ge c\, \cE[v(\tau)|\VD]\qquad\mbox{for all $\tau>0 $}\,.
\end{equation}
This represents one of the most delicate aspects of the problem, since the emergence of new critical exponents necessitates different strategies across distinct parameter ranges. Thanks to the entropy -- entropy production inequality, we deduce (possibly in a weak sense) a differential inequality, that implies the exponential decay of the entropy towards zero:
\[
\frac{\d}{\d\tau} \cE[v(\tau)|\VD] \leq - c\, \cE[v(\tau)|\VD]\qquad\mbox{hence}\qquad \cE[v(\tau)|\VD]\leq e^{-c\tau} \cE[\vo|\VD]\,.
\]
A (nontrivial) weak version of Gronwall's Lemma is needed in some cases.  Once exponential decay of the entropy is established, we can transfer the rates of convergence to $\LL^q$ or $C^\alpha$ distances, as discussed above.

Establishing inequality~\eqref{entropy-entropy-production-inq.intro} is a major difficulty of this method. In the study of~\eqref{entropy-entropy-production-inq.intro} a new exponent naturally enters into this panorama, namely
\begin{equation}\label{pD}
  p_D:=\frac{2N+1}{N+1}\in\left(p_2, 2\right)\,.
\end{equation}
In the range $p_D\leq p < 2$, inequality~\eqref{entropy-entropy-production-inq.intro} holds not only for solutions to~\eqref{eq:v} but for any function for which the relative entropy and relative Fisher information are finite. It is, indeed, an equivalent form of a Gagliardo--Nirenberg--Sobolev inequality, already considered in~\cite{Cordero2004,Agueh2008}, see also Section~\ref{ssec:conv-above-p_D}.
We stress that when $p_c\le p < p_D$ such a clean result is missing. Indeed, $\cE$ and $\cI$ are well-defined only when the solution is sufficiently close to the Barenblatt profile $\VD$, in the spirit of assumption~\eqref{stronger.GHP.intro}. In this range, in order to establish~\eqref{entropy-entropy-production-inq.intro} we compare $\cE$ and $\cI$ with linearized quantities around the steady state. Inequality~\eqref{entropy-entropy-production-inq.intro} finally follows from a weighted Hardy--Poincar\'e inequality (Proposition~\ref{prop:hp} and~\cite{Chlebicka2022}) together with a delicate comparison between nonlinear and linearized inequalities (Lemmata~\ref{lem:I-est},~\ref{lem:ABC} and~\ref{lem:linE-leq-linI-eps}).

Our next result provides sufficient conditions for explicit exponential decay rates, optimal in some cases, of the entropy functional in the whole good fast diffusion range.
\begin{theo}\label{theo:Entropy-decay} Under the same assumptions of Theorem~\ref{theo:RECR}, there exists an explict $\lambda>0$ such that
\begin{equation}\label{entropy-decay}
\cE[v(\tau,\cdot)|\VD] \leq C\,  e^{-\lambda \tau}\quad \forall \tau>0\,,
\end{equation}
and the constant $C$ depends on $p$, $N$, $\vo$ and, when $p_c<p\leq p_M$ also on $D_1, D_2$ of~\eqref{stronger.GHP.intro}. Assume moreover that $0\le \vo\in C^2(\R^N) $ is radially symmetric and decreasing ($\partial_r \vo\leq 0$), it satisfies~\eqref{stronger.GHP}, and that one of the following sets of assumptions holds:
\begin{enumerate}[{\it (i)}]
\item $p_M<p<2$ and there exist $A>0$ and $R_0>0$ satisfying
\begin{equation*}
  \partial_r \vo(r)\leq 0\quad\mbox{and}\quad |\partial_r \vo(r)|\leq A\, r^{-\frac{2}{2-p}}\ \ \forall r\ge R_0\,,
\end{equation*}
\item $p_c<p\leq p_M$ and there exist $D_1>D_2>0$  such that
\begin{equation}\label{H-derivatives}
    \partial_r \V_{D_2}(r) \leq \partial_r \vo(r) \leq \partial_r \V_{D_1}(r)\quad \forall\,\, r\ge0\,.
\end{equation}
\end{enumerate}
then inequality~\eqref{entropy-decay} holds true for $\Lambda=\Lambda^\star$ where $\Lambda^\star$ is given by
\begin{equation*}
\Lambda^\star= \begin{cases}
p-N(2-p)=\frac{1}{\beta} &\quad\mbox{when}\quad p_M<p<2\,,\\
\frac{\left[p-N(2-p)(p-1)\right]^2}{4p(p-1)(2-p)} &\quad\mbox{when}\quad p_c<p\leq p_M\,.
\end{cases}
\end{equation*}
\end{theo}
The proof is inspired by the methods used in~\cite{Agueh2009} and~\cite{Blanchet2009}, where the authors opened the way for establishing the needed relations between the relative entropy and the Fisher information, together with their linearized counterparts. We also take into account recent improvement established in~\cite{Bonforte2020b} in order to simplify some assumptions and proofs.\newline

\noindent\textit{Optimality and related results. }The above theorem provides the optimal rates of convergence for radially decreasing functions. The rate $\Lambda^\star$ has a precise meaning explained below. For non-radial functions, we always provide explicit and computable rates (as in Theorem \ref{theo:RECR}), whose optimality is not known.

To the best of our knowledge,  only a few  closely related results are present in literature, and none of them provides a proof of convergence rates in the entire range $p_c<p<2$, whether sharp or not. The first pioneering result, valid in the range $p_D<p<2$, is due to Del Pino and Dolbeault, announced in~\cite{DelPino2002} and proven in detail in~\cite{DelPino2003}. Their method connects convergence rates to equilibrium with the optimal constant in suitable Gagliardo--Nirenberg inequalities, obtained in~\cite{Del_Pino_2003}, using a new entropy method, that inspired the one developed in the present paper.   On the other hand, Agueh in \cite{Agueh2003,Agueh2008} provides explicit convergence rates in the same range, and the proof relies on entropy methods based on the Wasserstein gradient flow formulation of the problem. It uses in an essential way displacement convexity of the entropy, so the range of the validity of the method is limited to $p_D<p<2$. In the whole range $p_c<p<2$, the paper \cite{Agueh2009} provides convergence rates under similar assumptions. Unfortunately, there are some gaps in the proofs that affect the validity of the results in the range $p_c<p<p_M$: some integrations by parts are not correct, namely \cite[inequality (4.10)]{Agueh2009} and some quantities that are treated as bounded there, may be unbounded. For instance, the Fisher information happens to be only $\mathrm{L}^1$ in time (not necessarily bounded), hence the Gronwall-type argument of \cite{Agueh2009} fails and it needs to be done differently, e.g. as we do in the proof of Proposition~\ref{convergence.lp.below.pm} using a weak version of Grownwall's Lemma, see Lemma~\ref{weak.gronwall}.

\textit{Let us now explain in which sense the rate $\Lambda^\star$ is optimal. }The optimal decay rate of the entropy functional is given by the optimal constant in inequality~\eqref{entropy-entropy-production-inq.intro} or equivalently to the infimum of a Rayleigh's type quotient
\[
  \mathcal{C}=\inf \frac{\cI[v|\VD]}{\cE[v|\VD]}\,,
\]
where the infimum is taken among all non-negative and regular enough functions $v$ such that $\int_{\mathbb{R}^N}v(y)\dy=\int_{\mathbb{R}^N}\VD(y)\dy$. As we previously explained, the value of $\mathcal{C}$ in the range $p_D\leq p <2$ is known, as it is the optimal constant in a suitable Gagliardo--Nirenberg--Sobolev-type inequality, whose value was obtained for the first time in~\cite{Del_Pino_2003}, see also~\cite{Cordero2004} for an alternative proof.  However, when $p_c<p<p_D$, this correspondence is not available, and it remains unclear whether $\mathcal{C}\ge 0$ is finite and how to determine its value under general conditions. By taking inspiration from~\cite{Blanchet2009, Bonforte2010}, it is possible to compute this constant, at least for radial functions that are close to the Barenblatt profile $\VD$. Namely, $\Lambda^\star$ can be obtained as
\begin{equation}\label{optimal.rate}
  \Lambda^\star=\liminf_{\varepsilon\rightarrow0}\inf_{v\in\mathcal{S}_\varepsilon}\frac{\cI[v|\VD]}{\cE[v|\VD]}\,,
\end{equation}
where $\mathcal{S}_\varepsilon$ is the set of radial, smooth functions which satisfy~\eqref{stronger.GHP.intro}, inequality~\eqref{H-derivatives}, and the two conditions:
\[
  \left\|\frac{v-\VD}{\VD}\right\|_{\mathrm{L}^\infty(\rn)}<\varepsilon\,\quad \mbox{and}\quad \left\|\frac{\partial_rv-\partial_r\VD}{\partial_r\VD}\right\|_{\mathrm{L}^\infty(\rn)}<\varepsilon\,.
\]
Under suitable assumptions, radially decreasing solutions to~\eqref{eq:v} will eventually enter the set $\mathcal{S}_{\varepsilon}$ for any $\varepsilon\in\left(0,1\right)$, as we shall see below.  Therefore, it is perfectly natural to consider the limit~\eqref{optimal.rate} in such sets. As a consequence, the rate $\Lambda^\star$ is sharp in the following sense: for any $\lambda>\Lambda^\star$, one can always find an initial datum $\vo\in\mathcal{S}_\varepsilon$ for which the estimate $\cE[v(\tau)|\VD]\leq \cE[\vo|\VD]e^{-\lambda \tau}$ cannot hold for sufficiently large $\tau>0$.

In order to compute the value $\Lambda^\star$, we shall perform a delicate linearisation around the Barenblatt profile $\VD$. Last but not least, the value of $\Lambda^\star$ can be obtained (through a change of variables, see Section~\ref{ssec:proof_of_theorem_2}) as the optimal constant $\Lambda_{\textrm{opt}}$ in the following Hardy--Poincaré inequality, which we state here for radial functions for the sake of simplicity:
\begin{equation}\label{optimal-HP.intro}
\Lambda_{\textrm{opt}} \int_0^{\infty}   \frac{|\phi(r)|^2}{(1+r^2)^\frac{2-\gamma}{1-\gamma}} \d\mu(r)
\leq \int_0^{\infty}  \frac{|\phi'(r)|^2}{(1+r^2)^\frac{1}{1-\gamma}} \d\mu(r)\,,\qquad\mbox{where}\qquad\d\mu(r)=r^{\frac{2N(p-1)}{p}-1}\dr\,.
\end{equation}
This inequality holds under appropriate orthogonality conditions, see~\cite{Bonforte2017a,Bonforte2017b,Bonforte_2023}.  Being the same rate as the linearisation around the equilibrium, we do not see how this rate could be improved without imposing further conditions. \newline

Before proceeding with further results, we believe it is worthwhile to summarise what has been achieved so far in the good fast diffusion range $\frac{2N}{N+1}=p_c<p<2$ and clarify under what regime the entropy and the Fisher information are well-defined. Recall that the conservation of mass, i.e.~\eqref{p_c}, is true in this range and the convergence towards the Barenblatt profile $\VD$ holds as in~\eqref{mild-convergence-R}, and the parameter $D$ is uniquely determined by the mass of the initial datum $\vo$ as in~\eqref{Barenblatt-via-vD}. The entropy functional $\cE[v|\VD]$, defined in~\eqref{cal-E}, is finite along the flow of~\eqref{eq:v} under the sole assumption $\vo\in\X$ for $p_M<p<2$, see Lemma~\ref{lem:ent-lin-ent}. It is unbounded in general for $p_c<p\leq p_M$, but bounded in the range under extra assumption~\eqref{stronger.GHP}, see Lemma~\ref{lem:ent-lin-ent}. The Fisher information $\cI[v|\VD]$, defined in~\eqref{cal-I-def}, is $\mathrm{L}^\infty (\tau_0, \infty)\cap \mathrm{L}^1(\tau_0, \infty)$ for any $\tau_0>0$ along the flow of~\eqref{eq:v} in the range $p_M<p<2$ under the sole assumption $\vo\in\X$, see Lemma~\ref{lem:ent/ent-prod}. In the same range identity~\eqref{dE=-I} holds pointwise. In the range $p_c<p\leq p_M$ the Fisher information is only in $\mathrm{L}^1(\tau_0, \infty)$ for any $\tau_0>0$ along the flow~\eqref{eq:v} and identity~\eqref{dE=-I} holds in a weaker sense~\eqref{weak.entropy.production}, see Lemma~\ref{lem:ent/ent-prod-2}. The relation between the entropy and the Fisher information, i.e. inequality~\eqref{entropy-entropy-production-inq.intro}, is obtained via a Gagliardo--Nirenberg inequality when $p_D\leq p <2$, see~\cite{Agueh2008,Cordero2004}, and by a comparison with a Hardy--Poincar\'{e} inequality, see Proposition~\ref{prop:hp} and~\cite{Chlebicka2022}.  We summarise all the above information in the table below.

\begin{table}[h!]
\centering
\resizebox{\textwidth}{!}{
\begin{tabular}{|c|c|c|c|c|c|c|}
\hline
\textbf{Range} & $\cE[v(\tau)|\VD]$  & $\cI[v(v)| \VD]$ & $\tfrac{d}{d\tau}\cE[v(\tau)]=-\cI[v(v)]$ & \makecell{Inequality~\eqref{entropy-entropy-production-inq.intro} \\  $\cI[v(v)| \VD] \ge c\,\cE[v(\tau)|\VD] $} & \makecell{Entropy Decay\\ Rates} & \makecell{UCRE \\ Rates}  \\
\hline
$p_D\leq p<2$ & \makecell{$\mathrm{L}^\infty(\tau_0, \infty)$ \\ if $\vo\in\X$ \\ Lemma~\ref{lem:ent-lin-ent}} & \makecell{$\mathrm{L}^\infty\cap\mathrm{L}^1(\tau_0, \infty)$ \\ if $\vo\in\X$ \\ Lemma~\ref{lem:ent/ent-prod}} & Yes & \makecell{GNS inequality \\ see~\cite{Agueh2008,Cordero2004}} & \makecell{Theorem~\ref{entropy-decay} if $\vo\in\X$ \\ Optimal in $\X$: unknown \\ Optimal RD: $\Lambda^\star$}  & \makecell{Theorem~\ref{UCRE/Xp} if $\vo\in\X$ \\ Optimal in $\X$: unknown \\ Optimal expected: $\mathcal{O}(e^{-\tau})$}  \\
\hline
$p_M\leq p< p_D$ & \makecell{$\mathrm{L}^\infty(\tau_0, \infty)$ \\ if $\vo\in\X$ \\ Lemma~\ref{lem:ent-lin-ent}}  & \makecell{$\mathrm{L}^\infty\cap\mathrm{L}^1(\tau_0, \infty)$ \\ if $\vo\in\X$ \\ Lemma~\ref{lem:ent/ent-prod}} & Yes & \makecell{Hardy-Poincaré inequality \\ Prop.~\ref{prop:hp} and~\cite{Chlebicka2022}} & \makecell{Theorem~\ref{entropy-decay} if $\vo\in\X$ \\ Optimal in $\X$: unknown \\ Optimal RD: $\Lambda^\star$} & \makecell{Theorem~\ref{UCRE/Xp} if $\vo\in\X$ \\ Optimal in $\X$: unknown \\ Optimal expected: $\mathcal{O}(e^{-\frac{\tau}{\beta}})$}  \\
\hline
$p_c< p\leq p_M$ & \makecell{$\mathrm{L}^\infty(\tau_0, \infty)$ \\ if~\eqref{stronger.GHP.intro} \\ Lemma~\ref{lem:ent-lin-ent}}  & \makecell{$\mathrm{L}^1(\tau_0, \infty)$ \\ if~\eqref{stronger.GHP}  \\ Lemma~\ref{lem:ent/ent-prod-2}}  & \makecell{As in~\eqref{weak.entropy.production} \\ Lemma~\ref{lem:ent/ent-prod-2}} & \makecell{Hardy-Poincaré inequality \\ Prop.~\ref{prop:hp} and~\cite{Chlebicka2022}} & \makecell{Theorem~\ref{entropy-decay} if~\eqref{stronger.GHP.intro} \\ Optimal if~\eqref{stronger.GHP}: unknown \\ Optimal RD: $\Lambda^\star$} & \makecell{Theorem~\ref{UCRE/Xp} if $\vo\in\X$ \\ Optimal in $\X$ or~\eqref{stronger.GHP.intro}: unknown \\ Optimal expected: $\mathcal{O}(e^{-\frac{\tau}{\beta}})$} \\
\hline
\end{tabular}
}
\caption{Summary of the main properties of solutions to~\eqref{eq:v} in the good fast diffusion range. Recall that $\beta=(p-N(2-p))^{-1}$.}
\label{tab:gfdr}
\end{table}
Lastly, let us comment on the final two columns of Table~\ref{tab:gfdr}. As previously explained, Theorem~\ref{entropy-decay} provides the explicit and computable rate for the decay of the entropy functional under the assumption $\vo\in\X$ in the range $p_M < p < 2$, and, for $p_c < p \leq p_M$, under the stronger assumption~\eqref{stronger.GHP.intro}. In the range $p_D \leq p < 2$, the optimal decay rates for solutions with finite $|y|^{\frac{p}{p-1}}$ moment have been computed in~\cite{Agueh2008}. It remains unclear whether improved rates should be expected when the initial datum belongs to $\X$, as it has been proven for similar equations to~\eqref{eq:v} in~\cite{Bonforte2020a}. For $p_c < p < p_D$, even identifying the optimal decay rate within the class $\X$ appears to be currently out of reach. Nevertheless, for radial and decreasing solutions satisfying~\eqref{H-derivatives}, the optimality is achieved in the sense discussed above.

Regarding decay rates for the convergence in the relative error (UCRE, the last column of Table~\ref{tab:gfdr}), constructive and computable rates are established under the assumption $\vo\in\X$ for $p_M < p < 2$, and under assumption~\eqref{stronger.GHP.intro} for $p_c < p \leq p_M$. It remains an open problem to derive constructive rates for $p_c < p \leq p_M$ relying solely on the assumption $\vo\in\X$. Optimal rates are also unknown, but they may be conjectured based on specific examples, as discussed below question~\eqref{Q1}.\newline

\noindent\textit{Convergence of radial derivatives}. The equation \eqref{CPLE} is a gradient-driven diffusion problem, therefore, it is quite natural to investigate whether or not the convergence in relative error may hold also for derivatives of solutions. This has been established in~\cite[Theorem 1.4]{Bonforte2020b} in the case of radial derivatives for a radially decreasing initial datum and under some assumptions on the spatial decay of the radial derivative. The result reads as
\begin{equation}\label{gradient.decay.convergence}
\Big\|\frac{\partial_r v(\tau,\cdot)-\partial_r \VD(\tau, \cdot)}{\partial_r \VD(\tau, \cdot)}\Big\|_{L^\infty(\RN)}\xrightarrow[{\tau\to +\infty}]{} 0 \,.
\end{equation}
Interestingly, at least in dimension $N=1$, having a symmetric and decreasing initial datum is a necessary condition for the convergence of the spatial derivative. Counterexamples are provided in~\cite[Remark 1.5]{Bonforte2020b}. Even if counterexamples were not constructed in higher dimensions, we believe that for $N\ge2$, having a radially symmetric and decreasing initial datum is also a necessary condition.  We are interested in the following refinement of \eqref{gradient.decay.convergence}:
\begin{equation}\label{Q2}\tag{Q-2}
\begin{array}{c}
    \textit{Is it possible to prescribe an explicit rate for the uniform converge in relative error of radial derivatives?}
\end{array}
\end{equation}
As for~\eqref{Q1}, the rate should be polynomial for the original problem~\eqref{CPLE} and exponential in rescaled variables~\eqref{eq:v}. Again this can be inferred simply by considering a time-delayed Barenblatt profile: $u(t,x)=\B_M(t+T,x)$, in perfect analogy to what we observed while attempting question~\eqref{Q1}. We provide the first answer in the {\it good range} as follows.
\begin{theo}[\textbf{Convergence in relative error for radial derivatives when $p\in (p_c,2)$}]\label{Thm:gradient.convergence}
Let $N\ge 2$, $p_c<p<2$, and let $0\le \vo\in \LL^1(\rn) \cap C^2(\R^N) $ be radially symmetric with $M:=\|\vo\|_{L^1(\RN)}>0$, $D=D(M)$ as in~\eqref{Barenblatt-via-vD}. Let $v$ be a weak solution to~\eqref{eq:v} with datum $\vo$, and assume that one of the following set of assumptions holds:
\begin{enumerate}[{\it (i)}]
\item $p_M<p<2$ and there exist $A>0$ and $R_0>0$ satisfying
\begin{equation}\label{assumptions.thm.derivative}
  \partial_r \vo(r)\leq 0\quad\mbox{and}\quad |\partial_r \vo(r)|\leq A\, r^{-\frac{2}{2-p}}\ \ \forall r\ge R_0\,,
\end{equation}
\item $p_c<p\leq p_M$ and there exist $D_1>D_2>0$  such that
\begin{equation}\label{derivatives.supercritical}
    \partial_r \V_{D_2}(r) \leq \partial_r \vo(r) \leq \partial_r \V_{D_1}(r)\quad \forall\,\, r\ge0\,.
\end{equation}
\end{enumerate}
Then there exist $\tau_\star>0$, $k_\star>0$ and $\lambda=\lambda(p,N)>0$ such that
\begin{equation*}
\Big\|\frac{\partial_r v(\tau,\cdot)-\partial_r \VD(\tau, \cdot)}{\partial_r \VD(\tau, \cdot)}\Big\|_{L^\infty(\RN)}\le k_\star\,e^{-\lambda\,\tau}\,\quad\forall\,\tau\ge \tau_\star\,,
\end{equation*}
where $\partial_r v$ (resp. $\partial_r \V_D$) is the radial derivative of $v$ (resp. $\VD$). When $p_M<p<2$ then $ \tau_\star=\tau_\star(\vo, A, R_0, M, p, N)$  and $k_\star=k_\star(\vo, A, R_0, M, p, N)$, while when $p_c<p\leq p_M$ then $ \tau_\star=\tau_\star(\vo, D_1,D_2, T, p, N)$  and $k_\star=k_\star(\vo,D_1,D_2, T,  p, N)$.
\end{theo}
\begin{rem}
    \rm We notice that the regularity assumption on the initial datum, i.e., $\vo\in C^2(\rn)$ is imposed solely for the sake of simplicity in the exposition.
    More precisely, the regularity threshold for our proof to work is $\vo\in C^{1, \alpha}(\rn)$ for a suitable $\alpha>0$ depending on $p$.
\end{rem}

The proof of the above theorem does not rely on the techniques employed for Theorem~\ref{theo:RECR}. Instead, it relies on a clever idea introduced in \cite{Iagar2008}, that puts in 1 to 1 correspondence radial solutions to equation \eqref{CPLE} with radial solutions to a fast diffusion type equation (a density-driven diffusion), as thoroughly discussed in Section~\ref{sec:conv-radial}.

\subsection{The very fast diffusion range: \texorpdfstring{$1<p\leq p_c$}{1<p≤p_c}}

The long-time behaviour of solutions to~\eqref{CPLE} in the very fast diffusion range $1<p\leq p_c$ remains poorly understood. In fact, this regime presents several interesting open questions, many of which are still unanswered after more than 50 years of deep study. While solutions exist when $0\leq \uo \in \LL^1_{\textrm{loc}}(\rn)$ for every $p>1$, see~\cite{DiBenedetto1990}, several key properties of them that hold in the good range, are lost in the very fast diffusion regime, see Section~\ref{sec:prelim}. Let us begin with the failure of the mass conservation. As already mentioned, the value $p=p_c$ is critical for the conservation of mass, i.e.~\eqref{p_c}, of  solutions that are in $\LL^1(\rn)$. It is known that in this case solutions whose initial datum is in $\LL^1(\rn)$ still conserve mass, see for instance the recent survey~\cite{Vazquez:2023aa}. However, when $1<p<p_c$ the mass is not conserved anymore. More surprisingly, for a large class of data (i.e., $\uo\in L^{r}(\rn)$ with $r=n(2-p)/p$ and/or $\uo(x)\leq c \VD(x)$ for some $D>0$) the corresponding solution extinguishes at  $T>0$ (i.e., $u(t,\cdot)=0$, for all $t>T$ a.e. in $\rn$), cf.~\cite[Chapter 11]{Vazquez2007}.

In this range, the availability of the fundamental solution is lost: the Barenblatt profile~\eqref{beta} is no longer in $\LL^1(\rn)$ and does not represent the self-similar solution that takes a Dirac $\delta_0$ as its initial datum. However,  a \emph{pseudo}-Barenblatt solution can still be defined in this range as a self-similar profile by the formula
\begin{equation}\label{pseudo-Barenblatt}
\B_{D,T}(t, x)=R^{-N}_T(t)\VD\left(\frac{x}{R_T(t)}\right)\,,
\end{equation}
where $\VD$ is as in~\eqref{vD} and $R_T$ is given by
\begin{equation}\label{time.rescaling-pc}
  R_T(t):=\left(\frac{T-t}{|\beta|}\right)_+^\beta\quad\mbox{for}\ \ 1<p<p_c\qquad\text{ and }\qquad R_T(t):=\exp{\left\{ \ell(T+t)\right\} }\quad\mbox{if}\ \ p=p_c\,,
\end{equation}
where $\beta$ is as in~\eqref{beta}, while $\ell>0$ is a free parameter. We point out that $\beta<0$ for $1<p<p_c$. Notice also that, since $\B_{D,T}(t, \cdot)\notin \LL^1(\rn)$ for any $0<t<T$, the parameter $D$ in $\VD$ is a ``free'' parameter which does not represent the mass anymore. Furthermore, for $1<p<p_c$ \emph{pseudo}-Barenblatt solutions are positive until they vanish at a finite time $T>0$, but for $p=p_c$ they are positive for all times and $T\ge 0$ is a free parameter.

It is unclear what is the role of pseudo-Barenblatt profiles in the very fast diffusion range, more specifically if they are still attractors of solutions to~\eqref{CPLE}, i.e., if they represent the extinction profile of solution in the original variables and in what topology. More precisely, solutions corresponding to integrable non-negative initial data can converge to a pseudo-Barenblatt profiles, but there are also several other potential attractors known. Let us consider the case of the exponent $p_Y$ defined as
\begin{equation}\label{pY}
  p_Y:=\frac{2N}{N+2}\in\left(1, p_c\right)\,.
\end{equation}
For $p=p_Y$, in~\cite{Iagar2008} the authors find an explicit solution obtained by the separation of variables (see~\cite[(8.11)]{Iagar2008}) that can be written as
\begin{equation*}
  U(t, x)=(T-t)^\frac{1}{2-p}F(x)\quad\mbox{where}\quad F(x)=F(|x|)\sim |x|^{-\frac{N^2}{N-1}}\quad\mbox{as}\quad |x|\to\infty\,.
\end{equation*}
In this case, by using the results of~\cite{Iagar2008} and the more geometric setting of~\cite{delPinoSaez2001}, it can be proven that $U$ is an attractor for a certain class of solutions when $p=p_Y$. However, a broad picture is far more complex: due to~\cite[Theorem 8.2]{Iagar2008} it is known that for any $p\in\left(1, p_c\right)$ a family of self-similar solutions with different tails from $U$ and $\B_{D,T}$ exist. To the best of our knowledge the following natural question has not been addressed before:
\begin{equation}\label{Q3}\tag{Q-3}
\begin{array}{c}
    \textit{What is the basin of attraction of the pseudo-Barenblatt solution in the very fast diffusion regime? }
\end{array}
\end{equation}
Since solutions may extinguish in finite time and the extinction time depends on the initial datum in an implicit way, convergence toward the pseudo-Barenblatt solutions is subtle and needs to be treated with care in this very fast diffusion regime. Let us explain why. Suppose that a given solution $u$ converges to the pseudo-Barenblatt solution $\B_{D, T}$ in $\mathrm{L}^\infty(\mathbb{R}^n)$. As first preliminary step, in order to meaningfully measure the convergence of $u$ to $\B_{D, T} $, it is necessary to assume that both $u$ and $\B_{D, T}$ extinguish at the same time $T$. Indeed, if $u$ vanishes at some later time $T'>T$, since $\B_{D, T}(t)=0$ for $t>T$, it follows that for $t \in (T, T') $ the difference $ u(t,x) - \B_{D, T}(t, x)$ remains strictly positive in $\LL^\infty$ norm, reflecting only the discrepancy in extinction times, not the convergence rate of $u$ to $\B_{D, T}$ as $t \to T^{-}$. Unfortunately, there are no clear methods to compute the extinction time of a given solution in terms of the initial datum, only upper and lower bounds are known, see~\cite{Bonforte_2010,Vazquez2007}. The only exception concerns the case $p=1$ and one spatial dimension: then there is an explicit formula for the extinction time for the \emph{Total Variation flow}, see \cite[Proposition 2.10]{Bonforte2012}. In the known literature for similar equations as~\eqref{CPLE}, a way to bypass this problem is to consider an initial datum $\uo$ appropriately close to the stationary state. In the same spirit of~\cite{Blanchet2009}, we require the initial datum to satisfy the following inequality:
\begin{equation}\label{strong.GHP.veryfast}
  \B_{D_1, T}(x)\leq \uo(x) \leq \B_{D_2, T}(x)\qquad\mbox{for all}\quad x\in\mathbb{R}^N\,.
\end{equation}
We also stress the fact that both the above upper and lower bounds are necessary in order to be able ensure that $u(t)$ extinguish exactly at $T>0$, as well as its corresponding asymptotic pseudo-Barenblatt $\B_{D, T}(t)$.  Lastly, the above inequality should hold in the whole $\rn$, not only in some large set; counter-examples are known, see~\cite{Daskalopoulos2008}. In the case $p=p_c$, assumption~\eqref{strong.GHP.veryfast} is also convenient, since it allows us to address~\eqref{Q3} by using the entropy method. In order to do so, we introduce a natural rescaling, similar to that used in~\eqref{eq:v} which transforms this delicate problem into the study of the convergence to stationary solutions of a nonlinear Fokker--Planck type equation as~\eqref{eq:v}. Indeed, let us assume that the initial datum $\uo$ satisfies~\eqref{strong.GHP.veryfast} and consider the function
\begin{equation}\label{change.variables}
  v(\tau, y)={R_T^N(t)}u(t,x)
\end{equation}
with $T$ as in~\eqref{strong.GHP.veryfast} and the couple $(\tau, y)$ as in~\eqref{change:rho-u}. Then the function $v$ satisfies the problem~\eqref{eq:v} with an initial datum $\vo$ which is obtained accordingly to~\eqref{change.variables}. Lastly, we remark that condition~\eqref{strong.GHP.veryfast} in the new variables is nothing else but the assumption~\eqref{stronger.GHP.intro} on the initial datum $\vo$. \newline

\noindent\textit{New critical exponents in higher dimensions and the limit $p\to1$}. Another threshold appears when one considers integrability properties of the difference of two Barenblatt solutions  $\B_{D_2, T}-\B_{D_1, T}$,   defined in~\eqref{pseudo-Barenblatt}, for the same $T>0$, but for different $D_1$ and $D_2$. As it is clear from  Lemma~\ref{decay.lemma}, the difference is integrable when
\begin{equation}\label{integrability.condition}
    1<N<\frac{p}{(2-p)(p-1)}\,.
\end{equation}
  Note that condition~\eqref{integrability.condition} is satisfied in low dimensions ($1<N<6$) for all $p\in(1,2)$, while in high dimensions ($N\geq 6$) for $p$ close to $1$ and close to $2$, precisely for $p\in (1,p_1)\cup(p_2,2)$ where
 \begin{equation}
     \label{p1p2}
 p_1:=\tfrac32-\tfrac1{2N}-\tfrac{\sqrt{N^2-6N+1}}{2N}\qquad\text{and}\qquad p_2:=\tfrac32-\tfrac1{2N}+\tfrac{\sqrt{N^2-6N+1}}{2N}\,.
 \end{equation}
One of the quantities we need to control in the entropy method is the difference of masses of solutions often referred to as {\it relative mass}, in which these exponents play a key role. Even if for $p\in (1,p_1)\cup(p_2,2)$ two solutions are not integrable, when the difference is indeed integrable, such relative mass is preserved along the flow (under natural assumptions). This is crucial in the selection of the asymptotic profile, i.e. the right parameter $D$. We explain it in detail in Remark~\ref{rem.hyp.intro}~(ii).

Assumption~\eqref{integrability.condition} allows to address \eqref{Q3} using the entropy method via a weighted linearization, and the weighted Hardy--Poincar\'e inequalities, following the strategy of~\cite{Blanchet2009,Bonforte2010a}, but facing a number of extra difficulties that naturally arise due to the geometrical reasons in our complex panorama. In high dimensions the picture becomes more involved: two new critical exponents appear and a new surprising behaviour is found for $p\in (1,p_1)$. Indeed, this behaviour is unexpected since it is dramatically different from the limiting case $p=1$, often called the Total Variation Flow, see \cite{MazonBook2004}. In this ``extreme'' case the fine properties of solutions drastically differ from $p>1$. For instance bounded solution are not even continuous if the initial datum is not, but new discontinuities cannot be created, see \cite{Bonforte2012}. In the limit case the asymptotic behaviour is fully understood only in one spatial dimension, where convergence to the characteristic function of an interval is proven (in all $L^q$ norms and in the relative error). However, the convergence rates to stationary solutions are simply not possible to find in this case. In \cite{Bonforte2012}, via the explicit constructions, it is proven that for every possible rate, there exists a solution that converges with a faster rate, and also another one which converges more slowly. In higher dimension it is known that stationary solutions are characteristic functions of Cheeger sets satisfying certain geometric conditions. The convergence (with or without rates) to the stationary profiles is an intriguing  problem when $p=1$ and $N\ge 2$ that, to the best of our knowledge, remains open. \newline

\noindent\textit{Our main results in the very fast diffusive range.} Before stating our results, we recall that, when $1<p<p_c$ the natural rescaling \eqref{change.variables}  maps $(0,T)$ into $(0,\infty)$ and makes the pseudo-Barenblatt solutions stationary. In this way the parameter $T$ is fixed in the change of variables, and we transform this delicate problem into the analysis of the convergence to equilibrium of solutions to the nonlinear Fokker--Planck equation analogous to~\eqref{eq:v}. Namely we study now, as in the good fast diffusion range, the convergence as $\tau\to \infty$ to the pseudo-Barenblatt profiles $\V_{D}$, which now do not depend on $T$. In the critical case $p=p_c$, the situation does not change much. In~\eqref{change.variables} and~\eqref{change:rho-u} we shall use $R_T$ defined in~\eqref{time.rescaling-pc} instead of the one defined in~\eqref{R(t)}. At this moment, we shall not specify the value of $\ell$ in~\eqref{time.rescaling-pc}. This will be done when needed.

The main results in this range are two. The first one is a {\it general theorem} that contains a set of hypotheses under which convergence in relative error to a pseudo-Barenblatt solution is guaranteed. We shall explain thoroughly the role of each assumption, and their validity in different settings. The second result refines the first result by showing exponential convergence rates, essentially under the same assumptions.

\begin{theo}\label{theo:meta1}
Suppose $N\geq 2$ and and $p\in(1,2)\setminus\{\frac{3}{2}\}$ satisfy condition~\eqref{integrability.condition}. Let $v$ be a non-negative weak solution to~\eqref{eq:v} with non-negative initial datum $\vo\in \LL^1_{ \rm loc}(\rn)\cap\LL^\infty(\rn)$. Moreover, we assume what follows.
\begin{enumerate}[(i)]\item   Suppose that there exist $D_1>D_2>0$, such that
\begin{equation}\label{stronger.GHP.fast.range}
\V_{D_1}(y) \le \vo(y) \le \V_{D_2}(y)\quad\forall y\in\rn\,;
\end{equation}
\item Let $D>0$ be such that for every $\tau>0$ it holds
\[
\int_{\rn}(v(\tau,y)-\V_{D}(y))\dy=\int_{\rn} (\vo (y)-\V_{D}(y))\dy=0\,;
\]
\item Let $D>0$ be as in {\it (ii)}. Assume there is $\tau_0>0$, such that   $ \tau\mapsto\cI[v(\tau\,,\cdot)|\V_{D}]\in \LL^1(\tau_0, \infty)$  and
\[
\cE[v(\tau_1)|\V_{D}] - \cE[v(\tau_2)|\V_{D}] = \int_{\tau_1}^{\tau_2}  \cI[v(s)|\V_{D}] \ds \quad\mbox{for all}\quad \tau_2\geq\tau_1\geq\tau_0>0\,;
\]
\item Let $\tau_0$ be as in {\it (iii)}. Assume there is $\alpha\in(0,1)$ such that $v$ is $C^{1,\alpha}$-regular locally in space and time with a uniform constant, after $\tau_0$: there exists $C=C(\alpha,\tau_0, \vo)>0$
such that  $ \|v\|_{C^{1,\alpha}([\tau, \tau+1]\times \rn)} \leq C $  for all $\tau>\tau_0\,$.
\end{enumerate}
Then, for $D$ is as in {\it (ii)}, and for any  $q\in\left(N\frac{(p-1)}{p}, \infty\right]$  the following limit holds
\begin{equation}\label{fin1}
      \left\|\frac{v(\tau)-\V_{D}}{\V_{D}}\right\|_{\LL^q(\rn)} \xrightarrow[\tau\to\infty]{}0\,.
\end{equation}
\end{theo}

\begin{rem}\textit{Range(s) of validity and open questions. }\rm In the case $2\leq N < 6$, the entire range $p\in(1,2)\setminus\{\frac{3}{2}\}$ is covered. When $N\ge6$, condition~\eqref{integrability.condition} is equivalent to $p\in (1, p_1)\cup (p_2, 2)$ with $p_1$, $p_2$ as in~\eqref{p1p2}. Therefore, Theorem~\ref{theo:meta1} not only extends the range of admissible values of $p$ close to $2$ from $(p_c,2)$ to $(p_2,2)$, but also covers a very singular range close to $1$. The range $[p_1,p_2]$ is not covered in high dimensions  and we believe that our methods cannot be extended to this range. This is a difficult and intriguing open question that requires new ideas.
\end{rem}

We shall explain the role of each assumption of Theorem~\ref{theo:meta1} in the next remark.

\begin{rem}\label{rem.hyp.intro}\rm
 Assumption {\it (i)} is needed to guarantee both the boundedness of the entropy functional and the integrability of the difference between the solution $v$ and some pseudo-Barenblatt profile $V_D$.

 Assumptions {\it (ii)} and {\it (iii)} are needed to select the limit profile. In the good range $p_c<p<2$, the selection results from the conservation of mass. In the very fast diffusion range $p\leq p_c$, it is a priori not obvious how to indicate to which profile a solution $v$ converges, even if a large time limit exists.  Assumption {\it (ii)} plays a surrogate of this property. Even if two solutions are not in $L^1$, but the difference is, then it is preserved along the evolution (possibly under further assumptions). Therefore, it is called the \textit{conservation of relative mass} \cite{Blanchet2009,Bonforte2010a,Bonforte2017b}.  This may seem surprising, or even false, at first, but we consider it natural since, in the case of the classical fast diffusion equation, it holds merely under the analogous of assumption {\it (i)}, see~\cite[Theorem 1]{Blanchet2009}. On the other hand, assumption {\it (iii)} is the basic assumption needed for the entropy method to work, as explained above.

Assumption {\it (iv)} is used to guarantee the existence of a limit profile. The proof relies on the Ascoli--Arzel\'a Theorem. While this may seem a lot to ask, let us attempt to convince the reader that this is natural. The first argument is that  {\it (iv)} always holds true in the good range $p_c<p<2$. In this case, the uniform $C^{1,\alpha}$-regularity properties of solutions to~\eqref{CPLE} can be easily transferred to solutions to~\eqref{eq:v}, since when $p_c<p<2$ the change of variable is not singular nor degenerate in time, see ~\cite[Lemmata~2.2 and~2.3]{Agueh2009} for a complete proof. In the range $1<p\leq p_c$, solutions  to~\eqref{CPLE} also enjoy good regularity properties (at least the class of solutions that we treat here), but it is not possible to transfer the regularity properties to the solutions to~\eqref{eq:v}, since the change of variables is singular. We consider this a technical difficulty, not a strong impediment. Indeed, in order to obtain {\it (iv)}, one should provide regularity estimates for solutions to~\eqref{CPLE} that are uniform up to the extinction time. We believe this to be true, but could not find a relevant reference in the literature.
\end{rem}

\begin{rem} \rm
The proof of Theorem~\ref{theo:meta1} is based on the entropy method, and its assumptions can be considered the minimal set of hypotheses for our method to work. In the good regime some of them are automatically satisfied, while in the very fast regime they narrow the class of data, indeed compactly supported functions are not allowed. In the present range, inequality~\eqref{entropy-entropy-production-inq.intro} holds only for solutions close to the chosen Barenblatt profile $\VD$ and its proof is based on a quantitative comparison with the linearized counterpart, i.e. the Hardy--Poincaré type inequalities, see Proposition~\ref{prop:hp}.
\end{rem}

\begin{rem}\rm  We want to stress the fact that, when considering radial solutions which satisfies \textit{(i)} of Theorem~\ref{theo:meta1} and inequality~\eqref{H-derivatives}, then assumptions \textit{(ii)} and \textit{(iii)} are automatically verified, see Section~\ref{sec:jus_radial}.
\end{rem}

We conclude the consideration on Theorem~\ref{theo:meta1} by commenting on the norms used in the convergence result~\eqref{fin1}.
\begin{rem}\rm
Let us point out that in the good fast diffusion range $p_c<p<2$, both Theorem~\ref{theo:meta1} and Theorem~\ref{theo:RECR} hold. Theorem~\ref{theo:meta1}  may seem stronger than  Theorem~\ref{theo:RECR}, since we can obtain $\LL^q$ convergence (not only $\LL^\infty$) of the relative error $\frac{v(\tau, \cdot)-V_D(\cdot)}{V_D(\cdot)}$, which means a certain further control on the tails of the relative error. Nonetheless, this is not exactly the case. In fact, for $p_D<p<2$ assumption {\it (i)} of Theorem~\ref{theo:meta1} is much stronger than {\it (ii)} of Theorem~\ref{theo:RECR} due to Lemma~\ref{decay.lemma}. By the comparison principle, we can infer that $\frac{v(\tau,\cdot)-V_D(\cdot)}{V_D(\cdot)}$ has an integrable tail in this range. This is impossible for $p_c<p<p_D$, since the relative error is not integrable in general.  On the other hand, for $p_c<p\leq p_D$ hypothesis {\it (i)} of Theorem~\ref{theo:meta1} is exactly {\it (ii)} in Theorem~\ref{theo:RECR} written after the change of variables~\eqref{change:rho-u}. In this range, we do not need to impose  {\it (ii)--(iv)} since they are automatically verified along the flow. In turn, in the range $p_c<p\leq p_D$, the convergence~\eqref{fin1} yields a slight improvement compared to Theorem~\ref{theo:RECR}, in terms of $\LL^q$-control of the tails of the relative error. We find it interesting to state the result precisely in the current way to stress that {\it (i)--(iv)} are the only ingredients needed for the entropy method to work in the entire range~\eqref{integrability.condition} except for $p=3/2$, when $\cE$ is simply not defined ($\gamma=0$) and another entropy should be employed.
\end{rem}

Once one obtains the convergence result of~\eqref{fin1}, the natural question is whether the convergence holds with a~rate. We answer this question with the following result.
\begin{theo}\label{theo:meta}
Under the same assumptions of Theorem~\ref{theo:meta1}, assume furthermore that, for $D$ is as in {\it (ii)} of Theorem~\ref{theo:meta1},  there exist $\varepsilon_0,\tau_0,\kappa>0$ such that
\begin{equation}\label{gradient.assumption}
|\nabla v^{\gamma-1}(\tau,y)|\le \kappa\left(\varepsilon_0 + |\nabla \V_{D}^{\gamma-1}(y)|\right)\qquad\text{for all } y \in\rn\ \text{ and }\ \tau\geq \tau_0\,.
\end{equation}
Then for any $q\in\left(N\frac{(p-1)}{p}, \infty\right]$ exists $\tau_1=\tau_1(\vo, q)>0$ and $\mathfrak{K}=\mathfrak{K}(\vo, N, p, q)>0$ such that
\begin{equation}\label{fin}
      \left\|\frac{v(\tau,\cdot)-\V_{D}(\cdot)}{\V_{D}(\cdot)}\right\|_{\LL^q(\rn)} \le e^{-\mathfrak{K}\, \tau} \qquad\text{
for all $\tau\ge \tau_1$.}
\end{equation}
\end{theo}

\begin{rem}
\rm We believe that assumption~\eqref{gradient.assumption} does not need to be imposed under the conditions of Theorem~\ref{theo:meta1}, as it does hold at the level of equation~\eqref{CPLE}, see inequality~\eqref{space.gradient.decay}. However, in the range $p<p_c$, we were not able to obtain estimates like~\eqref{gradient.assumption} from properties of solutions to~\eqref{CPLE} due to the lack of the regularity estimates that depend on the extinction time.
\end{rem}

\noindent\textit{Convergence of radial derivatives in the very fast diffusion range}. As in the good fast diffusion range, when~\eqref{Q3} is (at least partially) resolved, it is natural to ask whether radial derivatives of solutions converge towards those of the pseudo-Barenblatt profile, i.e. question~\eqref{Q2} in the context of the very fast diffusion range. Generically, once a result as Theorem~\ref{theo:meta} is proven, this question could be settled by applying regularity results and interpolation inequalities (as in~\cite{Blanchet2009}). However, as it is explained in Remark~\ref{rem.hyp.intro}, in the present setting good regularity results are missing. We have, therefore, decided to use a different technique which uses a correspondence between radial derivatives of~\eqref{eq:v} (or equivalently~\eqref{CPLE}) with those of a weighted diffusion equation of porous medium type, see Section~\ref{sec:conv-radial}. In our analysis there appear again the exponent $p_Y$, defined in~\eqref{pY}, this time as a threshold exponent. We have decided to call this exponent after Yamabe since, in this case, the radial derivative of~\eqref{CPLE} solves the equation $\partial_t \Phi = \Delta \Phi^\frac{N-2}{N+2}$ (see Section~\ref{sec:conv-radial} for more information). This equation is related to the Yamabe flow for a conformally flat metric, see~\cite{delPinoSaez2001,Daskalopoulos2008, Ye94, Blanchet2009}.  The value $p_Y$ is also a sharp threshold for the gradient regularity of solutions when no extra assumptions are imposed, see Section~\ref{ssec:reg} and~\cite[Section 21.3]{DiBenBook}. We stress also that $\max\{p_Y,p_2\}<p_c$, but the relation between $p_Y$ and $p_2$ depends on the dimension. Recall that $p'=p/(p-1)$. Let us now present our main result.

\begin{theo}[\textbf{Convergence in relative error for radial derivatives}]\label{relative.error.radial.derivatives}
Let $N\ge2$, $1 < p \leq p_c$, and let $v$ be a solution to~\eqref{eq:v} with an initial datum $0\leq \vo\in  C^2(\RN)$ that is radially decreasing. Suppose that there exist $D_1, D_2>0$ for which it holds \eqref{stronger.GHP.fast.range}. Assume further that one of the following conditions {\it (i)--(iv)} is satisfied
\begin{enumerate}[{\it (i)}]
\item $N=2$ and $1=p_Y<p\leq p_c$,
\item $2< N < 6$ and $p_Y\leq p \leq p_c$,
\item $N\geq 6$ and $p_2< p \leq  p_c $,
\item $N\geq 6 $, $p_Y\leq p \leq p_2$, and there exist $\wt{D}>0$ and $f\in \LL^1((0,\infty), r^{n-1}\dr)$ with $n=2(1+N/p')$, such that
\begin{equation}\label{D.eqv}
  \partial_r \vo(r) = \partial_r \V_{\wt{D}}(r) + r^\frac{1}{p-1}\, f(r^\frac{p}{2(p-1)})\quad\forall r\ge0\,.
\end{equation}
\end{enumerate} Then there exists  $D=D(\vo)>0$,  $\tau_\diamond=\tau_\diamond(\vo) >0$, $C_\diamond=C_\diamond(\vo, D_1, D_2, N, p)>0$  and $\lambda=\lambda(N,p)>0$  such that
\begin{equation*}
\Big\|\frac{\partial_r v(t)-\partial_r \VD}{\partial_r \VD}\Big\|_{L^\infty(\RN)}+\quad\Big\|\frac{v(t)-\VD}{\VD}\Big\|_{L^\infty(\RN)}\leq C_\diamond\, e^{-\lambda\,\tau}\,\,\quad\forall\,\, \tau>\tau_\diamond\,.
\end{equation*}
Moreover, if $N\geq 6 $ and $p_Y\leq p \leq p_2$, then $D=\wt{D}$.
\end{theo}

We notice that, when $N\ge7$, we have that $p_Y<p_2$, so the result of Theorem~\ref{relative.error.radial.derivatives} extends the results of Theorem~\ref{theo:meta1} for radially decreasing solutions. The proof of Theorem~\ref{relative.error.radial.derivatives} does not make use of the entropy method directly for the solution itself but rather for its radial derivatives, which solves a similar, thought weighted diffusion equation. We refer again to Section~\ref{sec:conv-radial} for more information.\newline

\begin{figure}[ht]
    \centering
\begin{tikzpicture}
\draw [thick, ->] (0,0) -- (15,0);
\draw (1,-0.1)--(1,.1);
\node at (1,-0.4) {1};
\draw (14,-0.1)--(14,.1);
\node at (14,-0.4) {2};
\draw (8,-0.1)--(8,.1);
\node at (8,-0.4) {$\frac{3}{2}$};
\draw (4,-0.1)--(4,.1);
\node at (4,-0.4) {$p_1$};
\draw (8.4,-0.1)--(8.4,.1);
\node at (8.4,-0.4) {$p_Y$};
\draw (8.9,-0.1)--(8.9,.1);
\node at (8.9,-0.4) {$p_2$};
\draw (9.6,-0.1)--(9.6,.1);
\node at (9.6,-0.4) {$p_c $};
\draw (10.5,-0.1)--(10.5,.1);
\node at (10.5,-0.4) {$p_M$};
\draw (12,-0.1)--(12,.1);
\node at (12,-0.4) {$p_D$};
\end{tikzpicture}
    \caption{{\bf Special values of parameter $p$ when $N\geq 7$.} For $p>p_D$ the entropy is displacement convex, for $p>p_M$ Barenblatts have finite weighted ${p'}$-moments, for $p>p_c$ Barenblatt solutions are integrable and the mass is conserved, for $p\in(1,p_1)\cup(p_2,2)$ difference of two Barenblatt solutions is integrable, $p_Y$ is a gradient regularity threshold, for $p=\frac{3}{2}$ we have $\gamma=0$ and the entropy functional is not defined.}
    \label{fig:special-p}
\end{figure}
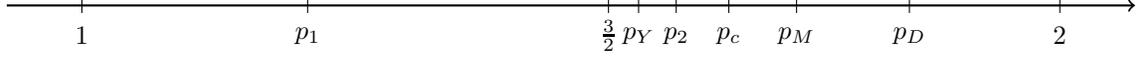

\noindent\textit{Final remarks for the very fast diffusion case.}
As we have seen in the previous section, the global panorama for the very fast diffusion range is quite complex. We recall here the main relation between the dimension $N$ and the critical exponents that appear in this range, see Figure~\ref{fig:special-p}:
\begin{enumerate}
    \item if $N=2$, then $1=p_Y<p_c<\frac 32<p_M<p_D<2$,
    \item if $N=3$, then $1<p_Y<p_c=\frac 32<p_M<p_D<2$,
    \item if $N\in\{4,5\}$, then $1<p_Y<\frac 32<p_c<p_M<p_D<2$,
    \item if $N=6$, then $1<p_1<\frac 32=p_Y=p_2<p_c<p_M<p_D<2$,
    \item if $N\geq 7$, then $1<p_1<\frac 32<p_Y<p_2<p_c<p_M<p_D<2$.
\end{enumerate}
\color{black}

We notice again that exponent $p_Y$ plays an important role in any dimension, while $p_1$ and $p_2$ appear only when $N\ge6$. Therefore it makes sense in our analysis to distinguish between these cases. When $2\leq N < 6$, we recall that Theorem~\ref{theo:meta1} gives the sufficient condition for an answer to question~\eqref{Q3}, and under the additional assumption~\eqref{gradient.assumption}, rates are provided in Theorem~\ref{theo:meta}. Such conditions are verified in the radial case, see Section~\ref{sec:jus_radial}. When it comes to the to convergence of radial derivatives in relative error, i.e. Theorem~\ref{relative.error.radial.derivatives}, the results hold up to $p_Y$. The second important exponent is $p=\frac{3}{2}$, when the entropy (at least in the form~\eqref{cal-E}) is not defined anymore and we cannot adapt the results of Theorem~\ref{theo:meta}. It is important to notice that $p_c$ could be smaller than $\frac{3}{2}$ (in $N=2$), equal to $\frac{3}{2}$ (in $N=3$), or larger than $\frac{3}{2}$ (in $N\in\{4,5\}$), which complicates the graphical representation. We also stress that, for any $D_2,D_1>0$ we have $\V_{D_2}-\V_{D_1}\in\mathrm{L}^1(\mathbb{R}^N)$ for any $1<p<2$ and $N\in\{2,3,4,5\}$. In such a case, and under assumption~\eqref{stronger.GHP.fast.range} the relative entropy $\cE[v|\VD]\in\mathrm{L}^\infty(\tau_0, \infty)$ for any $\tau_0>0$. The conservation of relative mass, i.e. property \textit{(ii)} in Theorem~\ref{theo:meta}, is expected to hold under other assumptions of Theorem~\ref{theo:meta}, as it holds for solutions whose initial data is radial-decreasing, see Section~\ref{sec:jus_radial}. It is unknown whether the Fisher information is bounded along the flow, or not. However, this property holds in the case of radially decreasing initial data, see Section~\ref{sec:jus_radial}. Inequality~\eqref{entropy-entropy-production-inq.intro}, i.e. the fact that the Fisher information controls the entropy, can be obtained by comparison with the Hardy--Poincaré inequality, see Proposition~\ref{prop:hp}, under the additional assumption~\eqref{gradient.assumption}.  When $N\ge6$, the situation is similar but a major difference: the exponents $p_1$ and $p_2$ appears and $ \V_{D_2}-\V_{D_1}\in\mathrm{L}^1(\mathbb{R}^N)$ holds only when $p\in\left(1, p_1\right)\cup \left(p_2, 2\right)$. Consequently, the conservation of relative mass and the boundedness of the entropy, could hold only in this last range. At the same time, when $p\in\left[p_1,p_2\right]$ the Hardy--Poincaré inequality of Proposition~\ref{prop:hp} does not hold any more. Instead, when $p\in\left(p_1, p_2\right)$ a Hardy-type inequality holds, see~\cite{Chlebicka2022}. In Table~\ref{tab:veryfast} we resume such properties for $N\ge7$.  We notice that the acronym ``RDDI'' stand for ``radial, decreasing whose derivatives satisfy a specific inequality''. Such a class of data will be fully described in Section~\ref{sec:jus_radial}. The case $N=6$ is very similar with the only caveat that $p_2=p_Y=\frac{3}{2}$.

\begin{table}[h!]
\centering
\resizebox{\textwidth}{!}{
\begin{tabular}{|c|c|c|c|c|c|c|}
\hline
\textbf{Range} & $\V_{D_2}-\V_{D_1}\in \mathrm{L}^1(\mathbb{R}^N)$ ? & $\cE[v(\tau)|\VD]$  & $\cI[v(v)| \VD]$ & $\tfrac{\d}{\d\tau}\cE[v(\tau)]=-\cI[v(v)]$ & \makecell{Inequality~\eqref{entropy-entropy-production-inq.intro} \\  $\cI[v(v)| \VD] \ge c\,\cE[v(\tau)|\VD] $}  & \makecell{Convergence in relative \\ error for radial derivatives}  \\
\hline
$p_c\geq p > p_2$ & Yes, Lemma~\ref{decay.lemma} & \makecell{$\mathrm{L}^\infty(\tau_0, \infty)$ \\ if \textit{(ii)}\\ of Theorem~\ref{theo:meta1}} & \makecell{Expected $\mathrm{L}^1(\tau_0, \infty)$ \\ Holds for RDDI} & \makecell{Required as~\textit{(iii)} \\ in Theorem~\ref{theo:meta1}} & \makecell{Hardy--Poincaré inequality \\ Prop.~\ref{prop:hp} and~\cite{Chlebicka2022}} & \makecell{Theorem~\ref{relative.error.radial.derivatives} \\under~\eqref{D.eqv}}  \\
\hline
$p_2\geq p\geq p_Y$ & No & $\notin\mathrm{L}^\infty(\tau_0, \infty)$  & Unknown & Unknown & \makecell{ Only Hardy's inequality is \\ available, see~\cite{Chlebicka2022}}  & \makecell{Theorem~\ref{relative.error.radial.derivatives} \\under~\eqref{D.eqv}}  \\
\hline
$p_Y> p> \frac{3}{2}$ & No & $\notin\mathrm{L}^\infty(\tau_0, \infty)$  & Unknown & Unknown & \makecell{ Only Hardy's inequality is \\ available, see~\cite{Chlebicka2022}} & Unknown \\
\hline
$\frac{3}{2}> p \geq p_1$ & No & $\notin\mathrm{L}^\infty(\tau_0, \infty)$ & Unknown & Unknown & \makecell{ Only Hardy's inequality is \\ available, see~\cite{Chlebicka2022}}  & Unknown \\
\hline
$p_1> p> 1$  & Yes, Lemma~\ref{decay.lemma} & \makecell{$\mathrm{L}^\infty(\tau_0, \infty)$ \\ if \textit{(ii)}\\ of Theorem~\ref{theo:meta1}}  & \makecell{Expected $\mathrm{L}^1(\tau_0, \infty)$ \\ Holds for RDDI }  & \makecell{Required as~\textit{(iii)} \\ in Theorem~\ref{theo:meta1}}  & \makecell{Hardy--Poincaré inequality \\ Prop.~\ref{prop:hp} and~\cite{Chlebicka2022}} & Unknown \\
\hline
\end{tabular}
}
\caption{Summary of the main properties of solutions to~\eqref{eq:v} in the very fast diffusion range for $N\ge7$.}
\label{tab:veryfast}
\end{table}

\smallskip
Let us now focus on the organization of the present article. \newline

\noindent{\bf Organization. } In Section~\ref{sec:results_in_original_variables}, we present some of our results proven for~\eqref{eq:v} in the original variables, i.e. for solutions to~\eqref{CPLE}.  In Section~\ref{sec:prelim}, we present basic background and known results related to problems like~\eqref{CPLE}. Section~\ref{sec:conv-mild} is devoted to proving convergence in the relative error, under the assumption that the solutions converge in $\LL^1$. The proof of one of our main results, namely Theorem~\ref{theo:RECR}, is given in Section~\ref{sec:conv-rel-er}. In Section~\ref{sec:conv-radial}, we address question~\eqref{Q2} and focus on results for radial solutions. This section also contains the proofs of Theorem~\ref{Thm:gradient.convergence} and Theorem~\ref{relative.error.radial.derivatives}. Section~\ref{ssec:proof_of_theorem_2} includes both the proof of Theorem~\ref{theo:Entropy-decay} and the main intermediate steps required for it. The convergence results for small values of $p$, namely Theorems~\ref{theo:meta1} and~\ref{theo:meta}, are presented in Section~\ref{sec:comm}. In Section~\ref{sec:jus_radial}, we show that the main assumptions for Theorems~\ref{theo:meta1} and~\ref{theo:meta} are automatically satisfied for radially decreasing data. We conclude with Section~\ref{sec:open}, where we list several open problems. Lastly, in Appendix, we prove some auxiliary lemmas and summarize the role of the most important parameters used throughout the paper.

\section{Results in the original variables}\label{sec:results_in_original_variables}

While it is easier to state our main results in the variables introduced in~\eqref{change:rho-u} (or equivalently in~\eqref{change.variables}), it may be not simple to understand such results for the original equation~\eqref{CPLE}. This is why we have decided to state such results (with the exception of Theorems~\ref{theo:Entropy-decay},~\ref{theo:meta1}, and~\ref{theo:meta}) for the original problem~\eqref{CPLE}. We notice that the excluded theorems could also be stated for the original problem. However, they make an extensive use of the entropy and the Fisher information which are objects that, for their intrinsic nature, are much better understood in terms of problem~\eqref{eq:v}. For instance, taking the relative entropy with respect to $\B_D(t)$ instead of $\VD$ could be done, but it would make all the computations much more delicate and unnecessarily complicated (for instance one should also differentiate $\B_D(t)$ in time when investigating the derivative of the entropy).

Let us now comment on our results in the original variables. In the good fast diffusion range $p_c<p<2$, the translation of Theorems~\ref{theo:RECR} and~\ref{Thm:gradient.convergence} do not offer any particular difficulty. Indeed, Theorem~\ref{theo:RECR} reads then as follows.
\begin{theo}
Let $N\ge3$, $0\le \uo \in \LL^1(\rn)\cap\X$,  and $M:=\|\uo\|_{\LL^1(\rn)}>0$. Assume $u$ is a weak solution to~\eqref{CPLE} with initial datum $\uo$. Suppose one of the following holds:
\begin{enumerate}[{\it (i)}]
\item $p_M<p<2$;
\item $p_c<p\le p_M$ and there exists $M_2>M_1>0$ and $T>0$ such that
\begin{equation}\label{stronger.GHP}
\B_{M_1}(T,x) \le \uo(x) \le \B_{M_2}(T, x)\qquad \forall\,x\in\rn\,.
\end{equation}
\end{enumerate}
Then there exist $T_\star=T_\star(p,N,M, \|\uo\|_{\X})>0$, $K_\star=K_\star(p,N,M, \|\uo\|_{\X})>0$ and $\sigma=\sigma(p,N)>0$ such that
\begin{equation*}
  \left\|\frac{u(t,\cdot)-\B_M(t,\cdot)}{\B_M(t,\cdot)}\right\|_{\LL^\infty(\rn)} \le K_\star\,t^{-\sigma}\qquad\forall\, t\ge T_\star\,.
\end{equation*}
\end{theo}

\begin{rem} \rm
We notice that in the case $N=2$, Theorem~\ref{theo:RECR} still holds as long as $p\neq \frac{3}{2}$.
\end{rem}

While, in the case of Theorem~\ref{Thm:gradient.convergence}, we have the following.
\begin{theo}\label{original-variables:p>pc}
Let $N\ge 2$, $p_c<p<2$, $0\le \uo\in \LL^1(\rn) \cap C^2(\R^N) $ be radially symmetric and $M:=\|\uo\|_{L^1(\RN)}>0$ and $u$ is a solution to\eqref{CPLE} with datum $\uo$. Suppose one of the following holds:
\begin{enumerate}[{\it (i)}]
\item $p_M<p<2$ and there exist $A>0$ and $R_0>0$ satisfying
\begin{equation*}
  \partial_r \uo(r)\leq 0\quad\mbox{and}\quad |\partial_r \uo(r)|\leq A\, r^{-\frac{2}{2-p}}\ \ \forall r\ge R_0\,.
\end{equation*}
\item $p_c<p\leq p_M$ and there exist $M_2>M_1>0$  and $T>0$ such that
\begin{equation*}
    \partial_r \B_{M_2}(T,r) \leq \partial_r \uo(r) \leq \partial_r \B_{M_1}(T,r)\quad \forall\,\, r\ge0\,.
\end{equation*}
\end{enumerate}
Then there exist $t_\star>0$, $k_\star>0$ and $\lambda=\lambda(p,N)>0$ such that
\begin{equation*}\
\Big\|\frac{\partial_r u(t,\cdot)}{\partial_r \B_M(t,\cdot)}-1\Big\|_{L^\infty(\RN)}\le k_\star\,t^{-\lambda}\,\quad\forall\,t\ge t_\star\,,
\end{equation*}
where $\partial_r u$ (resp. $\partial_r \B_M$) is the radial derivative of $u$ (resp. $\B_M$). When $p_M<p<2$ then $ t_\star=t_\star(\uo, A, R_0, M, p, N)$  and $k_\star=k_\star(\uo, A, R_0, M, p, N)$, while when $p_c<p\leq p_M$ then $ t_\star=t_\star(\uo, M_2, M_1, T, p, N)$  and $k_\star=k_\star(\uo, M_2, M_1, T,  p, N)$.
\end{theo}
\begin{rem}
    \rm We notice that the regularity assumption on the initial datum, i.e., $\uo\in C^2(\rn)$ is imposed for the sake of simplicity of exposition. More precisely, the regularity threshold for our proof to work is $C^{1, \alpha}(\rn)$ for $\alpha>0$ depending on $p$. The same remark applies also to Theorems~\ref{radial.derivatives.pc} and~\ref{relative.error.radial.derivatives}.
\end{rem}

Things are more complicated when we consider the very fast diffusion range $1<p\leq p_c$ and notably Theorem~\ref{relative.error.radial.derivatives}. Indeed, since when $p=p_c$ solutions and the pseudo-Barneblatt of~\eqref{pseudo-Barenblatt} are non-negative for all times $t\in\left(0, \infty\right)$ while when $1<p<p_c$ they extinguish in a finite time $T>0$, Theorem~\ref{relative.error.radial.derivatives} could not be translated easily. Therefore, we have chose to divide into two results, one when $p=p_c$ and the other when $1<p<p_c$. In the former case, we have the following.

\begin{theo}[\textbf{Convergence in relative error for radial derivatives for $p=p_c$}]\label{radial.derivatives.pc}
Let  {$N\ge2$}, $p= p_c$, and $u$ be a solution to~\eqref{CPLE} with an initial datum $0\leq \uo\in  C^2(\RN)$, which is radial and decreasing. Suppose that there exist $D_1, D_2>0$  and $T>0$, such that
\begin{equation}\label{hp.derivative}
  \partial_r \B_{D_2,T}(0,r) \leq \partial_r \uo(r) \leq \partial_r \B_{D_1,T}(0,r)\,\quad\forall r\ge0\,.
\end{equation}
Then there exists  $D=D(\uo)>0$,  $T_\diamond=T_\diamond(\uo) >0$, $C_\diamond=C_\diamond(\uo, D_1, D_2, N, p)>0$  and $\lambda=\lambda(N,p)>0$ such that
\begin{equation}\label{convergence.rate.pc}
\Big\|\frac{\partial_r u(t,\cdot)}{\partial_r \B_{D,T}(t, \cdot)}-1\Big\|_{L^\infty(\RN)}+\Big\|\frac{u(t,\cdot)}{\B_{D,T}(t, \cdot)}-1\Big\|_{L^\infty(\RN)}\leq C_\diamond\, t^{-\lambda}\,\,\quad\forall\,\, t>T_\diamond\,.\end{equation}
\end{theo}
While, in the case $1<p<p_c$ we have the following.
\begin{theo}[\textbf{Convergence in relative error for radial derivatives for $p<p_c$}]\label{relative.error.radial.derivatives.original}
Let $N\ge2$, $1 < p < p_c$, and let $u$ be a solution to~\eqref{CPLE} with an initial datum $0\leq \uo\in  C^2(\RN)$ radial and decreasing. Suppose that there exist $D_1, D_2>0$ and $T>0$ it holds \eqref{hp.derivative}. Assume further that one of the following conditions {\it (i)--(iv)} is satisfied
\begin{enumerate}[{\it (i)}]
\item $N=2$ and $1=p_Y<p<p_c$,
\item $2< N < 6$ and $p_Y\leq p<p_c$,
\item $N\geq 6$ and $p_2< p <  p_c $,
\item $N\geq 6 $, $p_Y\leq p \leq p_2$, and there exist $\wt{D}>0$ and $f\in \LL^1((0,\infty), r^{n-1}\dr)$ with $n=2(1+N/p')$, such that
\begin{equation}\label{D}
  \partial_r \uo(r) = \partial_r \B_{\wt{D},T}(0,r) + r^\frac{1}{p-1}\, f(r^\frac{p}{2(p-1)})\quad\forall r\ge0\,.
\end{equation}
\end{enumerate} Then there exists  $D=D(\uo)>0$,  $T_\diamond=T_\diamond(\uo) >0$, $C_\diamond=C_\diamond(\uo, D_1, D_2, N, p)>0$  and $\lambda=\lambda(N,p)>0$  such that
\begin{equation}\label{convergence.rata.below.pc}
\Big\|\frac{\partial_r u(t,\cdot)}{\partial_r \B_{D,T}(t, \cdot)}-1\Big\|_{L^\infty(\RN)}\leq C_\diamond\, (T-t)^{-\lambda }\quad\mbox{and}\quad\Big\|\frac{u(t,\cdot)}{\B_{D,T}(t, \cdot)}-1\Big\|_{L^\infty(\RN)}\leq C_\diamond\, (T-t)^{-\lambda}\,\,\quad\forall\,\, t_\diamond<t<T\,.
\end{equation}
Moreover, if $N\geq 6 $ and $p_Y\leq p \leq p_2$, then $D=\wt{D}$.
\end{theo}

\section{Preliminary information}\label{sec:prelim}

\subsection{Notation } Following a usual custom, we denote by $c$ a general positive constant. Different occurrences from line to line will be still denoted by $c$, while special occurrences will be denoted by $c_1, c_2,  \tilde c$ or similar. Relevant dependencies on parameters will be emphasized using parentheses, i.e., $c= c(p,M)$ means that $c$ depends on $p$ and $M$. We define $(u)_{+}:=\max\{0, u\}$.

We also recall the definitions of the constants $b_1$ and $b_2$ which appear in the definition of the Barenblatt function~\eqref{beta}, see also~\cite{Bonforte2020b} for more information:
\begin{equation}
    \label{b2} b_2:=\tfrac{2-p}{p}\left(p-N(2-p)\right)^{-\frac{1}{p-1}},
\end{equation}
while $b_1$ is such a positive constant that
\begin{equation}
    \label{b1}\int_\rn \big(b_1+b_2|x|^{p'}\big)^{-\frac{p-1}{2-p}}=1\,.
\end{equation}

\subsection{ Existence and uniqueness}

Let us first introduce the concept of non-negative weak solutions that we shall use throughout the present work and comment on their well-posedness.
\begin{defi} We say that $u$ is a non-negative weak solution  to problem \eqref{CPLE}  on $(0, \infty)\times\RN$ for $1<p < 2$  with non-negative initial data $\uo \in \LL^1_{{\rm loc}}(\RN)$ if $u\in \LL^p((0,\infty); W^{1,p}_{\rm loc}(\RN)) \cap C((0, \infty);\LL^1_{{\rm loc}}(\RR^N))$ and
\begin{align*}
\int_{\RN}u(s,x)\phi(s,x)\dx = \int_{\RN}u(t,x)\phi(t,x)\dx
&+ \int_s^t\int_{\RN}\big(-u(\tau,x)\partial_\tau \phi(\tau,x) + |\nabla u(\tau,x)|^{p-2}\nabla u(\tau,x) \cdot \nabla \phi(\tau,x) \big)\dx\dtau\,,
\end{align*}
for all $t>s>0$ and for all functions $\phi \in C^\infty([0,+\infty)\times \RN)$, such that the support of the maps $x\mapsto \phi(t,x)$ is compact for any $t\ge0$. The initial data is attained in the following sense
 $$\lim_{t\to 0} \int_{\RN}u(t,x)\varphi(x) \dx =\int_{\RN} \uo(x) \varphi(x) \dx \quad \forall \varphi \in C_c^{\infty}(\RN)\,. $$
\end{defi}

\subsection{ Boundedness and regularity of solutions} \label{ssec:reg} Since the literature concerning the regularity properties of solutions to equations like~\eqref{CPLE} is abundant, we do not aim to describe the state of art exhaustively. Instead, we shall restrict ourselves to presenting only the background needed for our study.

Here we shall mainly focus on the concept of local weak solution (not necessarily non-negative ones) which differ from our definition above by assumptions on the integrability properties of $u$, namely for a solution defined on $(0, T)\times \Omega$ one asks typically that $u\in\LL^\infty_{\rm loc}((0,T);\LL^2_{\rm loc}(\Omega))$ and $|\nabla u|\in \LL^p_{\rm loc}((0,T);\LL^p_{\rm loc}(\Omega))$. To the best of our knowledge, these kind of solutions have been studied first in~\cite{zbMATH03836584} where continuity of $\nabla u$ has been proven, with an explicit modulus of continuity, for $p>\max\{1,p_Y\}=\max\{1,\frac{2N}{N+2}\}$. In the same range of $p$, the H\"older continuity of the gradient has been obtained in~\cite{zbMATH03875840} (with some mistakes in the computations, as it was pointed out in~\cite{Yazhe_1986}, which have been fully solved in~\cite{zbMATH03957779}), see also the review~\cite{zbMATH04046258}. The threshold value $p_Y=\frac{2N}{N+1}$ is sharp {for the gradient regularity} under no extra assumptions on the solution. Indeed, when $p\leq p_Y$ weak solutions are not bounded, see the discussion~\cite[Section 21.3]{DiBenBook}. However, in the whole range $1<p<2$ it is known, that \emph{bounded} weak solutions are H\"older continuous, see~\cite{Ya_Zhe_1992,Ya_zhe_1988}. In our case, it is known that solutions to~\eqref{CPLE} are bounded provided the initial datum $\uo \in \LL^q_{\rm loc}(\RN)$ for $q> N\frac{(2-p)}{p}$, see for instance~\cite[Section III]{DiBenedetto1990} and~\cite[Theorem 2.1]{Bonforte_2010}.  In conclusion, when $p\in(p_Y,2)$ and the initial datum is integrable enough, weak solutions to~\eqref{CPLE} are bounded and, therefore, the function $(t,x)\mapsto u(t,x)$ is $C^{1,\alpha}_{\rm loc}((0, \infty)\times\RN)$. However, when $p<p_c$, the coefficient $\alpha$ may depend on the function itself other than on $p$ and $N$.

In the present paper, we shall also consider the case $p\leq p_Y$, so let us comment on how to obtain the $C^{1,\alpha}_{\rm loc}(\RN)$ regularity for solutions in this case. Let us stress, however, that the following reasoning holds in the whole range $1<p<2$. The main idea is to use the concept of viscosity solutions (see~\cite{Juutinen_2001} for a precise definition). In our case, when a weak solution is continuous, then it is also a viscosity solution, see for instance~\cite {Juutinen_2001,Siltakoski_2021,Feng_2023,Feng_2022}. Since we shall consider only bounded weak solutions (which are H\"older continuous), the above discussion proves that solutions to problem~\eqref{CPLE} are, indeed, viscosity solutions. The main advantage of employing this notion of solutions is that, in the last decade, there has been a growing interest in obtaining regularity results for viscosity solutions to equations related to~\eqref{CPLE}, see for instance~\cite{Feng_2023} for a detailed bibliography. For what concerns our investigation, the needed result is~\cite[Theorem 1.1]{Imbert_2017}, where the authors prove the function $(t,x)\mapsto u(t,x)$ is $C^{1,\alpha}_{\rm loc}((0, \infty)\times\RN)$ for some $\alpha>0$, which depends on $p$, $N$ and the solution itself.

Several considerations on the integrability of time derivative $\partial_t u$ and of $\Delta_p u = \textrm{div}\left(|\nabla u|^{p-2} \nabla u\right)$ are in order. It is known that for a continuous weak solution both $\partial_t u$ and $ \textrm{div}\left(|\nabla u|^{p-2} \nabla u\right)$ belong to $\LL^2_{\rm loc}(Q)$, where $Q$ is space-time cylinder. Furthermore, equation~\eqref{CPLE} is satisfied almost everywhere in $(t,x)$. We refer to~\cite[Corollary]{Bonforte_2010} and to~\cite{Feng_2022}. We also remark that more is known about the integrability properties of derivatives of $ (t, x)\mapsto |\nabla u|^\frac{p-2+s}{2}\nabla u$ (where $s$ is chosen appropriately), for which we refer to~\cite{Feng_2022,Feng_2023}. When it comes to explicit estimates of the continuity of the gradient of solutions, we refer to~\cite{Min5,  Min3, Min2, Min1}. In those results, the authors obtain the $C^{1,\alpha}$-regularity by exploiting very interesting connections with nonlinear potential estimates. While they are valid mainly when $p>p_c$, those results also apply when the equation~\eqref{CPLE} has a measure as right-hand-side. We refer to~\cite{KuMi-guide} for a general overview.

Since the uniform convergence in relative error is related to Harnack inequalities, let us conclude this subsection with some considerations on them.  The problem of obtaining a precise form of those inequalities  has been a long-standing quest. Indeed, in this nonlinear setting, the intrinsic cylinders depend on the solution itself, showing several differences between the case $p_c<p<2$ and $1<p\leq p_c$.  In the good range, we refer to the paper~\cite{DiBen2008}, while in the whole range $1<p<2$, it has been proven in~\cite{Bonforte_2010}. The Harnack inequalities considered in those two papers are valid for local solutions, i.e., no assumption on boundary data is made. For boundary Harnack inequalities, we refer to~\cite{KuMiNy}. Nowadays, several related results are available, see the monograph~\cite{DiBen2008} and references therein.

\subsection{Comparison principles}
By the comparison principle,  we mean that, in some sense, ordered data generate ordered solutions at all times.
Such results for solutions to~\eqref{CPLE} seem to be well known by experts in the field, cf.~\cite[Sections 3 and 4]{Feo_2021} or~\cite[Section 4.5]{Bogelein:2023aa}. Nonetheless, we could not find references with complete proofs in the case of the Cauchy problem within the whole range $1<p<2$. One of the main difficulties in the proofs of comparison principles is that, at least when $p<p_c$, in general solutions are  not integrable. Thus, a priori the quantity $(u_1-u_2)_{+}$
(where $u_1$ and $u_2$ are two solutions to~\eqref{CPLE}) cannot be used as a test function. However, in our case when solutions are regular and bounded, one is equipped with two comparison principles. We decided to include them with the sketches of the proofs for completeness.

The first comparison principle we present reads: for $u_1$ and $u_2$ being two solutions (that are regular enough, cf. Section~\ref{ssec:reg}) with initial data $u_{1,0}$ and $u_{2,0}$, respectively, we have
\begin{equation}\label{comparison-principle}
  \mbox{if}\quad  u_{1,0}\leq u_{2,0}\,,\qquad\mbox{then}\quad u_1(t,x)\leq u_2(t,x)\quad \forall\,\,t>0\ \,\forall\,\,x\in\RN\,.
\end{equation}
 It can be proven via the construction inspired by~\cite[Chapter~3]{DiBenedetto1990} where the authors constructed solutions to~\eqref{CPLE} by approximation by solutions to the Dirichlet problem. The comparison principle for the Dirichlet problem goes back to~\cite{Lions69}. For a more recent proof we refer to~\cite[Chapter~7, Corollary~1.1]{DiBenBook}, cf. also references therein. We notice that for the Dirichlet problem, there is no restriction on $p>1$. We have $u_1$ and $u_2$ obtained by approximation with solutions to the Dirichlet problem with initial data $u_{1,0}\phi_R$ and $u_{2,0}\phi_R$ (where $\phi_R$ is a cut-off function supported in the ball of radius $R>0$). The approximate solutions are ordered on $B_{2R}$. This relation holds in the limit $R\rightarrow \infty$.

One is not deprived from comparison results in the class of less regular solutions. In fact the following \emph{$\LL^1_{\rm loc}$-comparison principle} holds: let $u_1$ and $u_2$ be two solutions (from the class the class $\Sigma^\star$, see~\cite[Chapter II]{DiBenedetto1990}):
\begin{equation}\label{local.comparison}
  \mbox{if}\quad (u_1-u_2)_{+}\rightarrow 0\,\,\mbox{in $\LL^1_{\rm loc}(\RN)$}\ \ \mbox{as}\ \  t\rightarrow0\,,\qquad\mbox{then}\quad u_1(t,x)\leq u_2(t,x)\quad \forall\,\,t>0\,\,\forall\,\,x\in\RN\,.
\end{equation}
The main advantage of this local comparison is that it avoids using global integrability, even if it assumes an order in $\LL^1_{\rm loc}(\RN)$. It might be proven by following the lines of~\cite[Proposition~II.3.1 and Theorem~II.1.1]{DiBenedetto1990}. Indeed, a careful inspection of the proof of~\cite[Proposition II.3.1]{DiBenedetto1990} shows that under the assumption of~\eqref{local.comparison} one gets that, for all $T>0$, $(u_1(t,x)-u_2(t,x))_+\in\LL^\infty\left((0, T); \LL^q_{\rm loc}(\RN)\right)$ for all $q\in\left[1, \infty\right)$, and there exists $C=C(N, p, q)>0$ such that for all $R>0$, $t\in(0,T)$, and $\sigma\in(0,1)$ one has
\begin{equation}\label{inequality.dib}
   \int_{B_R} \Big(u_1(t,x)-u_2(t,x)\Big)_+^q \dx \leq \frac{C}{(\sigma\,R)^p}\, \int_0^t\int_{B_{(1+\sigma)R}} \Big(u_1(t,x)-u_2(t,x)\Big)_+^{q+p-2} \dx\dt\,.
\end{equation}
Indeed, we also notice that, while originally inequality~\eqref{inequality.dib} is stated for the absolute value of the difference (i.e. $|u_1-u_2|$), its proof is done for the positive part of $u_1-u_2$. Once~\eqref{inequality.dib} is obtained, using the same argument as in the proof of~\cite[Theorem II.1.1, p. 257]{DiBenedetto1990} one obtains that, for every $q>1$, all $t\in\left(0, T\right)$, and $C=C(N,p,q)$ independent of $R>0$:
\begin{equation*}
  \int_{B_R} \Big(u_1(t,x)-u_2(t,x)\Big)_+^q \dx \leq C\,
  t^\frac{q}{2-p} R^{N-q\frac{p}{2-p}}\,.
\end{equation*}
By choosing $q$ so large such that $qp/(2-p)>N$, one obtains in the limit $R\rightarrow0$ that $\int_{\RN} (u_1(t,x)-u_2(t,x))_+^q\dx=0$ for all $t>0$. Consequently, $u_1(t,x)\leq u(t_2,x)$ for all $t>0$ and $x\in\RN$. We  acknowledge that, despite the comparison principle has not been stated in the  form~\eqref{local.comparison} in~\cite{DiBenedetto1990}, it has been used in this form in~\cite[Proposition III.7.1]{DiBenedetto1990}.  Hence, we do not claim any originality for the above result.

\section{Convergence in relative error under a priori convergence in Lebesgue space}\label{sec:conv-mild}
The goal of this section is to obtain an explicit convergence rate towards the Barenblatt profile in the uniform relative error, provided that we know a priori a convergence rate in a weaker norm. In what follows we shall use the $\LL^1$-norm. One can prove a counterpart of this result involving $\LL^q$-norm with $1\le q \le \infty$, by interpolation arguments. We recall here that $u_0\in \X$ if
\[
\|u_0\|_{\X}:=\sup_{R>0}\,R^{\frac{p}{2-p}-N}\,\int_{|x|\ge R} u_0(x)\dx <\infty\,.
\]
\begin{theo}\label{theo:conv-L1-RE}
Let $N\ge1$, $p_c<p<2$, $0\le \uo \in \LL^1(\rn)\cap\X$ and $M=\int_{\rn}\uo(x)\dx>0$. Assume $u$ is a solution to~\eqref{CPLE} with initial datum $\uo$. Suppose that for some $\wt{T}>0$, $\wt{K}>0$ and $N+1\ge\nu>0$,  we have that
\begin{equation}\label{given.convergence}
\|u(t,\cdot)-\B_{\Mo}(t,\cdot)\|_{\LL^{1}(\rn)}\le \wt{K}\,t^{-\nu}\qquad\forall\, t\ge \wt{T}\,.
\end{equation}
Then there exist $K_\star=K_\star(p,N, \Mo, \|\uo\|_{\X}, \wt K)>0$  and $T_\star=T_\star(p,N,\Mo, \|\uo\|_{\X},\wt{T})>0$ such that  we have
\begin{equation*}
\left\|\frac{u(t,\cdot)-\B_{\Mo}(t,\cdot)}{\B_{\Mo}(t,\cdot)}\right\|_{\LL^\infty(\rn)} \le\, K_\star \,  t^{-\frac{\nu(2-p)}{N+1}}\qquad\forall\, t\ge T_\star\,.
\end{equation*}
\end{theo}

As we shall see, the strategy of the proof of the above theorem is to consider separately two different regions in the $(t,x)$ plane: \emph{inner cylinders}, i.e. $\{|x|\le C\, t^\beta\}$ for a constant $C>0$, and \emph{outer cylinders}, i.e. $\{|x|\ge C\,t^\beta\}$. Assumption~\eqref{given.convergence} plays a major rôle in the inner cylinders, while in the case of the outer cylinders it is the global Harnack principle~\eqref{ghp.inq} to imply the wanted result.

\subsection{Properties of solutions to \eqref{CPLE} for $p_c<p<2$ and initial datum $\uo\in\X$}

We stress that in a significant part of our paper (i.e., Sections~\ref{sec:conv-mild}, \ref{sec:conv-rel-er}, and \ref{sec:conv-radial}) we shall consider the exponent $p\in\left(p_c,2\right)$. We notice that within this range of $p$ the results of \cite[Theorem~1.1 and~1.3]{Bonforte2020b} imply the following consequences having fundamental meaning in our reasoning.
\begin{enumerate}[{\it (i)}]
    \item If $0\leq \uo\in\X\setminus\{0\}$, where $\X$ as in~\eqref{UCRE/Xp}, then for any $t_0 >0$ there exist (explicit) constants $\tau_1, \, M_1, \, \tau_2, \, M_2$
such that for all $x\in \RN$ and $t \ge  t_0$ the following upper and lower bounds hold true
\begin{equation}\label{ghp.inq}
  \B_{M_1}(t-\tau_1,x)\, \le \,u(t,x) \,\le \, \B_{M_2}(t+\tau_2,x)\,;
  \end{equation}
  The above inequality is also known as the \emph{global Harnack principle}.
\item Let $\Mo=\|\uo\|_{\LL^1(\rn)}>0$ and $\uo\geq 0$. Then
\begin{equation*}
\lim\limits_{t\rightarrow \infty} \Big\|\frac{u(t,\cdot)}{\B_{\Mo}(t,\cdot)}-1\Big\|_{\LL^\infty(\RN)}= 0\,
\end{equation*}
 if and only if $\uo\in\X\setminus\{0\}$.
\item If $0\leq \uo\in\LL^1(\rn)$, then $\LL^\infty$-norm of the gradient decays in time. More precisely, there exists a~constant $c_1=c_1(p,N)>0$ such that
\begin{equation}\label{time.grandient.decay}
\|\nabla u (t,\cdot) \|_{L^\infty(\RN)}\le\, c_1\, \frac{\|\uo\|_{L^1(\RN)}^{2\beta}}{t^{(N+1)\beta}} \qquad\mbox{for any }t>0\,.
\end{equation}
\item If $0\leq \uo\in\X$, then we can say more about the spacial decay of the gradient. More precisely, there exists a constant $c_2=c_2(N, p)>0$ such that
\begin{equation}\label{space.gradient.decay}
|\nabla u(t,x)| \le c_2 \,\frac{   \|\uo\|_{\LL^1(\RN)}^{2\beta} + \|\uo\|_{\X}^{2\beta}  +t^\frac{2\beta}{2-p} }{\left(1+|x|\right)^{\frac{2}{2-p}}t^{(N+1)\beta}}\qquad\mbox{for any $x\in \RN$ and $t>0$}\,.
\end{equation}
\end{enumerate}

\subsection{Convergence in outer cylinders}
\begin{prop}\label{prop:u-geq-B}
Let $N\ge1$, $p_c<p<2$, $0\le \uo \in \LL^1(\rn)\cap\X$ and $M=\int_{\rn}\uo(x)\dx>0$. Assume $u$ is a weak solution to problem~\eqref{CPLE} with initial datum $\uo$.  Then for any $\varepsilon\in (0,1)$ there exists $\underline{T}(\varepsilon)>0$ and $\underline{\vr}(\varepsilon)>0$ such that
then
\begin{equation*}
u(t,x)\geq \, (1-\varepsilon) \,\B_{\Mo}(t,x)\quad \forall\, |x|\geq \underline{\vr}(\varepsilon)\, t^{\beta}\quad \forall t> \underline{T}(\varepsilon)\,.
\end{equation*}
\end{prop}
The strategy of the proof of the above proposition follows closely the proof of~\cite[Proposition 4.6]{Bonforte2020a}.
\begin{proof}
By inequality~\eqref{ghp.inq}, for $t\geq t_0 = 1$, we have $u(t,x)\geq \B_{M_1}(t-\tau_1,x)$. By integrating~\eqref{ghp.inq}, we find as well that $M_1\le M$. If $M_1=M$, then we conclude that $u=\B_{M}$ and the proposition is proven. Therefore we restrict our attention to the case $M_1<M$. Let us define $\underline{\varepsilon}\in\left(0,1\right)$ such that $(1-\underline{\varepsilon}) M^\beta = M_1^\beta$ for $\beta$ from~\eqref{beta} and let $0<\varepsilon<\underline{\varepsilon}$. In order to prove the claim we need to prove  that for $|x|$ and $t$ large enough it holds
\begin{equation}\label{BB-lower-bound}
\frac{\B_{M_1}(t-\tau_1,x)}{{\B}_{\Mo}(t,x)}\geq 1-\ve\,.
\end{equation}
Let us notice that for $b_1,b_2$ from \eqref{b1} and \eqref{b2}, respectively, it holds that
\begin{align*}
 \frac{\B_{M_1}(t-\tau_1,x)}{{\B}_{\Mo}(t,x)}
&=\left(\frac{ t-\tau_1}{t}\right)^\frac{1}{2-p} \left(\frac{ b_1\left(\frac{t}{\Mo^{2-p}}\right)^{\beta p'} }{ b_1\left(\frac{t-\tau_1}{{M_1}^{2-p}}\right)^{\beta p'}}\right)^\frac{p-1}{2-p}\left(\frac{1 +\frac{b_2}{b_1}\Mo^{(2-p)\beta p'} |x|^{p'}t^{-\beta p'}}{ 1 +\frac{b_2}{b_1}{M_1}^{(2-p)\beta p'} |x|^{p'}(t-\tau_1)^{-\beta p'}}\right)^\frac{p-1}{2-p}\,,
\end{align*}
where $p'=p/(p-1)$.  Upon setting
\begin{equation}
    \label{eta-s-C}
\eta(t):=\left(\frac{t}{t-\tau_1}\right)^\beta,\qquad s(t,x)=|x|^{p'}t^{-\beta p'},\qquad\text{and}\qquad c=\frac{b_2}{b_1}{\Mo}^{(2-p)\beta p'}
\end{equation}
and recalling that $\beta=(p-N(2-p))^{-1}$, the left-hand-side of~\eqref{BB-lower-bound} becomes
\[\frac{\B_{M_1}(t-\tau_1,x)}{{\B}_{\Mo}(t,x)}=\left(\frac{t}{ t-\tau_1}\right)^{N \beta}  \left(\frac{1 +c\,s(t,x)}{\left(1-\overline{\varepsilon}\right)^\frac{p-2}{p-1} +c\,\eta^{p'}(t)\,s(t,x)}\right)^\frac{p-1}{2-p}=\eta^N( t) \left(\frac{1 +cs(t,x)}{\left(1-\overline{\varepsilon}\right)^\frac{p-2}{p-1} +c\eta^{p'}(t)s(t,x)}\right)^\frac{p-1}{2-p}\,.\]
Therefore, inequality~\eqref{BB-lower-bound} is equivalent to
\begin{equation}\label{BB-lower-bound-equivalent}
s(t,x)\geq \frac{ \left(\frac{1-\ve}{1-\underline{\varepsilon}}\right)^\frac{2-p}{p-1}\eta^{-N\frac{2-p}{p-1}}(t)-1}{c\left(1-\eta ^{\frac{1}{\beta(p-1)}}(t)(1-\ve)^\frac{2-p}{p-1}\right)}=:\underline{s}(t, \ve)\,,
\end{equation}
provided that $1>\eta ^{\frac{1}{\beta(p-1)}}(t)(1-\ve)^\frac{2-p}{p-1}$. We restrict our attention to $t>\underline{T}(\ve)$, where \begin{equation}
    \label{under-T}
\underline{T}(\ve):=\max\left\{\frac{\tau_1}{1-(1-\ve)^{2-p}}\,, t_1\,, t_0  \right\}\,
\end{equation}
and $t_1>0$ is such that $\eta(t_1)^\frac{1}{\beta(p-1)}={2}{(1+(1-\ve)^\frac{2-p}{p-1})^{-1}}$. Observe that for $t> \frac{\tau_1}{1-(1-\ve)^{2-p}}$, we have
$\eta ^{\frac{1}{\beta(p-1)}}(t)(1-\ve)^\frac{2-p}{p-1}<1$.
Then~\eqref{BB-lower-bound} will follow from~\eqref{BB-lower-bound-equivalent} and $t> \underline{T}(\ve)$.  Inequality~\eqref{BB-lower-bound-equivalent} holds true as long as $s(t,x)\ge \underline{\vr}^{p'}(\varepsilon)$ for $\underline{\vr}(\varepsilon)$ defined below. Indeed, since $\eta\geq1$ is decreasing and, for  $t\geq t_1$, it holds that $\eta(t)\leq \eta(t_1)$, we have
\begin{align}\label{under-rho}
    \underline{s}(t, \ve)&\leq \frac{\left(1-\underline{\ve}\right)^\frac{p-2}{p-1} (1-\ve)^\frac{2-p}{p-1}}{c\left(1-(1-\ve)^\frac{2-p}{p-1}\eta^\frac{1}{\beta(p-1)}(t_1)\right)}\leq \frac{\left(1-\underline{\ve}\right)^\frac{p-2}{p-1}(1-\ve)^\frac{2-p}{p-1}}{c\left(1-\frac{2(1-\ve)^\frac{2-p}{p-1}}{1+(1-\ve)^\frac{2-p}{p-1}}\right)}=\frac{(1-\ve)^\frac{2-p}{p-1}(1+(1-\ve)^\frac{2-p}{p-1})}{c\,\left(1-\underline{\ve}\right)^\frac{2-p}{p-1}\left(1-(1-\ve)^\frac{2-p}{p-1}\right)}=:\underline{\vr}^{p'}(\varepsilon)\,.
\end{align}
\end{proof}

\begin{prop}
\label{prop:u-leq-B} Let $N\ge1$, $p_c<p<2$, $0\le \uo \in \LL^1(\rn)\cap\X$ and $M=\int_{\rn}\uo(x)\dx>0$. Assume $u$ is a weak solution to problem~\eqref{CPLE} with initial datum $\uo$. Then for any $\varepsilon\in (0,1)$ there exists $\overline{T}(\varepsilon)>0$ and $\overline{\vr}(\varepsilon)>0$ such that
\[
u(t,x)\leq (1+\varepsilon)\,\B_{\Mo}(t,x)\,\quad\forall\, |x|\geq \overline{\vr}(\varepsilon)\,t^\beta\quad \forall\,t\geq\,\overline{T}(\varepsilon)\,.
\]
\end{prop}
\begin{proof}
We shall proceed as in the proof of Proposition~\ref{prop:u-geq-B}. By inequality~\eqref{ghp.inq}, for $t\geq t_0 = 1$, we have $u(t,x)\leq \B_{M_2}(t+\tau_2,x)$. By integrating~\eqref{ghp.inq}, we find as well that $M\le M_2$. As previously, we can assume $M<M_2$. Let us define $\overline{\varepsilon}>0$ such that $(1+\overline{\varepsilon}) M^\beta = M_2^\beta$ and let $0<\varepsilon< \min\{\overline{\varepsilon},1\}$. In order to prove the claim we need to prove that for $|x|$ and $t$ large enough it holds
\begin{equation}
    \label{BB-upper-bound}
\frac{\B_{M_2}(t+\tau_2,x)}{{\B}_{\Mo}(t,x)}\leq 1+\ve\,.
\end{equation}
Let us notice that for $\eta(t):=\left(\frac{t}{t+\tau_2}\right)^\beta$, $s$, and $c$ as in~\eqref{eta-s-C}, we have
\[
\frac{\B_{M_2}(t+\tau_2,x)}{{\B}_{\Mo}(t,x)}=\eta^N( t) \left(\frac{1 +cs(t,x)}{\left(1+\overline{\varepsilon}\right)^\frac{p-2}{p-1} +c\eta^{p'}(t)s(t,x)}\right)^\frac{p-1}{2-p}\,,
\]
and inequality~\eqref{BB-upper-bound} is equivalent to
\begin{equation}\label{BB-upper-bound-equivalent}
s(t,x)\geq \frac{1-\left(1+\overline{\varepsilon}\right)^\frac{p-2}{p-1}(1+\ve)^\frac{2-p}{p-1}\eta^{-N\frac{2-p}{p-1}}(t)}{c\left(\eta ^{\frac{1}{\beta(p-1)}}(t)(1+\ve)^\frac{2-p}{p-1}-1\right)}=:\overline{s}(t, \ve)\,,
\end{equation}
provided $\eta ^{\frac{1}{\beta(p-1)}}(t)(1+\ve)^\frac{2-p}{p-1}>1$. We restrict our attention to $t>\overline{T}(\ve)$, where
\begin{equation}\label{T-overline}
\overline{T}(\ve):=\max\left\{\frac{\tau_2}{(1+\ve)^{2-p}-1},\, t_2,\, t_0\right\}\,
\end{equation}
and $t_2>0$ is such that $\eta(t_2)^\frac{1}{\beta(p-1)}=2(1+(1+\ve)^\frac{2-p}{p-1})^{-1}$. We observe that for $t> \frac{\tau}{(1+\ve)^{2-p}-1}$, we have $\eta ^{\frac{1}{\beta(p-1)}}(t)(1+\ve)^\frac{2-p}{p-1}>1$. Then~\eqref{BB-upper-bound} will follow from~\eqref{BB-upper-bound-equivalent} and $t>\overline T(\ve)$. Inequality~\eqref{BB-upper-bound-equivalent} holds true as long as $s(t,x)\geq \overline{\vr}^{p'}(\ve)$ (for $\overline{\vr}$ defined below). Indeed, since $\eta\leq 1$ is increasing, and, for  $t\geq t_2$, it holds that $\eta(t)\leq \eta(t_2)$, we find that
\begin{align}\label{rho-overline}
    \overline{s}(t, \varepsilon)&\leq \frac{1}{c\left(1-(1+\ve)^\frac{2-p}{p-1}\eta^\frac{1}{\beta(p-1)}(t_2)\right)}= \frac{1}{c\left(1-\frac{2}{1+(1+\ve)^\frac{2-p}{p-1}}\right)}=\frac{1+(1+\ve)^\frac{2-p}{p-1}}{c\left((1+\ve)^\frac{2-p}{p-1}-1\right)}=:\overline{\vr}^{p'}(\ve)\,.
\end{align}
\end{proof}

\subsection{Convergence in inner cylinders}

\begin{prop}\label{prop:conv-inside.2}
Let $N\ge1$, $p_c<p<2$, $0\le \uo \in \LL^1(\rn)\cap\X$ and $M=\int_{\rn}\uo(x)\dx>0$. Assume $u$ is a weak solution to problem~\eqref{CPLE} with initial datum $\uo$.  Suppose that for some $\wt T>0$ and $\wt K>0$ inequality~\eqref{given.convergence} holds. Then there exist $\vr_0=\vr_0(p,N,\Mo)>0$, $\overline{K}_\star=\overline{K}_\star(p,N, \Mo, \wt K)>0$ such that for any $\vr\geq \vr_0$ we have
\begin{equation}\label{inq.conv-inside}
\left|\frac{u(t,x)-\B_{\Mo}(t,x)}{\B_{\Mo}(t,x)}\right| \le\, \overline{K}_\star\,\vr^\frac{p}{2-p}\,  t^{-\frac{\nu}{N+1}}\qquad\forall\,|x|\leq\, \vr\, t^{\beta}\quad\forall\, t\ge \wt T\,,
\end{equation}
where $\nu$ and $\wt T$ are as in~\eqref{given.convergence}.
\end{prop}
\begin{proof} For any $t>0$ and $|x|\leq \vr\, t^\beta$, the  relative error satisfies the following inequality
\[
\left|\frac{u(t,x)-\B_{\Mo}(t,x)}{\B_{\Mo}(t,x)}\right| \le \|u(t,\cdot)-\B_{\Mo}(t,\cdot)\|_{\LL^\infty(\rn)} \sup_{\{|x|\le \vr\, t^\beta\}}\, \frac{1}{\B_{\Mo}(t,x)}\,.
\]
Since for any $t>0$, the function $|x|\mapsto\B_{\Mo}(t,|x|)$ is decreasing, we find that the supremum in the above inequality is attained at $|x|=\vr\,t^{\beta}$. Using the expression of the Barenblatt profile~\eqref{beta} and $b_1,b_2$ from \eqref{b1} and \eqref{b2}, respectively, a simple computation shows that, for any $\vr\ge \vr_0:=b_1^\frac{p-1}{p}/\Mo^{\beta\,(2-p)}$, we have that
\[
\sup_{\{|x|\le \vr\, t^\beta\}}\,\frac{1}{ \B_{\Mo}(t,x)} \le \vr^\frac{p}{2-p}\, t^{N\beta}\, \left(1+b_2\right)^\frac{p-1}{2-p}= C(p,N)\,\vr^\frac{p}{2-p}\, t^{N\beta}\,.
\]
Combining the two estimates we find
\[
\left|\frac{u(t,x)-\B_{\Mo}(t,x)}{\B_{\Mo}(t,x)}\right| \le C\,\vr^\frac{p}{2-p}\,t^{N\beta}\,\|u(t,\cdot)-\B_{\Mo}(t,\cdot)\|_{\LL^\infty(\rn)}\,.
\]
Before continuing, let us recall the Gagliardo--Nirenberg inequality
\[
\|f\|_{\LL^\infty(\rn)} \le C_{N}\, \|\nabla f\|_{\LL^\infty(\rn)}^\frac{N}{N+1}\,\| f\|_{\LL^1(\rn)}^\frac{1}{N+1}\,,
\]
holding for any $f$ regular enough for which all the involved quantities are finite, see~\cite[Lemma 3.5] {Bonforte2020a}. Combining the above inequality with~\eqref{given.convergence}, the time decay of gradient of solutions~\eqref{time.grandient.decay}, and the triangle inequality, we get that for any $t\ge \w T$:
\begin{equation*}\begin{split}
    t^{N\beta}\,\|u(t,\cdot)-\B_{\Mo}(t,\cdot)\|_{\LL^\infty(\rn)} &\le C_N\, t^{N\beta}\,  \|\nabla u(t, \cdot) - \nabla \B_{\Mo}(t,\cdot)\|_{\LL^\infty(\rn)} ^\frac{N}{N+1}\, \|u(t,\cdot)-\B_{\Mo}(t,\cdot)\|_{\LL^1(\rn)} ^\frac{1}{N+1}\\
    &\le C_N\,t^{N\beta}\,c_1^\frac{N}{N-1}\left(2\,M t^{-(N+1)\beta}\right)^\frac{N}{N+1} \left(\wt{K}\,t^{-\nu}\right)^\frac{1}{N+1}\\
    &\le C_\star\, t^{N\beta - \frac{N}{N+1} (N+1)\beta - \frac{\nu}{N+1} } \le C_\star\, t^{-\frac{\nu}{N+1}}\,,
\end{split}\end{equation*}
where $C_\star>0$ is a constant depending on $N, p$, $\Mo$ and $\wt K$. Combining all the above estimates, we can pick $ \overline{K}_\star:=C(p,N)\,C_\star$.
\end{proof}

\subsection{Proof of convergence in relative error under a priori convergence in Lebesgue's space}
After establishing the convergence inside and outside the cylinders in the last two sections, the only remained task is to link them.
\begin{proof}[Proof of Theorem~\ref{theo:conv-L1-RE}]
    From Propositions~\ref{prop:u-geq-B}, \ref{prop:u-leq-B} and~\ref{prop:conv-inside.2}  we know that for fixed $\ve\in(0,1/2)$, there exist
    \[
    T(\varepsilon)=\max\{\underline{T}(\ve),\overline{T}(\ve), \wt T\}>0\qquad\text{and}\qquad\vr(\varepsilon)=\max\{\underline{\vr}(\ve),\overline{\vr}(\ve)\}>0\,,
    \]
where $\underline{T}, \overline{T}$ are defined in~\eqref{under-T} and~\eqref{T-overline}, respectively, and $\underline{\vr},\overline{\vr}$ are defined in~\eqref{under-rho} and \eqref{rho-overline}, respectively, such that
\begin{equation*}
    \left|\frac{u(t,x)-\B_{\Mo}(t,x)}{\B_{\Mo}(t,x)}\right| \le \varepsilon \qquad\forall |x| \geq \vr(\varepsilon)\, t^{\beta}\quad \forall t\ge T(\varepsilon)\,.
\end{equation*}
 In the same way, using~\eqref{inq.conv-inside}, we obtain that for $t\ge \wt T$ it holds
\[
\left|\frac{u(t,x)-\B_{\Mo}(t,x)}{\B_{\Mo}(t,x)}\right| \le\,\overline{K}_\star\,\vr(\varepsilon)^\frac{p}{2-p}\,  t^{-\frac{\nu}{N+1}}\, \quad \quad\forall |x|\le 2\,\vr(\varepsilon)\, t^{\beta}\,.
\]
By a simple computation one finds that there exists a constant $\kappa_0(p)>0$, such that $\vr(\varepsilon)\le \kappa_0(p)\,\varepsilon^\frac{1-p}{p}$. Therefore, we have that
\[
\overline{K}_\star\,\vr(\varepsilon)^\frac{p}{2-p}\,  t^{-\frac{\nu}{N+1}}\le \varepsilon\qquad\mbox{for}\quad t\ge C_p\, \varepsilon^{-\frac{1}{2-p}\frac{N+1}{\nu}}\,,
\]
where the constant $C_p=C_p(N, p, \overline{K}_\star)>0$ is independent of $\varepsilon$. Therefore
\begin{equation}\label{relative.error.epsilon}
\left\|\frac{u(t,\cdot)-\B_{\Mo}(t,\cdot)}{\B_{\Mo}(t,\cdot)}\right\|_{\LL^\infty(\rn)}  \le \varepsilon\qquad\mbox{for any}\qquad t\ge \max\left\{C_p\, \varepsilon^{-\frac{1}{2-p}\frac{N+1}{\nu}}, T(\ve)\right\}\,.
\end{equation}
From a careful analysis of the proofs of Propositions~\ref{prop:u-geq-B} and~\ref{prop:u-leq-B}, we find that there exist $\kappa_1,\kappa_2$, independent of $\ve$, such that $\kappa_1\, \ve^{-1}\le T(\ve) \le \kappa_2\,\ve^{-1}$.  Since $(N+1)(2-p)^{-1}\nu^{-1}>1$, we find that the left inequality in~\eqref{relative.error.epsilon} holds for any $ t\ge C_p\, \varepsilon^{-\frac{N+1}{(2-p)\nu}}$, when $\varepsilon$ is small enough. Let us take a positive integer $m$ such that $\varepsilon\in \left[2^{-(m+1)}, 2^{-m}\right]$, then for any $s=\frac {t}{C_p}\in\left[2^\frac{m(N+1)}{(2-p)\nu},2^\frac{(m+1)(N+1)}{(2-p)\nu}\right]$ we have that
\[
\left\|\frac{u(t,\cdot)}{\B_{\Mo}(t,\cdot)}-1\right\|_{\LL^\infty(\rn)}\le \varepsilon \le 2^{-m} \le 2\, 2^{-(m+1)} \le 2\, C_p^{\frac{(2-p)\nu}{N+1}}\, t^{-\frac{(2-p)\nu}{N+1}}\,.
\]
The above computation holds for any $\ve\in(0,1/2)$, so that we conclude that inequality~\eqref{convergence.relative.error} holds for any $t\ge T_\star:=2^\frac{N+1}{(2-p)\nu} C_p$ and $K_\star:=2\,C_p^{\frac{(2-p)\nu}{N+1}}$.  The proof is complete.
\end{proof}

 \section{Convergence in relative error with rates}\label{sec:conv-rel-er}
The main goal of this section is to prove Theorem~\ref{theo:RECR}. Our strategy starts with providing a convergence rate in a weaker norm (i.e., a Lebesgue norm) and then applying Theorem~\ref{theo:conv-L1-RE}. The \emph{intermediate asymptotics}, i.e. the behaviour of solutions to~\eqref{CPLE} for large $t$, is much better understood when the equation is written in a different set of variables. Hence, we shall study the behaviour of solutions to the rescaled problem~\eqref{eq:v}. Our main tool will be the entropy functional introduced in~\eqref{cal-E}.

\subsection{Relation between the relative entropy and the Fisher information along the flow}\label{ssub:entropy-fisher}

In this section we shall pave the way for the use of the entropy method. The first step is to establish a relation between the entropy functional $\cE$ defined in \eqref{cal-E} and the relative Fisher information $\cI$ given by~\eqref{cal-I-def}. But prior to that we should prove that both quantities are well defined under our assumptions. In order to do so, we introduce a~few sufficient conditions for a given function $v:\rn\rightarrow[0, \infty)$ to have finite  entropy and Fisher information.
\begin{description}[style=multiline]
\item[\namedlabel{A0}{(A0)}]
    There exist $C_1, C_2>0$ and $D>0$  such that
\begin{equation*}
   C_1\le \frac{v(y)}{\VD(y)} \le C_2\quad \forall y\in \rn\,;
\end{equation*}
\item[\namedlabel{A1}{(A1)}]
     There exist  $\ve\in(0,1)$ and $D>0$   such that
\begin{equation*}
    1-\varepsilon\le \frac{v(y)}{\VD(y)} \le 1+\varepsilon\quad \forall y\in \rn\,;
\end{equation*}

    \item[\namedlabel{A2}{(A2)}]  There exist  $D_1,D_2>0$  such that
\[
\V_{D_1}(y)\le v(y)\le \V_{D_2}(y)\quad \forall y\in \rn\,.
\]
\end{description}
If $v$ is a solution to~\eqref{eq:v}, condition~\ref{A0} just rewriting  the global Harnack principle~\eqref{ghp.inq} in the new variables~\eqref{change.variables}. Condition~\ref{A1} is a way of quantifying the  closeness in relative error~\eqref{UCRE/Xp}. Lastly, condition~\ref{A2} is nothing but assumption~\eqref{stronger.GHP} in selfsimilar variables defined in~\eqref{change.variables}. We also notice that~\ref{A2} always implies~\ref{A1}; it implies~\ref{A2} when $D_1$ and $D_2$ are sufficiently close to $D$. \newline

Here we collect some observations that explain why the above conditions are expected and useful in our reasoning.
\begin{rem}\label{rmq} \rm

 We suppose that $p\in(p_c,2)$ and $0\le \vo \in \LL^1(\rn)$ is an initial datum for $v$ being a weak solution to~\eqref{eq:v}.
\begin{enumerate}[{\it (i)}]
    \item If $0\le \vo\in\X\setminus\{0\}$, then for any $\tau_0>0$ there exist $C_1(\tau_0), C_2(\tau_0)>0$ for which $v$ satisfies~\ref{A0}.
    \item If $0\le \vo\in\X\setminus\{0\}$,  then for any $\ve\in(0,1)$ there exists $\tau_\ve>0$ such that the solution $v$ satisfies~\ref{A1} (with the same $\ve$ and for some $D>0$) for all $\tau>\tau_\ve$.
    \item If $0\le \vo\in \LL^1_{{\rm loc}}(\rn)$ satisfies~\ref{A1} for some $\ve\in(0,1)$ and some $D>0$, then there exists $\tau_\ve>0$ such that the solution $v$ satisfies~\ref{A1} (with the same $\ve$ and $D$) for all $\tau>\tau_\ve$.
    \item If $\vo$ satisfies~\ref{A2} for some $D_1, D_2>0$, then by the comparison principle $v(\tau, \cdot)$ satisfies~\ref{A2} with $D_1$ and $D_2$ for all $\tau>0$.
   \item Let $0\le \vo\in\X\setminus\{0\}$ and $D>0$ be such that $\int_{\rn} \vo \dy = \int_{\rn} \V_D\dy$. If $p_c<p\le p_M$, we additionally suppose that $\vo$ satisfies~\ref{A2} for some $D_1,D_2>0$. Then,  by applying the change of variable~\eqref{change.variables} to solutions to~\eqref{CPLE} and~\cite[Theorem 1.1]{Bonforte2020b}, we deduce that
\[
\left\|\frac{v(\tau, \cdot)-\V_D}{\V_D}\right\|_{\LL^\infty(\rn)}\longrightarrow 0\quad\mbox{as}\quad \tau\rightarrow\infty\,.
\]
As a consequence we also have that $\cE[v(\tau,\cdot)|\V_D]\rightarrow0$ as $\tau\rightarrow\infty$.
\end{enumerate}
\end{rem}

Let us now provide some sufficient conditions for the entropy functional to be finite. It will be clear then that such conditions are fulfilled by solutions to~\eqref{eq:v} and, consequently, the entropy functional is well-defined along the flow. Prior to that, let us recall the exponent $p_M=(3(N+1)+\sqrt{(N+1)^2+8})/(2(N+2))$ defined in~\eqref{p_M}. For $p_M<p<2$, solutions to~\eqref{eq:v} (under the assumption $\vo\in\X$) have finite $|y|^\frac{p}{p-1}$-moments.
This is a necessary condition for the entropy functional $\cE$ to be well-defined. Contrary to this case, in the range $p_c<p\leq p_M$, solutions to~\eqref{eq:v} do not have a finite $|y|^\frac{p}{p-1}$-moment anymore and we need to invoke a stronger assumption to make the entropy functional finite along the flow.

\begin{lem}\label{lem:ent-lin-ent}
 Let $N\ge1$, $p_c<p<2$, $0\leq v\in\LL^1(\rn)$, and $D>0$ be such that $\int_{\rn}\V_D(y)\dy=\int_{\rn}v (y)\dy$. Suppose $v$ satisfies~\ref{A0} for some $C_1,C_2>0$. In the case $p_c<p\leq p_M$ we additionally assume that $v$ satisfies~\ref{A2} for some $D_1,D_2>0$. Then
 \[
 \cE[v|\VD]<\infty\,.
 \]
\end{lem}
Before the proof, let us emphasize that, at least when $p$ is close to $2$, condition~\ref{A0} is not necessary to infer a finite relative entropy. We refer more to~\cite{DelPino2003,Agueh2008,Agueh2009} for further discussion. We have chosen to use condition~\ref{A0} since it is quite practical in our setting and solutions to~\eqref{eq:v} will automatically fulfil it for any $\tau>0$.

\begin{proof}
Let us consider the case $p_M<p<2$. We shall prove that $y\mapsto \left[v^\gamma(y)-\V_D^\gamma(y)\right]-\gamma \V_D^{\gamma-1}(y)[v(y)-\VD(y)] \in \LL^1(\rn)$. We start by proving that $\V_D^{\gamma-1}[v-\VD] \in \LL^1(\rn)$. Notice that $\V_D^{\gamma-1}(y)=D+\frac{p-2}{p}|y|^\frac{p}{p-1}$, so that it is sufficient to show that, under assumption~\ref{A0}, both $v$ and $|y|^\frac{p}{p-1}v$ are integrable. From~\ref{A0}, we know that there exists a constant $C_2$ such that $v(y) \le C_2\, \VD(y)$ for all $y\in\rn$. Since $p\in(p_M,2)$, i.e. $\VD$ has a finite weighted $|y|^{\frac{p}{p-1}}$-moment, the previous inequality shows that both $\int_{\rn} v(y)  \dy $ and $\int_{\rn} |y|^{\frac{p}{p-1}} v(y)  \dy$ are finite. We deduce than that $\V_D^{\gamma-1}(y)[v(y)-\VD(y)] \in \LL^1(\rn)$.  Consequently, to conclude we only need to prove that $v^\gamma, \V^\gamma_D$ are integrable. Since both $v$ and $\VD$ are integrable and have the   $|y|^{p'}$-moments finite, we are allowed to conclude by using the Carlson--Levin inequality (see~\cite[Lemma 5]{Carrillo_2019} and references therein):
\[
\mathcal{C}_{s, q, N}\,\left(\int_{\rn} |f(y)|^q\dy\right)^\frac{1}{q}\le \left(\int_{\rn} |y|^s |f(y)|\dy\right)^\frac{N(1-q)}{q\,s}\, \left(\int_{\rn}|f(y)|\dy\right)^{1-\frac{N(1-q)}{q\,s}}\,,
\]
which holds for any $s>0$ and $q$ such that $N/(N+s)<q<1$. We choose $s:=p'$ and $q:=\gamma$. Note that they satisfy the assumption of the above Carlson--Levin inequality, because $\gamma=\frac{2p-3}{p-1}$ and $p_M<p<2$ (cf.~\eqref{p_M}).

Let us consider now the case $p_c<p\leq p_M$. Fix $y\in\rn$ and consider identity~\eqref{MVT}, with $t=v(y)$, $s=\VD$ and $v\le\xi\le \VD$. We find from that inequality, by using assumption~\ref{A0}, for all $y\in\rn$,
\begin{equation}\label{eq1}
\frac{\V_D^{\gamma-2}(y)}{2C_2^{2-\gamma}}\,(v(y)-\VD(y))^2\le \frac{v^\gamma(y)-\V_D^\gamma(y)-\V_D^{\gamma-1}(y)\left(v(y)-\VD(y)\right)}{\gamma(\gamma-1)} \le \frac{\V_D^{\gamma-2}(y)}{2C_1^{2-\gamma}}\,(v(y)-\VD(y))^2\,.
\end{equation}
Under the stronger assumption~\ref{A2}, by Lemma~\ref{decay.lemma}, we also have that $|v(y)-\VD(y)|\le C\,|y|^{-\frac{p}{(2-p)(p-1)}}$ for a constant $C>0$ and $|y|\ge 1$. This decay is enough to prove that the first and last terms in inequality~\eqref{eq1} are integrable and so $\cE[v|\VD]$ is finite.
\end{proof}

  In order to show  that the Fisher information is the derivative in time of the entropy functional, we make use of the space decay rate of the gradient of solutions  to rescaled problem~\eqref{eq:v}. To get it we adapt inequality~\eqref{space.gradient.decay} for \eqref{CPLE}.

\begin{lem} Let $N\ge 1$, $p_c<p<2$, and $v$ be a weak solution to~\eqref{eq:v} with $0\le \vo \in \LL^1(\RN)\cap\X(\rn)$. Then there exists a constant $c_3=c_3(\tau,p,N)>0$ such that
\begin{equation}\label{grad-v-est}
|\nabla_y v(\tau, y) | \le c_3\, {\max\left\{1, \|\vo\|_{\LL^1(\rn)}^{2\beta}+\|\vo\|_{\X}^{2\beta}\right\}} |y|^{-\frac{2}{2-p}}\quad\mbox{for any}\quad |y|\ge1\ \text{ and }\ \tau>0 \,.
\end{equation}
Constant $c_3$ can be chosen in such a way that for all $\tau>0$ large enough it holds $c_3(\tau,p,N)\leq c(p,N)$ for some $c(p,N)$.
\end{lem}

\begin{proof}
We recall that, the problems~\eqref{CPLE} and~\eqref{eq:v} are related through the change of variables~\eqref{change.variables}. Therefore, the estimate~\eqref{space.gradient.decay} became
\[
\begin{split}
|\nabla_y v(\tau, y) | = e^{(N+1)\tau} |\nabla_x u(t,y e^\tau)| & \le c e^{(N+1)\tau} \frac{\left(\|\uo\|_{\LL^1(\rn)}^{2\beta}+\|\uo\|_{\X}^{2\beta}+e^{\tau\frac{2}{2-p}}\right)}{e^{\tau(N+1)}(1+e^\tau|y|)^\frac{2}{2-p}}\le C \frac{\left(\|\uo\|_{\LL^1(\rn)}^{2\beta}+\|\uo\|_{\X}^{2\beta}+e^{\tau\frac{2}{2-p}}\right)}{1+e^{\frac{2\,\tau}{2-p}}|y|^\frac{2}{2-p}}\,.
\end{split}
\]
Notice that we change a little bit the estimate from~\eqref{space.gradient.decay} to be valid for $t>1$, eventually with a~different constant. The constants above are independent of the initial datum. Having the above estimate,~\eqref{grad-v-est} follows from a direct computation.\end{proof}

Let us now focus on proving rigorously that the Fisher information is the derivative in time of the entropy functional along the flow defined by~\eqref{eq:v}. We concentrate now on the case $p_M<p<2$.

\begin{lem}\label{lem:ent/ent-prod}
 Let $N\ge1$, $p_M<p<2$, $0\leq \vo\in\LL^1(\rn)$, and $D>0$ be such that $\int_{\rn}\VD(y)\dy=\int_{\rn}\vo (y)\dy$. Let $v$ be a solution to~\eqref{eq:v} with $\vo$ as its initial datum. If $\vo$ satisfy~\ref{A0} holds, then for any $\tau_0>0$ we have that $\tau\mapsto\cI[v(\tau)|\VD] \in \LL^\infty(\tau_0, \infty)\cap \LL^1(\tau_0, \infty)$  and
\begin{equation}\label{entropy.production}
\tfrac{\d}{\d\tau}\cE[v(\tau)|\VD]= - \cI[v(\tau)|\VD]\quad\mbox{for almost every}\,\,\,\tau>\tau_0\,.
\end{equation}
\end{lem}
We notice that it is not necessary that $\vo$ satisfies~\ref{A0}, we could have just asked that $v(\tau)$ satisfies either~\ref{A0} or~\ref{A1} after some time $\tau_1>0$. Indeed, this condition will be satisfied by solutions to~\eqref{eq:v} once the initial datum is assumed to be in $\X$.

\begin{proof}
The formal proof goes through equation~\eqref{eq:v} and integration by parts, namely:
\begin{flalign*}
-\frac{\d}{\d\tau}(\cE[v(\tau)|\VD])&=\frac{1}{\gamma-1}\int_\rn \left(v^{\gamma-1}-\V_D^{\gamma-1}\right)\dv_y \Big(v(\tau) \cdot \ba[ v(\tau) ]+v(\tau)\cdot y \Big) \dy\\
&=\frac{1}{|\gamma-1|^p}\int_\rn v(\tau) \left(\nabla v^{\gamma-1}-\nabla \V_D^{\gamma-1}\right)\cdot\Big( \bb[v^{\gamma-1}(\tau)]-\bb[\V_D^{\gamma-1}]\Big) \dy\,,
\end{flalign*}
where $\bb[\phi]:=|\nabla \phi|^{p-2}\nabla \phi$ and $\ba=\bb[\phi]$. Let us justify it rigorously. For a smooth cut-off function $\phi_R$ such that $\phi_R=1$ in $B_R$ and $\phi_R=0$ outside $B_{2R}$, we define
\begin{flalign*}
\cE_{\phi_R}[v(\tau)|\VD]&:=\frac{1}{\gamma(\gamma-1)}\int_{\rn}\phi_R(y)\left\{ \left[v^\gamma(\tau,y)-\V_D^\gamma(y)\right]-\gamma \V_D^{\gamma-1}(y)[v(\tau,y)-\VD(y)]\right\} \dy\,,\\
\cI_{\phi_R}[v(\tau)|\VD]&:=\frac{1}{|\gamma-1|^p}\int_{\rn}\phi_R(y)\ v(\tau,y)\big(\nabla v^{\gamma-1}(\tau,y)-\nabla  \V_D^{\gamma-1}(y) \big) \cdot\Big(\bb[ v^{\gamma-1}(\tau,y)]-\bb[\V_D^{\gamma-1}(y)]\Big) \dy\,,\\
{\cal R}_{\phi_R}[v(\tau)|\VD]&:=-\frac{1}{\gamma-1}\int_\rn\nabla\phi_R(y)\ (v^{\gamma-1}(y)-\V_D^{\gamma-1}(y))(|\nabla v(y)|^{p-2}\nabla v(y)+y v(y))\dy.
\end{flalign*}
It will be clear in a few lines that the above quantities are well defined. In what follows, we would like to test equation~\eqref{eq:v} against $\gamma(v^{\gamma-1}-\V_D^{\gamma-1})\phi_R$. This is an admissible test function. Indeed, the assumption~\ref{A0} on $\vo$ implies (by the maximum principle) that $v(\tau)$ satisfies~\ref{A0} for all $\tau>0$. Moreover, we know that for a fixed $R>0$ function  $y\mapsto v^{\gamma-1}\phi_R$ is  $C^1$. Actually more is known: the function $y\mapsto v(\tau, y)$ is $C^{1,\alpha}$ locally in space for some $\alpha=\alpha(N, p)$, i.e., $\|\nabla v\|_{C^{0,\alpha}(B_R)}$ is finite for any $R>0$, see~\cite[Theorem III.8.1]{DiBenedetto1990}. Therefore, by a  standard approximation procedure, we can test~\eqref{eq:v} against $\gamma(v^{\gamma-1}-\V_D^{\gamma-1})\phi_R$ to get
\begin{align*}L:=\int_s^t\int_\rn\gamma \phi_R \left(v^{\gamma-1}-\V_D^{\gamma-1}\right)\partial_\tau v\dy\dtau=-\int_s^t \int_\rn\gamma\phi_R \left(v^{\gamma-1}-\V_D^{\gamma-1}\right)\dv_y \Big(v(\tau) \cdot \ba[ v(\tau) ]+v(\tau)\, y \Big) \dy\dtau=:K,
\end{align*}
where
\begin{align*}
    L&=\int_\rn \gamma v(t)(v^{\gamma-1}(t)-\V_D^{\gamma-1})\phi_R\dy-\gamma\int_\rn v(s)(v^{\gamma-1}(s) -\V_D^{\gamma-1})\phi_R\dy-\int_s^t\gamma\int_\rn v(\gamma-1)v^{\gamma-2}\partial_\tau v\,\phi_R\dy\dtau\\
    &=:L_1+L_2+L_3\,.
\end{align*}
By Fubini's theorem we infer that
\begin{align*}
    L_3=(\gamma-1)\int_s^t \int_\rn\partial_\tau(v^{\gamma})\,\phi_R\dy\dtau=(\gamma-1)\left[\int_\rn v^\gamma(t)\phi_R\dy-\int_\rn v^\gamma(s)\phi_R\dy\right]\,.
\end{align*}
Let us notice that
\begin{align*}
    L&=\int_\rn v^{\gamma}(t)\phi_R\dy-\int_\rn \gamma v(t)\V_D^{\gamma-1}\phi_R\dy-\int_\rn v^{\gamma}(s) \phi_R\dy+\int_\rn \gamma v(s)\V_D^{\gamma-1}\phi_R\dy\\
    &\quad-(\gamma-1)\left[\int_\rn \V_D^\gamma\phi_R\dy-\int_\rn \V_D^\gamma\phi_R\dy\right]\,\\
    &=\int_{\rn}\phi_R\ \left[v^\gamma(t,y)-\V_D^\gamma(y)\right]-\gamma \V_D^{\gamma-1}(y)[v(t,y)-\VD(y)] \dy\\
    &\quad -\int_{\rn}\phi_R\ \left[v^\gamma(s,y)-\V_D^\gamma(y)\right]-\gamma \V_D^{\gamma-1}(y)[v(s,y)-\VD(y)] \dy\\
    &=\gamma(\gamma-1)\left(\cE_{\phi_R}[v(t)|\VD]-
\cE_{\phi_R}[v(s)|\VD]\right)\,.
\end{align*}
On the other hand,
\begin{align*}
    K&=\gamma \int_s^t\int_\rn \phi_R\ v(\tau) \left(\nabla v^{\gamma-1}-\nabla \V_D^{\gamma-1}\right)\cdot\Big( \bb[v^{\gamma-1}(\tau)]-\bb[\V_D^{\gamma-1}]\Big) \dy\dtau\\
    &\quad+\gamma \int_s^t\int_\rn \nabla\phi_R\ (v^{\gamma-1}-\V_D^{\gamma-1})(|\nabla v|^{p-2}\nabla v+y v)\dy\dtau\\
    &= -\gamma(\gamma-1)\int_s^t\cI_{\phi_R}[v(\tau)|\VD]\dtau+\gamma(\gamma-1)\int_s^t\mathcal{R}_{\phi_R}[v(\tau)|\VD]\dtau\,.
\end{align*}
We find therefore, by posing $s=\tau_0$ and $t=\tau_0+h$, for $h>0$, that
\begin{align}\label{E-difference}
\cE_{\phi_R}[v(\tau_0+h)|\VD]-
\cE_{\phi_R}[v(\tau_0)|\VD]=-\int_{\tau_0}^{\tau_0+h} \cI_{\phi_R}[v(\tau)|\VD]\dtau-\int_{\tau_0}^{\tau_0+h} {\cal R}_{\phi_R}[v(\tau)|\VD]\dtau\,.
\end{align}

We need to prove that, for almost every $\tau>0$, it holds
\begin{align}
    \cE_{\phi_R}[v(\tau)|\VD]&\xrightarrow[R\to\infty]{}\cE[v(\tau)|\VD]\,,\label{EphiR-conv}\\
\int_{\tau_0}^{\tau_0+h} \cI_{\phi_R}[v(\tau)|\VD]\dtau&\xrightarrow[R\to\infty]{} \int_{\tau_0}^{\tau_0+h} \cI[v(\tau)|\VD]\dtau\,,\label{IphiR-conv}\\
\int_{\tau_0}^{\tau_0+h}\big|\mathcal{R}_{\phi_R}[v(\tau)|\VD]\big|&\xrightarrow[R\to\infty]{}0\,.\label{RphiR-conv}
\end{align}
This would allow to pass to the limit in~\eqref{E-difference}, from which we will get that
\begin{equation}\label{needed.step}
\cE[v(\tau_0+h)|\VD]-
\cE[v(\tau_0)|\VD]=-\int_{\tau_0}^{\tau_0+h} \cI[v(\tau)|\VD]\dtau\,.
\end{equation}
Then, by using the  Lebesgue Differentiation Theorem,  we will conclude with~\eqref{entropy.production} and the proof would be complete.

Let us first deal with the entropy $\cE$. By identity~\eqref{MVT} (with $t=v(\tau,y)$ and $s=\VD(y)$) we know that the integrand of $\cE[v(\tau, \VD)]$ is positive, and, by the Monotone Convergence Theorem, we can say that~\eqref{EphiR-conv} holds. Let us now focus on the relative Fisher information. In order to justify~\eqref{IphiR-conv} and to prove that $\tau\mapsto\cI[v(\tau)|\VD] \in \LL^\infty(\tau_0, \infty)$ we will make use of the Dominated Convergence Theorem. Let us note that that the integrand of $\cI_{\phi_R}[v(\tau)|\VD]$ can be estimated pointwise by using Young's inequality and assumption~\ref{A0} (on the solution $v(\tau)$) as follows
\begin{align*}
    I_\tau(y)&:=v(\tau,y)\big(\nabla v^{\gamma-1}(\tau,y)-\nabla  \V_D^{\gamma-1}(y) \big) \cdot\Big(\bb[ v^{\gamma-1}(\tau,y)]-\bb[\V_D^{\gamma-1}(y)]\Big)\\
    &\leq c(p,\gamma)| v(\tau,y)| \,\big(|\nabla  v^{\gamma-1}(\tau,y)|+|\nabla \V_D^{\gamma-1}(y)|\big)^{p} \\
    &\leq c(p,\gamma)| v(\tau,y)|  \,\big(|\nabla  v^{\gamma-1}(\tau,y)|^p+|\nabla \V_D^{\gamma-1}(y)|^p\big) \\
    &\leq c(p,\gamma)| v(\tau,y)|  \,\big( v^{p(\gamma-2)}(\tau,y)|\nabla  v(\tau,y)|^p+ v^{p(\gamma-2)}_D(\tau,y)|\nabla \VD(y)|^p\big)\\ &\leq c(p,\gamma, C_2)| \VD(y)|^{-\frac{1}{p-1}}  \,\big( |\nabla  v(\tau,y)|^p+ |\nabla \VD(y)|^p\big)\,,
\end{align*}
Notice that in the above computation the value of the constant may change from line to line and in the last step we have used assumption~\ref{A0} so that $c$ depends on the value of the constant $C_2$ appearing in~\ref{A0}.  For $|y|\leq 1$ we note that the right-hand side of the above inequality is bounded by a constant. Indeed, the only element that cannot be explicitly computed is $|\nabla v(\tau, y)|$. However notice that, by applying inequality~\eqref{time.grandient.decay} (after the change of variables~\eqref{change.variables}), we get that $|\nabla v(\tau, \cdot)|\in\LL^\infty(\rn)$ for any $\tau>0$. On the other hand,  for $|y|\geq 1$, we can estimate $|\nabla v(\tau, y)|$ by inequality~\eqref{grad-v-est}, to obtain that:
\begin{align*}
    |I_\tau(y)|    &\leq c(\uo,p,\gamma)| \V_D(\tau,y)||y|^{\frac{p^2}{(2-p)(p-1)}}\,|y|^{-\frac{2p}{2-p}}= c(\uo,p,\gamma)| \V_D(\tau,y)||y|^{p'}\,,
\end{align*}
where the right-hand side is integrable outside a ball as long as $p>p_M$. Consequently, we have proven that $\tau\mapsto\cI[v(\tau)|\VD] \in \LL^\infty(\tau_0, \infty)$. At the same time, it is clear from the above that~\eqref{IphiR-conv} follows from the Dominated Convergence Theorem.

Lastly, let us now consider the error term $\mathcal{R}$. We notice that $\mathcal{R}_{\phi_R}[v(\tau)|\V_D]=\int_{\rn} F_R(y)\dy$ with $F_R(y)\rightarrow 0$ almost everywhere as $R\rightarrow \infty$, since the term $\nabla \phi_R(y)$ is supported only in $A_R:=B_{2R}\setminus B_{R}$. So, to prove~\eqref{RphiR-conv},  we only need to find an integrable function $G$ such that $|F_R(y)|\le G(y)$ (uniformly in $R$ and $y$) and then to invoke the Dominated Convergence Theorem. We already know that $F_R=0$ in $A_R^c$, so that we only need to estimate it in $A_R$. Since $\vo$ satisfies~\ref{A0}, so it does~$v(\tau, \cdot)$ for any $\tau>0$. Therefore, by applying condition~\ref{A0} with inequality~\eqref{grad-v-est}, we find that, for $|y|$ large enough
\[\begin{split}
|F_R(y)|& \le \tfrac{C_1}{R}\,|\V_D(y)|^{\gamma-1}\left(|\nabla v|^{p-1} + |y|\,|v
\V_D(y)|\right) \le \tfrac{C_2}{R}\, |v
\V_D(y)|^{\gamma-1}\,|y| \left({|y|^{-\frac{p}{2-p}}} +  |\V_D(y)| \right) \\
&\le C_3\, |\V_D(y)|^{\gamma-1} \left({|y|^{-\frac{p}{2-p}}} +  |\V_D(y)| \right) \le C_4\,|\V_D(y)|^{\gamma}\,,
\end{split}
\]
where in the third inequality we used the fact that $|y|/R\le 2$ on $A_R$, while in the fourth one simply applies the fact that $|y|^{-\frac{p}{2-p}}\le \kappa\, \V_D(y)$, for a $\kappa>0$ independent of $y$. We recall that in the proof of Lemma~\ref{lem:ent-lin-ent} we have already proven that the function $y\mapsto|v_D(y)|^\gamma$ is integrable whenever $p\in\left(p_M,2\right)$. Therefore~\eqref{RphiR-conv} holds by the Dominated Convergence Theorem. Consequently, identity~\eqref{needed.step} is proven.

It only remains to prove that $\tau\mapsto \cI[v(\tau)|\VD]$ is in $\LL^1(\tau_0, \infty)$ for any $\tau_0>0$. This is easily done by observing that,  we have $\cE[v(\tau+h)|\VD]\rightarrow0$ as $h\rightarrow\infty$. Therefore, we find that $\int_{\tau_0}^\infty \cI[v(\tau)|\VD]\dtau = \cE[v(\tau_0)|\VD]$. We can conclude by observing that $\cI[v(\tau)|\VD]\ge0$ and $\cE[v(\tau_0)|\VD]<\infty$.
\end{proof}

We shall now look at the relation between the relative entropy and the Fisher information when $p_c<p\leq p_M$. The main difficulty here is that the Fisher information might be an unbounded function of time. We can still establish that the Fisher information is the derivative of the entropy along the flow within this range, but in a weaker form and under stronger assumptions.

\begin{lem}\label{lem:ent/ent-prod-2}
 Let $N\ge1$, $p_c<p\le p_M$, $0\leq \vo\in\LL^1(\rn)\cap\X$, and $D>0$ be such that $\int_{\rn}\V_D(y)\dy=\int_{\rn}\vo (y)\dy$. Let $v$ be a solution to~\eqref{eq:v} with $\vo$ as its initial datum. If $\vo$ satisfies both~\ref{A1} and~\ref{A2}, then $\tau\mapsto\cI[v(\tau)|v_D] \in \LL^1(\tau_0, \infty)$ for any $\tau_0>0$ and
\begin{equation}\label{weak.entropy.production}
\cE[v(\tau_0)|v_D] = \int_{\tau_0}^\infty \cI[v(\tau)|v_D]\dtau\quad\forall\,\tau_0>0\,.
\end{equation}
\end{lem}
\begin{proof}
Let us first  notice that, under the current assumption we have that $\cE[v(\tau)|\VD]<\infty$ for any $\tau>0$. Indeed, if $\vo$ satisfies~\ref{A1} then $\vo\in\X$ and, by Remark~\ref{rmq}~{\it (i)}, for any $\tau>\tau_0>0$,  the solution $v$ satisfies~\ref{A0} for some constants $C_0, C_1>0$ which depend on $\tau_0$. Then, thanks to Lemma~\ref{lem:ent-lin-ent}, we know that $\cE[v(\tau)|\VD]<\infty$ for any $\tau>0$.

We start with proceeding along the lines of the proof of Lemma~\ref{lem:ent/ent-prod}. What we need to motivate differently is that~\eqref{needed.step} still holds in the current regime. The Monotone Convergence Theorem ensures that~\eqref{EphiR-conv} and~\eqref{IphiR-conv} hold true. However, it does not directly imply that $\int_{\tau_0}^{\tau_0+h} \cI[v(\tau)|\V_D]\dtau<\infty$ for any $h, \tau_0>0$, since~\eqref{RphiR-conv} is not known yet. With the same notation as in the proof of Lemma~\ref{lem:ent/ent-prod}, recall that $\mathcal{R}_{\phi_R}[v(\tau)|\V_D]=\int_{\rn} F_R(y)\dy$ with $F_R(y)=-\frac{1}{\gamma-1}\nabla \phi_R(y)\cdot (v^{\gamma-1}-\V_D^{\gamma-1})(|\nabla v|^{p-2}\nabla v+y\, v)$ supported in $A_R=B_{2R}\setminus B_R$. We restrict attention to $R>1$. By using assumption~\ref{A2} and the Mean Value Theorem we have that
\[
|v^{\gamma-1}(y)-\V_D^{\gamma-1}(y)|\le |\gamma-1|\xi^{\gamma-2}|v(y)-\V_D(y)|\,,
\]
where $\min\{\V_{D_1}(y),\V_D(y)\}\le \xi \le \max\{\V_{D_2}(y),\V_D(y)\}$.  Note that $|\xi|\ge c_0\,|y|^{-\frac{p}{2-p}}$, for a constant $c_0>0$ independent of $y$, so $|\xi|^{\gamma-2}\le c_1\,|y|^{\frac{p}{(2-p)(p-1)}}$. By Lemma~\ref{decay.lemma}, for $|y|>1$ we also have that $|v(y)-\V_D(y)|\le c_2\,|y|^{-\frac{p}{(2-p)(p-1)}}$ for a constant $c_2>0$ independent of $y$. This is enough to conclude that for $|y|>1$ it holds $|v^{\gamma-1}(y)-V_D^{\gamma-1}(y)|\le c_3$ for a constant $c_3>0$, independent of $y$. Since $\nabla \phi_R$ is supported in $A_R$ and  $|\nabla \phi_R \cdot y|\le c_4$ for a constant $c_4>0$ independent of $y$, we can write that for all $R>1$ and a  constant $C>0$ it holds
\[
\Big|{\cal R}_{\phi_R}[v(\tau)|\V_D]\Big|\le C \left(\frac{1}{R}\int_{A_R} |\nabla v(\tau, y)|^{p-1} \dy + \int_{A_R} v(\tau, y) \dy\right)\,.
\]
On the right-hand side above, the second integral converges to $0$, as $R\rightarrow\infty$, since $v(\tau, \cdot)$ is an integrable function and $|A_R|\subset\{|y|>R\}$. It only remains to justify the convergence of the first integral. By inequality~\eqref{grad-v-est}, we have that for large enough $|y|$ and large enough $\tau>0$ it holds $|\nabla v(\tau, y)|^{p-1}\le c_5\,|y|^{-\frac{2(p-1)}{2-p}}$ for a constant $c_5>0$ independent of $y$ (a similar estimate holds also for small $\tau$ but with a constant $c_5$ depending on $\tau$). Therefore, we are left with
\[
\frac{1}{R}\int_{A_R} |\nabla v(\tau, y)|^{p-1} \dy \le \frac{c_6}{R}\,\int_R^\infty r^{N-1-\frac{2(p-1)}{2-p}}\dr\le c_7\,R^{N-\frac{p}{2-p}}\,,
\]
for a constant $c_7>0$ independent of $R$. This justifies~\eqref{RphiR-conv} and implies that $\int_{\tau_0}^{\tau_0+h} \cI[v(\tau)|\V_D]\dtau<\infty$ for any $h, \tau_0>0$. Therefore,~\eqref{needed.step} holds true. The proof is complete.
\end{proof}

\subsection{Linearised entropy functional and Fisher information}

Recall that the entropy functional $\cE$ is given by~\eqref{cal-E}. We define the \emph{linearised relative entropy} by
\begin{equation*}
\linE[v]:=\frac{1}{2}\int_\rn |v-\VD|^2 \V_D^{\gamma-2} \dy\,.
\end{equation*}
This functional will play an important role in the most challenging regime $p_c<p\leq p_M$. Let us justify that this functional is well defined.

\begin{lem}\label{lem:ent-lin-ent-2}
Let $N\ge1$, $p_c<p<2$, $0\leq v\in\LL^1(\rn)$, and $D>0$ be such that $\int_{\rn}\VD(y)\dy=\int_{\rn}v (y)\dy$. Suppose $v$ satisfies~\ref{A1} for some $\ve\in\left(0,1\right)$. If $p_c<p\leq p_M$, we additionally assume that $v$ satisties~\ref{A2} for some $D_1,D_2>0$. Then $\linE[v]<\infty$ and
\begin{equation*}
    \left(1+\varepsilon\right)^{\gamma-2} \linE[v] \le \cE[v |\VD] \le \left(1-\varepsilon\right)^{\gamma-2} \linE[v]\,.
\end{equation*}
\end{lem}

\begin{proof}
We notice that under the current assumptions inequality~\eqref{eq1} holds trues. By Lemma~\ref{decay.lemma}, we also have that $|v(y)-\VD(y)|\le C\,|y|^{-\frac{p}{(2-p)(p-1)}}$ for a constant $C>0$ and $|y|\ge 1$. This decay is enough to prove that the first and last terms in inequality~\eqref{eq1} are integrable. Then, by integrating inequality~\eqref{eq1}, we deduce the claim.
\end{proof}

Let us now introduce two different quantities which, in the sequel, will play the rôle of a \emph{linearised} version of the Fisher information $\cI$ defined in~\eqref{cal-I-def}. For any function $v:\rn\rightarrow\R$ and any $\eta\geq0$, consider
\begin{equation}\label{def-I-ve-gamma}
\linI^{(\eta)}_\gamma[v]:=\frac{1}{|\gamma-1|^p} \int_\rn  \left|\nabla v^{\gamma-1}-\nabla \V_D^{\gamma-1}\right|^2\,\VD\,\left(\eta+|\nabla \V_D^{\gamma-1}|\right)^{p-2} \dy
\end{equation}
and
\begin{equation}\label{def-I-ve}
\linI^{(\eta)}[v]:= \frac{1}{|\gamma-1|^p}\int_\rn  \left|\nabla \left(\V_D^{\gamma-2}(v- \V_D)\right)\right|^2\,\VD\,\left(\eta+|\nabla \V_D^{\gamma-1}|\right)^{p-2} \dy\,.
\end{equation}
Our first result is a control, along the flow, of the Fisher information $\cI[v(\tau)|\VD]$ from below  by the first quantity defined above, at least when the considered solution $v$ is close to Barenblatt profile $\VD$ in the sense of assumptions~\ref{A1} and~\ref{A2}.

\begin{lem}\label{lem:I-est}
Let $N\ge1$, $p_c<p < 2$, $0\leq \vo\in\LL^1(\rn)$, and $D>0$ be such that $\int_{\rn}\V_D(y)\dy=\int_{\rn}\vo (y)\dy$. Let $v$ be a solution to~\eqref{eq:v} with $\vo$ as its initial datum and assume that  $\vo$ satisfies~\ref{A1} for some $\ve\in(0,1)$. If $p_c<p<p_M$, we additionally assume that $\vo$ satisfies~\ref{A2} for some $D_1,D_2>0$. Then there exists $\tau_\ve>0$ and $C_\ve=C_\varepsilon(p,N,\vo,\varepsilon, D)>0$ such that for all $\tau>\tau_\varepsilon$ it holds
\begin{flalign}\label{cal-I-est}
\cI[v(\tau)|\VD]\geq C_\varepsilon\, \linI_\gamma^{(\varepsilon)} [v(\tau)]\,.
\end{flalign}
 \end{lem}

\begin{proof}
We recall that, by Remark~\ref{rmq} that, there exists $\tau_\ve>0$ such that $v(\tau, \cdot)$ satisfies~\ref{A1} for all $\tau>\tau_\ve$. By applying Lemma~\ref{lem:p-male} and assumption~\ref{A1} we find  for all $\tau>\tau_\ve$ that
\begin{flalign}
\nonumber\cI[v(\tau)|\VD]&\geq\frac{\min\{1/2, p - 1\}}{|\gamma-1|^p}\int_{\rn}v(\tau,y)\big(|\nabla v^{\gamma-1}(\tau,y)|+|\nabla  \V_D^{\gamma-1}(y)|\big)^{p-2}\big|\nabla  v^{\gamma-1}(\tau,y)-\nabla \V_D^{\gamma-1}(y)\big|^2 \dy \\
\label{tech-est-on-cI}&\geq \frac{\min\{1/2, p - 1\}}{|\gamma-1|^p}(1-\varepsilon)\int_{\rn}\VD(y)\big(|\nabla v^{\gamma-1}(\tau,y)|+|\nabla  \V_D^{\gamma-1}(y)|\big)^{p-2}\big|\nabla  v^{\gamma-1}(\tau,y)-\nabla \V_D^{\gamma-1}(y)\big|^2 \dy \,,
\end{flalign}
where in the second inequality we use condition~\ref{A2} for all $\tau>\tau_\ve$. We will to show that
\begin{equation}\label{est}
|\nabla v^{\gamma-1}(y)|\le C(\vo,p,N,\varepsilon,D)\left(\varepsilon+|\nabla \V_D^{\gamma-1}(y)|\right)\,.
\end{equation}
Let us notice that $|\nabla \V_D^{\gamma-1}(y)|=(1-\gamma)|y|^{2-\gamma}\ge 0$ and $\nabla v^{\gamma-1}=(\gamma-1)\,v^{\gamma-2}\nabla v$. By the fact that $\|\nabla v\|_{\LL^\infty(\rn)}<\infty$ (resulting from~\eqref{time.grandient.decay}) and condition~\ref{A2} we have that, for $|y|\le 1$ and $\tau>\tau_\ve$ it holds
\begin{flalign*}
|\nabla v^{\gamma-1}(\tau, y)|&\le (1-\gamma)(1-\varepsilon)^{\gamma-2}\,\V_D^{\gamma-2}(y)\,\|\nabla v\|_{\LL^\infty(\rn)}\le \frac{(1-\gamma)}{(1-\ve)^{2-\gamma}} \sup_{|y|\le1}\V_D^{\gamma-2}\, \frac{\|\nabla v\|_{\LL^\infty(\rn)}}{\varepsilon}\,\left(\varepsilon+|\nabla \V_D^{\gamma-1}(y)|\right)\\
&\le C(\vo,p,N,\ve,D)\left(\varepsilon+|\nabla \V_D^{\gamma-1}(y)|\right)\,.
\end{flalign*}
 On the other hand,  for $|y|\ge1$ and $\tau>\tau_\ve$, by using~\eqref{grad-v-est} we infer that for a constant $C(\vo,p,N)>0$ it holds
\begin{flalign*}
|\nabla v^{\gamma-1}(\tau, y)|&\le (1-\gamma)\,(1-\varepsilon)^{\gamma-2}\,\V_D^{\gamma-2}(y)\,C(\vo)|y|^{-\frac{2}{2-p}} = (1-\varepsilon)^{\gamma-2}\,C(\vo,p,N)\,|\nabla \V_D^{\gamma-1}(y)|\,\V_D^{\gamma-2}(y)\,|y|^{-\frac{p}{(2-p)(p-1)}}\\
&\le C(\vo,p,N,\ve,D) \,\left(\varepsilon+|\nabla \V_D^{\gamma-1}(y)|\right)\,,
\end{flalign*}
where we used the fact that $\V_D^{\gamma-2}(y)\,|y|^{-\frac{p}{(2-p)(p-1)}}\le C(D, N, p)$ for $|y|\ge1$ where $C(D, N, p)>0$ depends only on $D$, $N$ and $p$. Combining the two cases $|y|\le1 $ and $|y|\ge1$ together, we find that \eqref{est} holds for all $y\in\rn$.
Then \eqref{est} together with~\eqref{tech-est-on-cI} implies~\eqref{cal-I-est} with the constant
$C_\varepsilon =2^{p-2}\,C(\vo,p,N,\varepsilon, D)^{p-2}\,c_2 (1-\varepsilon)$.
 \end{proof}

The previous lemma show a relation between the Fisher information and one of its linearised versions. In what follows, we state an inequality that holds among all the linearised quantities introduced in this section. The following lemma is originally contained in~\cite[Claim~1 of Proposition~4.2]{Agueh2009}. In that paper, the authors also use some quantities very similar to ours~$\linI^{(\varepsilon)}$ and~$\linI^{(\ve)}_\gamma$, however their definition is slightly different. Here is their result written in our notation. Notice slightly different constants and extended range of $p$ in comparison with~\cite{Agueh2009}.

\begin{lem}
\label{lem:ABC}   Let $N\ge1$, $p\in(1,2)$, $N<\frac{p}{(2-p)(p-1)}$, $0\leq v\in\LL^1_{{\rm loc}}(\rn)$ such that $\nabla v\in\LL^2_{{\rm loc}}(\rn)$, and let  $D>0$ be such that $v-\VD\in\LL^1(\rn)$ and $\int_{\rn}\left(\V_D(y)-v(y)\right)\dy=0$. Suppose $v$ satisfies~\ref{A1} for some $\ve\in\left(0,1\right)$. If $1<p\leq p_M$, we additionally assume that $v$ satisties~\ref{A2} for some $D_1,D_2>0$.  Then, for any $\eta>0$, we have that
\begin{equation}\label{inq.linearised.quantities}
\linI^{(\eta)}[v] \le \kappa_1(\ve)\,\linI^{(\eta)}_\gamma[v] + \kappa_2(\ve)\,\linE[v]
\end{equation}
where
\begin{equation}\label{kappas}
\kappa_1(\ve)=\frac{(1+\ve)^{2(2-\gamma)}}{(1-\gamma)^2}\qquad \mbox{ and } \qquad \kappa_2(\ve)=\frac{\mathcal{C}_{p,N}}{(1-\gamma)^{p-1}}\left(\frac{(1+\ve)^{2(2-\gamma)}}{(1-\ve)^{2(2-\gamma)}}-1\right)\,,
\end{equation}
where $\mathcal{C}_{p,N}>0$ is a constant depending on $N$ and $p$.
\end{lem}
In what follows we only sketch the main steps of the reasoning, which is based on the proof of~\cite[Proposition~4.2]{Agueh2008}.
\begin{proof}
We first explain the outline of the proof and then we justify the key technicality, which is an integration by parts. Let us introduce $h_k(s):=s^{k-1}-1$, for $k\in\{2,\gamma\}$, and $\mu(y):=(\eta+(1-\gamma)|y|^{2-\gamma})^{p-2}\,\V_D(y)$ so that $\d\mu(y):=\mu(y)\dy$. In order to unify the notation, let us define
\[
I^\eta_k [w]:=\frac{1}{|1-\gamma|^p}\int_{\rn}\left|\nabla(\V_D^{\gamma-1} h_k(w))\right|^2\,\d\mu(y)\,,\qquad \text{for }\, k\in\{2,\gamma\}\,.
\]
Notice that with this definition and upon setting $w:=v/\VD$, we have that $\linI^{(\eta)}[v]=I^\eta_2[w]$ and $\linI^{(\eta)}_\gamma[v]=I^\eta_\gamma[w]$, where $\linI^{(\eta)}[v]$ and $\linI^{(\eta)}_\gamma[v]$ are defined in~\eqref{def-I-ve} and in~\eqref{def-I-ve-gamma}, respectively. Recall that $\nabla \V_D^{\gamma-1}=(1-\gamma)|y|^{1-\gamma}y$.  Then, by computing the gradient $\nabla(\V_D^{\gamma-1} h_k(w))$, expanding a square, and recognizing the form of $\nabla h_k^2(w)$, for any $k$, we have
\begin{equation*}
|1-\gamma|^pI^\eta_k[w]=\int_\rn|h'_k(w)|^2\,|\nabla w|^2\,\V_D^{\gamma-1}\,\d\mu(y) + (1-\gamma)^2\int_\rn|h_k(w)|^2\,|y|^{2(2-\gamma)}\,\d\mu(y) +  \int_\rn(\nabla h_k^2(w))\cdot(\nabla\V_D^{\gamma-1})\,\V_D^{\gamma-1}\d\mu(y)\,.
\end{equation*}
By integrating by parts the last term, that will be justified below, one gets that
\begin{equation}\label{int.by.parts}
\begin{split}
|1-\gamma|^pI^\eta_k[w]&=\int_\rn|h'_k(w)|^2\,|\nabla w|^2\,\V_D^{\gamma-1}\,\d\mu(y) + (1-\gamma) \int_\rn|h_k(w)|^2\,|y|^{2(2-\gamma)}\,\d\mu(y)\\
&\quad -(1-\gamma) \int_\rn h_k^2(w)\,\V_D^{\gamma}\ \dv\left(y\,|y|^{1-\gamma}\left(\eta+(1-\gamma)|y|^{2-\gamma}\right)^{p-2}\right) {\dy}\,.
\end{split}
\end{equation}
Under the current assumptions and for any $s\in\left(1-\ve, 1+\ve\right)$, we can deduce the following relations between $h_2$ and $h_\gamma$
\begin{equation}\label{inq.pointwise}
\left(1-\varepsilon\right)^{2(2-\gamma)}\,\frac{h_\gamma(s)^2}{(1-\gamma)^2}\leq h_2(s)^2 \leq \left(1+\varepsilon\right)^{2(2-\gamma)}\,\frac{h_\gamma(s)^2}{(1-\gamma)^2}\qquad\mbox{ and }\qquad (h_2'(s))^2\leq \frac{(1+\ve)^{2(2-\gamma)}}{(1-\gamma)^2}\,(h'_\gamma(s))^2\,.
\end{equation}
Combining the above inequalities with identity~\eqref{int.by.parts} we obtain that
\begin{flalign*}
|1&-\gamma|^p I^\eta_2 [w]\\
\leq&  \frac{(1+\ve)^{2(2-\gamma)}}{(1-\gamma)^2}\,|1-\gamma|^p I^\eta_\gamma[w]+(1-\gamma)\int_\rn \left(\frac{(1+\ve)^{2(2-\gamma)}}{(1-\gamma)^2}h_\gamma^2(w)-h_2^2(w)\right)\,\V_D^{\gamma}\ \dv\left(y\,|y|^{1-\gamma}\left(\eta+(1-\gamma)|y|^{2-\gamma}\right)^{p-2}\right){\dy}\\
\leq & \frac{(1+\ve)^{2(2-\gamma)}}{(1-\gamma)^2}\,|1-\gamma|^p I^\eta_\gamma[w] +(1-\gamma)\left(\frac{(1+\ve)^{2(2-\gamma)}}{(1-\ve)^{2(2-\gamma)}}-1\right) \int_\rn h_2^2(w)\,\V_D^{\gamma}\,\left|\dv\left(y\,|y|^{1-\gamma}\left(\eta+(1-\gamma)|y|^{2-\gamma}\right)^{p-2}\right)\right|{\dy}\,,
\end{flalign*}
where in the last step we have used the first inequality in~\eqref{inq.pointwise} in order to control $h_\gamma^2$ with $h_2^2$. It only remains to compute the divergence in the last term of the above inequality. By a simple, though long, computation, one finds that
\[
\dv\left(y\,|y|^{1-\gamma}\left(\eta+(1-\gamma)|y|^{2-\gamma}\right)^{p-2}\right) = \frac{|y|^{1-\gamma}}{\left(\eta(1-\gamma)|y|^{2-\gamma}\right)^{3-p}}\,\left[\eta(N+1-\gamma)+(1-\gamma)(N+p-\gamma-1)|y|^{2-\gamma}\right]\,.
\]
Since both $(N+1-\gamma)$ and $(N+p-1-\gamma)$ are non-negative, the divergence  above has a sign. Moreover, we notice that there exists a constant $\mathcal{C}_{p,N}>0$,  depending on $N$ and $p$ but not on $\eta$, such that
$\left|\dv\left(y\,|y|^{1-\gamma}\left(\eta+(1-\gamma)|y|^{2-\gamma}\right)^{p-2}\right) \right|\leq \mathcal{C}_{p,N}$ for all $y\in\rn$. Combining all the above estimates and noting that $2\linE[v]=\int_{\rn} h_2^2(w)\V_D^\gamma \dy$, one obtains~\eqref{inq.linearised.quantities}.

Let us give some details for the justification of the above integration by parts applied to~\eqref{int.by.parts}. At first we notice that if $\linI^{(\eta)}_\gamma[v]=\infty$, there is nothing to prove, so we may assume that $\linI^{(\eta)}_\gamma[v] <\infty$. To deal with the opposite case we notice that all the integrals in this proof are well defined while restricted to a ball $B_R$, for $R>0$. What is more, all the computations are exactly the same up to an error term which comes from the integration by parts applied to~\eqref{int.by.parts}. Therefore, what remains to show is that this error term indeed vanishes in the limit as $R\rightarrow\infty$.  In order to do so, let us define $U(y):=y\,(1-\gamma)|y|^{1-\gamma}\V_D^{\gamma-1}(y)\mu(y)$. Then, by Green's identity, we have that for any $R>0$ it holds
\begin{equation*}
\int_{B_R}(\nabla h_k^2(w))\cdot(\nabla\V_D^{\gamma-1})\,\V_D^{\gamma-1}\d\mu(y)\, = \int_{B_R} \nabla h_k^2(w)\cdot U \dy = - \int_{B_R} h_k^2(w)\,\dv\,U\,\dy + \int_{\partial B_R} h_k^2(w) \, U\cdot \widehat{y}\d\sigma(y)\,,
\end{equation*}
where $\widehat{y}:=y/|y|$ and $\d\sigma$ is the surface measure of $\partial B_R$.  Notice that \[\dv\left(U\right)=\V_D^\gamma\,\dv\left(y\,|y|^{1-\gamma}\left(\eta+(1-\gamma)|y|^{2-\gamma}\right)^{p-2}\right) - \gamma(1-\gamma)\,|y|^{2(2-\gamma)}\,\mu(y)\,.\] Consequently, in the limit $R\rightarrow\infty$ one recovers the two last terms in~\eqref{int.by.parts}. Hence, we only need to show that the remainder term $\int_{\partial B_R} h_k^2(w) \, U\cdot \widehat{y}\d\sigma(y)$ converges to zero as $R\rightarrow\infty$. When $p>p_M$, by assumption~\ref{A1}, we have that $\left|h_k^2(w)\,U\cdot \widehat{y}\right|\leq C |h_2(w)|^2\,|U| \leq C |y|\,\V_D^\gamma$. As shown in Lemma~\ref{lem:ent-lin-ent}, $\V_D^\gamma$ is integrable as long as $p>p_M$. In turn,  $|y|\, \V_D^\gamma \leq C |y|^{1-N-\delta}$ for some $\delta>0$, which implies that the term $\int_{\partial B_R} h_k^2(w) \, U\cdot \widehat{y}\d\sigma(y)$ vanishes in the limit in the considered range since. Let us now concentrate on the case of $p\leq p_M$ and assumption~\ref{A2}. Thanks to Lemma~\ref{decay.lemma}  we know that $|v-\V_D|\leq C\,|y|^{-\frac{p}{(p-1)(2-p)}}$ for $|y|\geq 1$, so that $\left|h_k^2(w)\,U\cdot \widehat{y}\right|\leq C\,|y|^{-\frac{p}{(p-1)(2-p)}+1} $. Since $p>N(p-1)(2-p)+\delta$ for some $\delta>0$, we conclude that also in this case the remainder term vanishes in the limit $R\rightarrow\infty$. The proof is complete.
\end{proof}

Before we establish the final entropy -- entropy production inequality~\eqref{entropy-entropy-production-inq.intro} along the flow, we shall prove its linearised version. With this aim we employ the Hardy--Poincar\'e inequality provided as~\cite[Example 3.1]{Chlebicka2022} for  $q=2$, $\gamma=0$, $\beta=p/(p-1)$ and $\alpha=-1/(2-p)$.

\begin{prop}[Hardy--Poincar\'e inequality]\label{prop:hp}  Let $N\ge2$, $p\in(1,2)$, and $N<\frac{p}{(2-p)(p-1)}$.   Then there exists a finite constant  $C_{HP}=C_{HP}(p,N)>0$, such that for every compactly supported $\vp\in W^{1,\infty}(\rn)$ the following inequality holds true
 \begin{equation}
    \label{inq:ex-hp}
\int_\rn \ |\vp-\overline{\vp}|^2 (1+|y|^\frac{p}{p-1})^{-\frac{1}{2-p}}\dy\leq C_{HP} \int_\rn \ |\nabla \vp|^2 |y|^{2}(1+|y|^\frac{p}{p-1})^{-\frac{1}{2-p}} \dy,\end{equation}
where $\overline\vp$ is the average of $\vp$ with respect to $ (1+|y|^\frac{p}{p-1})^{-\frac{1}{2-p}}$.
\end{prop}

Before continuing, let us also define $M_\star:=\beta^\frac{N}{p}b_1^{\frac{p-1}{\beta p(2-p)}}$, where $b_1$ is as in~\eqref{b1} and $\beta$ as in~\eqref{beta}. Then
\begin{equation}\label{def:Mstar}
M_\star^\frac{\beta p(2-p)}{p-1}=\beta^{N\beta\frac{2-p}{p-1}}\ b_1\qquad\mbox{ and }\qquad M_\star=\int_{\rn} \V_1(y) \dy\,,
\end{equation}
with $\V_1$ as in~\eqref{vD} and~\eqref{Barenblatt-via-vD}. In what follows, for the sake of simplicity of the exposition, we shall consider solutions of mass equal to $M_\star$.

\begin{lem} \label{lem:linE-leq-linI-eps}
 Let $N\ge2$, $p_c<p < p_2$ and $0\leq v\in\LL^1(\rn)$, such that $M_\star=\int_{\rn}v (y)\dy$. Assume $v$ satisfies~\ref{A1} for some $\ve\in(0,1)$. If $p_c<p<p_M$, we additionally assume that $v$ satisfies~\ref{A2}. Then there exists a constant $\mathcal{C}=\mathcal{C}(p, N)>0$ such that
\begin{equation}\label{inq:linE-linI}
\mathcal{C}\,\linE[v]\leq \linI^{(\ve)}[v]\,.
\end{equation}
\end{lem}
\begin{proof}
Let $\vp$ be any function for which both sides of~\eqref{inq:ex-hp} are well defined. We notice that, even if~\eqref{inq:ex-hp} is stated for regular and compactly supported functions, the same inequality holds true for a larger class of functions through a~standard approximation procedure. Let also ${\overline{\vp}}=\int_{\rn}\vp(y)\d\mu_1(y)/\mu_1(\rn)$ where $\d\mu_1(y):=(1+|y|^{p'})^{-\frac{1}{2-p}}\dy$. Since $p\,V_1^{\gamma-1}\ge (2-p)\,(1+|y|^\frac{p}{p-1}) $, we obtain that
\begin{equation}\label{lin-ent-inq-HP}
\left(\frac{2-p}{p}\right)^\frac{1}{2-p}\inf_{c\in\RR} \int_{\rn} |\vp(y)-c|^2\,\V_1^{2-\gamma}(y)\dy \le \int_{\rn} |\vp(y)-\overline{\vp}|^2\,\d\mu_1(y)\,.
\end{equation}
It is known that the infimum on the left-hand side of the above inequality is achieved when $c=Z^{-1}\,\int_{\rn}\vp(y)\,\V_1^{2-\gamma}(y)\dy $ where $ Z=\int_{\rn}\,\V_1^{2-\gamma}(y)\dy$. We apply inequality~\eqref{lin-ent-inq-HP} to the function $\vp=\left(v-\V_1\right)\V_1^{\gamma-2}$ (which satisfies $0=\int_{\rn}\vp(y)\,\V_1^{2-\gamma}(y)\dy$) and we find
\begin{equation*}
2\,\left(\frac{2-p}{p}\right)^\frac{1}{2-p} \linE[v] \le \int_{\rn} |\vp(y)-\overline{\vp}|^2\,\d\mu_1(y)\,.
\end{equation*}
  It only remain to estimate the right-hand side of~\eqref{inq:ex-hp} by the right-hand side of~\eqref{inq:linE-linI}.  In order to do this, we observe that for any $y\in\rn$ we have
\[
\frac{|y|^2}{\left(1+|y|^\frac{p}{p-1}\right)^\frac{1}{2-p}}\frac{\left(\ve+(1-\gamma)|y|^{2-\gamma}\right)^{2-p}}{\V_1(y)}\leq \left(\frac{2}{p}\right)^\frac{p-1}{2-p}\,\left(2-\gamma\right)^{2-p}\,,
\]
where we used the fact that $|\nabla \V_1^{\gamma-1}(y)|=(1-\gamma)|y|^{2-\gamma}$. This is enough to prove that
\[
\int_\rn \ |\nabla \vp|^2 |y|^{2}\d\mu_1(y) \leq (1-\gamma)^p\,\left(\frac{2}{p}\right)^\frac{p-1}{2-p}\,\left(2-\gamma\right)^{2-p} \, \linI^{(\ve)}[v] \,.
\]
We conclude, therefore, that inequality~\eqref{inq:linE-linI} holds with the constant
\[
\mathcal{C}:=\frac{2}{C_{HP}\,(2-\gamma)^{2-p}}\,\frac{p^\frac{p-1}{2-p}}{2^\frac{p-1}{2-p}\,(1-\gamma)^p}\,\left(\frac{2-p}{p}\right)^\frac{1}{2-p}\,.
\]\end{proof}

\subsection{Convergence in $\LL^1$, case $p_D \le p<2$}\label{ssec:conv-above-p_D}

There is a special value of parameter $p$, above which it is well-known that one is equipped with strong tools. This value is
\begin{equation*}
p_D=\frac{2N+1}{N+1}\,,
\end{equation*}
already mentioned above. Before presenting this result let us recall that $p_c<p_D$ and let us refer to Introduction for more comments on other special values of parameter $p$.
When $p\in [p_D,2)$, the relative entropy functional is displacement convex, i.e. the entropy functional is convex along geodesics in the space of probability densities equipped with the Wasserstein metric, cf.~\cite{McCann}. Following~\cite{Agueh2008}, displacement convexity implies a relation between $\cE$ and $\cI$, that in turn shows exponential $\LL^1$-rates of convergence of solutions to \eqref{eq:v}, via a Gronwall's argument.

The main accomplishment of our paper is establishing rates of convergence in the range $(1,p_D)$, providing also a new proof when $p\in [p_D,2)$. Indeed, the parameter $p_D$ plays no role in our proof of Theorem~\ref{theo:RECR}: $p_D$ is relevant only for the validity of the optimal transportation method, but it is not critical for the asymptotic behaviour.  We show here a very short proof of the exponential convergence when $p\ge p_D$, following~\cite{Agueh2008}, without using the linearisation of $\cE$ or $\cI$.

\begin{prop}
   Let $N\ge1$, $p_D\leq p<2$, $0\leq \vo\in\LL^1(\rn)\cap \X$, and $D>0$ be such that $\int_{\rn}\VD(y)\dy=\int_{\rn}\vo (y)\dy$. Assume that $v$ is a weak solution to \eqref{eq:v} with $\vo$ as initial datum. Then there exists $c=c(\vo,p,N,D)$, such that for any $\tau>0$
\begin{equation}
    \label{conv-above-pD}
    \|v(\tau)-\VD\|_{\LL^1(\rn)}\leq c\, e^{-\tau/2}\,.
\end{equation}
\end{prop}

\begin{proof}
   Since $\vo\in\X$, we have, thanks to Remark~\ref{rmq} {\it (i)}, that given any $\tau_0>0$ the solution $v(\tau, \cdot)$ satisfies~\ref{A0} for any $\tau>\tau_0$. Then, thanks to  Lemma~\ref{lem:ent-lin-ent}, we have that $\cE[v(\tau)|\VD]<\infty$. Thus, by Lemma~\ref{CK.inq} we know that $\|v-\VD\|^2_{\LL^1(\rn)}
\leq c(p,N,D)\cE [v(\tau)|\VD]$. Moreover, since $p_D>p_c$, Lemma~\ref{lem:ent/ent-prod} implies that $\frac{\d}{\d\tau}\cE[v(\tau)|\VD]= - \cI[v(\tau)|\VD]$ for almost every $\tau>0$. By~\cite[Theorem~2.2]{Agueh2008}, within the range $p_D\leq p<2$ and for $\vo\in \LL^1(\RN)$, it holds
    \[\cE[v(\tau)|\VD]\leq \cI[v(\tau)|\VD]=-\tfrac{\d}{\d\tau}\cE[v(\tau)|\VD]\,.\]
Then Gronwall's Lemma implies $\cE[v(\tau)|\VD]\le c e^{-\tau}\cE[v(\tau_0)|\VD]$ for all $\tau\ge \tau_0$. Collecting all information we get~\eqref{conv-above-pD}.
\end{proof}

\subsection{Convergence in $\LL^1$,  case $p_M<p<2$}

In this section we shall provide a proof of the convergence in the $\LL^1$-norm in the range $p_M<p<2$. Note that for $p_M<p_D$ this result is new, while for $p_D\leq p<2$ we give a different proof that is not involving the optimal transportation tools (applied in~\cite{Agueh2009}). We stress again that  the parameter $p_D$ plays no role in this reasoning. For the simplicity of the exposition, we shall consider initial data with a fixed mass $\int_{\rn}\vo(y)\dy=M_\star $ with $M_\star$ being defined in~\eqref{def:Mstar}. We will recover full generality in the proof of Proposition~\ref{convergence.L1norm}.
\begin{prop}\label{prop:convergence.lp}
 Let $N\ge2$, $p_M<p<p_2$, $0\leq \vo\in\LL^1(\rn)\cap \X$,  be such that $\int_{\rn}\vo (y)\dy=M_\star$. Assume that $v$ is a weak solution to \eqref{eq:v} with $\vo$ as initial datum. If $\vo$ satisfies~\ref{A1} for some $\ve\in(0,1)$ and $D=1$, then there exists $\tau_\ve>0$,$c=c(\vo,p,N,\ve)$, and $\vartheta=\vartheta(p, N, \ve)>0$ such that for any $\tau>\tau_\ve$ it holds
 \begin{equation}\label{convergence.l1.pMp2}
  \|v(\tau)-\V_1\|_{\LL^1(\rn)}\leq c\, e^{-\vartheta\,\tau/2}\,.
 \end{equation}
\end{prop}

\begin{proof} Thanks to Remark~\ref{rmq} {\it (iii)} and the fact that $M_\star=\int_{\rn}\V_1 (y)\dy$, there exists $\tau_\ve>0$  such that the solution $v(\tau, \cdot)$ satisfies~\ref{A1} with $D=1$ for any $\tau>\tau_\ve$. Then, thanks to  Lemmata~\ref{lem:ent-lin-ent} and~\ref{lem:ent/ent-prod}, we have that $\cE[v(\tau)|\V_1]$, $ \cI[v(\tau)|\V_1]<\infty$ for all $\tau>\tau_\ve$. In what follows, we assume that $\tau>\tau_\ve$. By Lemma~\ref{lem:ABC} for every $\ve$ there exist $\kappa_0,\kappa_2>0$ such that $
\linI^{(\ve)}[v] \le \kappa_1\,\linI^{(\ve)}_\gamma[v] + \kappa_2\,\linE[v]$ and $\kappa_2\to 0$ as $\ve\to0.$ Let $\mathcal{C}$ be the constant from Lemma~\ref{lem:linE-leq-linI-eps} for which
$\mathcal{C}\,\linE[v]\leq \linI^{(\ve)}[v]$. Moreover, we recall that by Lemma~\ref{lem:I-est} we have
$C_\varepsilon\, \linI_\gamma^{(\varepsilon)} [v(\tau)]\leq \cI[v(\tau)|\V_1]$. We restrict attention to $\ve\in(0,1)$ small enough to ensure that $\kappa_2\leq \frac{\mathcal{C}}{2}$. Summing up, we get that
\begin{equation}\label{Fisher-lin-ent}
\tfrac{\mathcal{C}}{2}\linE[v] \le \frac{\kappa_1}{{C}_\ve} \cI[v(\tau)|\V_1]\,.
\end{equation}
On the other hand, by using Lemmata~\ref{lem:ent-lin-ent-2} and~\ref{lem:ent/ent-prod}  we obtain
\[
\tfrac{d}{d\tau}\cE[v(\tau)|\V_1]\leq -{C}_\ve\,\tfrac{\mathcal{C}}{2\kappa_1}(1+\ve)^{\gamma-2}\cE[v(\tau)|\V_1] \,.
\]
By using Gronwall's Lemma we can deduce from the above inequality that $\cE[v(\tau)|\V_1]\le  e^{-\,\vartheta\tau}\,\cE[v(\tau_\ve)|\V_1]$ for some $\vt=\vt(p,N,\ve)>0$ and all $\tau\ge \tau_\ve$. Then~\eqref{convergence.l1.pMp2} follows from Lemma~\ref{CK.inq}.
\end{proof}

\subsection{Convergence in relative error with rate, case $p_c<p\le p_M$}
 We finally address the case of $p_c<p\le p_M$. The main difference here is that the relative Fisher information is not bounded anymore, but it is merely an $\LL^1$-function of time. This technical difficulty will be overcome by the use of a~different version of Gronwall's Lemma, namely Lemma~\ref{weak.gronwall}.

\begin{prop}\label{convergence.lp.below.pm}
Let $N\ge2$, $p_c<p\le p_M$, $0\leq \vo\in\LL^1(\rn)\cap \X$,  be such that $\int_{\rn}\vo (y)\dy=M_\star$. Assume that $v$ is a weak solution to \eqref{eq:v} with $\vo$ as initial datum. If $\vo$ satisfies~\ref{A1}, for some $\ve\in(0,1)$ and $D=1$, and $\vo$ satisfies~\ref{A2}, for some $D_1,D_2>0$, then there exists $\tau_\ve>0$, $c=c(\vo, p,N,\ve)$ and $\vartheta=\vartheta(p, N, \ve)>0$ such that for any $\tau>\tau_\ve$ it holds
\begin{equation}\label{convergence.l1.below.om}
  \|v(\tau)-\V_1\|_{\LL^1(\rn)}\leq c\, e^{-\vartheta\,\tau/2}\,.
\end{equation}
\end{prop}

\begin{proof} Thanks to Remark~\ref{rmq} $(iii)$ and the fact that $M_\star=\int_{\rn}\V_1 (y)\dy$, there exists $\tau_\ve>0$  such that the solution $v(\tau, \cdot)$ satisfies~\ref{A1} with $D=1$ for any $\tau>\tau_\ve$. Then, thanks to  Lemma~\ref{lem:ent-lin-ent} we have that $\cE[v(\tau)|\V_1]<\infty$ for all $\tau>\tau_\ve$. In this case, thanks to Lemma~\ref{lem:ent/ent-prod-2}, we know that the function $\tau\mapsto \cI[v(\tau)|\V_1]$ is in $\LL^1(\tau_0, \infty)$, for any $\tau_0>0$. This implies that $\cI[v(\tau)|\V_1]<\infty$ \emph{only for almost every} $\tau>\tau_\ve$.  Our goal is to prove the following inequality for some $\vartheta>0$:
\begin{equation}\label{goal}
2\,\vartheta\,\int_{\tau}^\infty \cE[v(\tau)|\V_1]\dtau \le \cE[v(\tau)|\V_1]\qquad\forall\, \tau_0>\tau_\ve\,.
\end{equation}
Indeed, by using~\eqref{goal} with Lemma~\ref{weak.gronwall}, we can easily prove that $\cE[v(\tau)|\V_1]\le c e^{-2\,\vartheta\tau}$ for all $\tau\ge \tau_\ve$. Then inequality~\eqref{convergence.l1.below.om} will follow from Lemma~\ref{CK.inq}. In order to prove~\eqref{goal} let us start from identity~\eqref{weak.entropy.production}, namely
\begin{equation*}
\cE[v(\tau)|\V_1] = \int_{\tau}^\infty \cI[v(\tau)|\V_1]\dtau\qquad\forall\,\tau>0\,.
\end{equation*}
Let $A=\{\tau\in(\tau_\varepsilon, \infty):\ \cI[v(\tau)|\V_1]=\infty \}$. Since the function $\tau \mapsto I[v(\tau)| V_1]$ is in $ \LL^1(\tau_\varepsilon, \infty)$, the $1$-dimensional measure of $A$ is zero. We notice that,  since $\cE[v(\tau)|\V_1]$ is finite for every $\tau>0$, for every $\tau\in A$ and some $c>0$ it holds $\cI[v(\tau)|\V_1]\ge c\, \cE[v(\tau)|\V_1] $. On the other hand, on the set $A^c$, by proceeding as in the case of Proposition~\ref{prop:convergence.lp}, we can prove inequality~\eqref{Fisher-lin-ent} and, again with the use of Lemmata~\ref{lem:ent-lin-ent-2} and ~\ref{lem:ent/ent-prod},  we obtain
\begin{equation*}
{C}_\ve\,\tfrac{\mathcal{C}}{2\kappa_1}(1+\ve)^{\gamma-2}\cE[v(\tau)|\V_1] \leq \cI[v(\tau)|\V_1]\qquad\forall \tau\in A^c\,.
\end{equation*}
From this inequality and consideration what happens in $A$, one easily deduces~\eqref{goal} with $\vartheta: = {C}_\ve\,\tfrac{\mathcal{C}}{4\kappa_1}(1+\ve)^{\gamma-2}$.
\end{proof}

\subsection{Proof of \texorpdfstring{Theorem~\ref{theo:RECR}}{Theorem 1}}
We focus now on the proof of polynomial rate of $\LL^1$-convergence of solutions to~\eqref{CPLE} towards a Barenblatt profile. In order to obtain such result we shall first prove the following proposition.
\begin{prop} \label{convergence.L1norm}
Under the same assumptions of Theorem~\ref{theo:RECR}, there exists $\wt{T}=\wt{T}(p,N,M, \|\uo\|_{\X})>0$, $\wt{K}=\wt{K}(p,N,M, \|\uo\|_{\X})>0$ and $\nu=\nu(p,N)>0$ such that
\begin{equation}\label{convergence.relative.error-prop}
  \left\|u(t,\cdot)-\B_M(t,\cdot)\right\|_{\LL^1(\rn)} \le \wt{K}\,t^{-\nu}\quad\forall\, t\ge \wt{T}\,.
\end{equation}
\end{prop}
\begin{proof}
We observe that it is enough to consider solutions whose initial datum $\uo$ has initial mass equal to $M_\star$, that is $\int_{\rn}\uo(x)\dx=M_\star$, for $M_\star$ as in~\eqref{def:Mstar}. Indeed, let $u$ be a weak solution to~\eqref{CPLE} with $\uo$ as initial datum and let $M=\int_{\rn}\uo(x)\dx$. Then,  as previously observed in~\cite[Preliminaries]{Bonforte2020b}, by defining
\begin{equation}\label{change.mass}
\widetilde{u}(t,x):=\tfrac{M_\star}{M}\, u\left(t\,\left(\tfrac{M_\star}{M}\right)^{2-p},x\right)
\end{equation}
one gets $\widetilde{u}$ being a solution to~\eqref{CPLE} with initial datum $\uo\,M_\star/M$ and mass $M_\star$. Once~\eqref{convergence.L1norm} will be obtained for $\widetilde{u}$, we can rescale it back with the use of the identity~\eqref{change.mass} and get the same inequality for $u$. This follows since the same mass changing formula~\eqref{change.mass} applies also to the family of Barenblatt solutions. After this computation we get~\eqref{convergence.L1norm}, where the constant $\wt{K}$ changes its value, but the rate $\nu$ remains the same.

In the same way, we observe that it is enough to consider inequality~\eqref{convergence.relative.error-prop} for any among the Barenblatt solutions $\B_{M_\star}(t+T,\cdot)$ (for $T\geq 0$) and not necessarily $\B_{M_\star}(t,\cdot)$. Indeed, we may notice that there exists a constant $C=C(T, M_\star)>0$ such that
\begin{equation}\label{change.profile}
\|\B_{M_\star}(t+T,\cdot)-\B_{M_\star}(t,\cdot)\|_{\LL^1(\rn)}\leq {C(T)}t^{-1}\qquad\forall t>T\,.
\end{equation}
For a moment, let us assume the above inequality. Then, from a convergence result with respect to the profile $\B_{M_\star}(t+T,\cdot)$, we retrieve
\begin{equation*}
\begin{split}
\|u(t,\cdot)-\B_{M_\star}(t,\cdot)\|_{\LL^1(\rn)}&\leq \|u(t,x)-\B_{M_\star}(t+T,\cdot)\|_{\LL^1(\rn)}+\|\B_{M_\star}(t+T,\cdot)-\B_{M_\star}(t,\cdot)\|_{\LL^1(\rn)}\\
&\leq \tfrac{\wt{K}}{t^\nu} +  \tfrac{C(T)}{t}\leq \tfrac{\wt{K}+C(T)}{t^\nu}\,,
\end{split}
\end{equation*}
for any $t>\max\{1, T, \wt{T}\}$. We have used above that $\nu\leq1$. Let us justify it. The rate obtained in~\eqref{change.profile} is optimal, which can be proven by a direct computation. Since Proposition~\ref{convergence.L1norm}  also covers the case considered in~\eqref{change.profile}, we conclude that $\nu\leq1$ by its optimality.

On the other hand, inequality~\eqref{change.profile} can be proven by considering the relative error between the two solutions $\B_{M_\star}(t+T,\cdot)$ and $\B_{M_\star}(t,\cdot)$. Indeed, for any fixed $t>0$, the supremum in $|x|$ of the quotient $\B_{M_\star}(t+T,\cdot)/\B_{M_\star}(t,\cdot)$ is attained either at $0$ or at $\infty$ (one can prove this through a simple, however lengthy, computation). From this observation one finds that
\begin{equation*}
\left\|\frac{\B_{M_\star}(t+T,\cdot)}{\B_{M_\star}(t,\cdot)}-1\right\|_{\LL^\infty(\rn)}\leq \frac{C(T, M_\star)}{t}\qquad\forall t>T\,.
\end{equation*}
It is also direct to see that the above inequality is optimal.  Then, inequality~\eqref{change.profile} is obtained through the following computation
\begin{equation}\label{inq.different.barenblatt}
\int_{\rn}|\B_{M_\star}(t+T,x)-\B_{M_\star}(t,x)|\dx\leq \left\|\frac{\B_{M_\star}(t+T,\cdot)}{\B_{M_\star}(t,\cdot)}-1\right\|_{\LL^\infty(\rn)}\,\int_{\rn} \B_{M_\star}(t,x) \dx \leq M_\star\,\frac{C(T, M_\star)}{t}\,.
\end{equation}
We also notice that the optimality in~\eqref{change.profile} can be deduced by a similar reasoning as above.

We will present arguments for ranges $p_M< p <2$ and $p_c<p\leq p_M$ separately. In both cases, as we have noticed before,  it is sufficient to consider the mass of the initial datum being equal to $M_\star$.

Let us consider the case $p_M< p <2$. As we explained, it is enough to consider  to compute the rate of convergence towards the Barenblatt profile $\B_{M_\star}(t+\beta, \cdot)$. Under the current assumption, we can perform the change of variables~\eqref{change.variables} and consider a solution $v$ to~\eqref{eq:v} with initial datum $\vo$ and mass $M_\star$. By Remark~\ref{rmq} (items {\it (ii)} and {\it (iii)}) we know that, for any $\ve\in(0,1)$ there exists $\tau_\ve>0$ such that   $v$ satisfies~\ref{A1} (with $\ve$ and $D=1$) for any $\tau>\tau_\ve$. Therefore, all hypotheses of Proposition~\ref{prop:convergence.lp} are satisfied and inequality~\eqref{convergence.l1.pMp2} holds true. By rescaling back inequality~\eqref{convergence.l1.pMp2}, one finds exactly~\eqref{convergence.relative.error-prop}. Hence, the claim is proven in this case.

Lastly, we are considering the case $p_c<p\leq p_M$. Under the additional hypothesis, namely~\eqref{stronger.GHP}, we shall employ the profile $\B_{M_\star}(t+T,\cdot)$ for the convergence result. Therefore, in the change of variable~\eqref{change.variables}, instead of using $R_\beta(t)$, we shall use $R_T(t)$ defined in~\eqref{R(t)} with $T$ from~\eqref{stronger.GHP}. Let us define also
\begin{equation*}
v(\tau,y):=R_T(t)^N\,u(t,x)\,,
\end{equation*}
where $y$ and $\tau$ are as in~\eqref{change:rho-u}.
Then $v$ is a solution to~\eqref{eq:v} with $\vo(y)={R}_T(0)^{N}\uo(y{R}_T(0))$. We notice that $\uo\in\X$, so  by Remark~\ref{rmq} {\it (ii)} we know that, for any $\ve\in(0,1)$, there exists $\tau_\ve>0$ such that $v(\tau,\cdot)$ satisfies~\ref{A1} for all $\tau>\tau_\ve$. At the same time, thanks to assumption~\eqref{stronger.GHP} and after the change of variables,  we get that the initial datum $\vo$ satisfies~\ref{A2}. Therefore, $v(\tau,\cdot)$ satisfies both~\ref{A1} and~\ref{A2} for any $\tau>\tau_\ve$. Eventually, by taking $v(\tau_\ve,\cdot)$ as the initial datum, the assumptions of Proposition~\ref{convergence.lp.below.pm} are satisfied and, hence,~\eqref{convergence.l1.below.om} holds true. By re-scaling back, we find that for all $t$ large enough it holds
\[
\|u(t,x)-\B_{M_\star}(t+T,\cdot)\|_{\LL^1(\rn)}\leq {\wt{K}}{t^{-\nu}}\,.
\]
By inequality~\eqref{inq.different.barenblatt} we find the wanted result.\end{proof}

We are in the position to justify our main accomplishment.
\begin{proof}[Proof of Theorem~\ref{theo:RECR}]
Once Proposition~\ref{convergence.L1norm} is proven, inequality~\eqref{convergence.relative.error} directly results from Theorem~\ref{theo:conv-L1-RE}.\end{proof}

\section{Convergence of derivatives of radial  solutions}\label{sec:conv-radial}

This section is devoted to providing an exhaustive answer to~\eqref{Q2} and prescribing an explicit rate for the uniform converge in relative error of radial derivatives of radial solutions in three cases separately. In particular we prove Theorems~\ref{Thm:gradient.convergence} and~\ref{relative.error.radial.derivatives} for problem~\eqref{eq:v} or, equivalently, Theorems~\ref{original-variables:p>pc},~\ref{radial.derivatives.pc} and~\ref{relative.error.radial.derivatives.original} for the original problem~\eqref{CPLE}. Here,  we exploit a stronger relation between radial derivatives of solutions to~\eqref{CPLE} and radial solutions to a weighted version of the FDE, established in~\cite{Iagar2008}. In fact, the radial formulation of~\eqref{CPLE} can be rewritten for $u$ being a function of $(r=|x|, t)$ as follows
\begin{equation}\label{ple.radial}
\partial_t u = r^{1-N} \, {\partial_r} \left(r^{N-1}\,|{\partial_r u}|^{p-2}\,\partial_r u \right)\,.
\end{equation}
Let us consider  $\Phi: \mathbb{R}^N\times(0, \infty)\rightarrow \mathbb{R}$ being a \emph{non-negative} function of $(\vr=|x|, t)$ and a solution to
\begin{equation}\label{fde.radial}
\partial_t \Phi= \vr^{1-n} \, \partial_{\vr}\left(\vr^{n-1}\, \partial_{\vr} \Phi^m\right)\,,\qquad m=p-1\,,\end{equation}
where $ n=N+2\tfrac{N}{p'}$ is a positive parameter. In Section~\ref{ssec:information_on_related_radial_classical_and_weighted_fast_diffusion_equations} we point out in more details in what sense a derivative of a solution to~\eqref{fde.radial} solves~\eqref{ple.radial} and comment on related results.

\subsection{Information on related radial classical and weighted Fast Diffusion Equations}\label{ssec:information_on_related_radial_classical_and_weighted_fast_diffusion_equations}
Below we present a radial equivalence due to~\cite[Theorem 1.2]{Iagar2008}. We point out, however, that we use a particular version of the result stated in~\cite{Iagar2008}. The original result is indeed much stronger and valid also for sign-changing solutions. Notice also that we consider a slightly different equation, so the constant $\mathcal{D}$ defined below differs slightly comparing with the one of~\cite[Theorem 1.2]{Iagar2008}. In what follows we shall denote by $r=|x|$ in the \eqref{ple.radial} case and by $\vr=|x|$ the coordinates for the~\eqref{fde.radial} equation.
\begin{theo}[{\cite[Theorem 1.2]{Iagar2008}}]\label{radial.transformation.thm} Suppose $2<n<\infty$. If $u$ is a  radially symmetric and decreasing solution  of equation \eqref{ple.radial}, then $\Phi$ being a non-negative solution of \eqref{fde.radial} is related to $u$ through the following transformation:
\begin{equation}\label{radial.transformation}
-\partial_r u(t,r)= {\mathcal{D}}\,\vr^\frac{2}{m+1}\,\Phi(t, \vr)\,, \quad \mathcal{D}= \left(\frac{2m}{m+1}\right)^\frac{2}{m-1}\,,
\end{equation}
 where $ r=\vr^{\frac{2m}{m+1}}$ and the correspondence of the parameters is given by
\begin{equation}\label{identity.parameters}
p=m+1,\,\quad N= (n-2)\frac{(m+1)}{2m}\,.
\end{equation}
\end{theo}
Let us note that in \cite{Iagar2008} the authors also analyze the case $0<n<2$, however, we shall not use those results. We also remark that, even if the solutions to~\eqref{CPLE} are at least $C^{1,\alpha}$, a priori we do not know whether $\Phi$ is well-defined at the origin by transformation~\eqref{radial.transformation}. We will address these issues below.

In fact, in Theorem~\ref{radial.transformation.thm} the following transformation is defined as
\begin{equation}\label{parameter.transformation}
  (p, N) \mapsto (m, n, \aa):=\left(p-1\,, 2+2\tfrac{N}{p'}\,, N-2-2\tfrac{N}{p'}\right)\,.
\end{equation}
Recall that $p'=(p-1)/p$ and note that $\aa+n=N$, so the above map is injective.  Let us recall that $p_c=\frac{2N}{N+1}$ and $p_Y=\frac{2N}{N+2}$, and  $p_Y<p_c$. We have the following ranges for different values of $p$:
\begin{enumerate}[{\it(i)}]
\item if $p_c<p<2$, then
  \begin{equation}\label{good.range}
     \tfrac{n-2}{n}<m<1\,,\quad N+1<n<N+2\,,\quad \mbox{and}\quad -2<\aa<-1\,;
  \end{equation}
  \item if $p=p_c$, then $ \   m=\tfrac{N-1}{N+1}=\tfrac{n-2}{n}\,,\quad n=N+1\,,\quad\mbox{and}\quad \aa=-1\,$;
  \item if $ p_Y<p<p_c$, then $\tfrac{n-2}{n+2}<m<\tfrac{n-2}{n}\,,\quad N< n < N+1\,,\quad \mbox{and} \quad  - 1 < \aa < 0\,$;
  \item if $p=p_Y$ (the Yamabe case), then $\
    m=\tfrac{N-2}{N+2}\,,\quad n=N\,,\quad\mbox{and}\quad \aa=0\,$.
    \end{enumerate}
In general, the artificial dimension $n$ is not an integer: this happens only in the limit case $p=p_c$ and $p=p_Y$. As already noticed in~\cite{Iagar2008}, the only case when $n=N$ (and also when the weight $\aa=0$) is the case $p=p_Y$. We also remark that, when $N=2$ the value $p_Y=1$ and it is excluded from our analysis. \newline

Note that the equation
\begin{equation}\label{ckn.fde.radial}\tag{WFDE}
\partial_t \Phi =
|x|^\aa\,\dv\left(|x|^{-\aa}\,\nabla \Phi^m\right)\,
\end{equation} written for radial solutions is exactly~\eqref{fde.radial} for the choice of parameters from~\eqref{parameter.transformation}. This equation is sometimes referred to as the Weighted FDE with Caffarelli--Kohn--Nirenberg weights, see~\cite{Bonforte2017a, Bonforte2017b, Bonforte2020}. We also stress that radial initial data produce radial solutions for~\eqref{ckn.fde.radial}. In order to infer the asymptotics of  derivatives of radial solutions we exploit the relation between radial solution to~\eqref{CPLE} and~\eqref{ckn.fde.radial} together with known properties of solutions to~\eqref{ckn.fde.radial}.
Therefore a solution to~\eqref{fde.radial} is a radial solution to~\eqref{ckn.fde.radial}, cf.~\cite{Iagar2008} and Proposition~\ref{from.Cauchy.to.Cauchy}. In the same spirit, we notice that $\Phi$ as a function of $(\vr,\,t)$  is a radial solution to the original (unweighted) FDE when $n=N$. We stress that in~\eqref{fde.radial}, parameter $n$ plays the role of an artificial dimension and is not an integer in general. It is unusual to consider equations in a continuous dimension, however, in the radial case, this allows us to unveil some unexpected features.  Let us remark that equation~\eqref{ckn.fde.radial} shares many features with~\eqref{CPLE}, as it was already observed in~\cite{Bonforte2020}. For values of $m\in(0,1)$ there are three different ranges where solutions to~\eqref{ckn.fde.radial} behave differently. As for~\eqref{CPLE}, we call the interval $\frac{n-2}{n}<m<1$ the \emph{good range} corresponding to $p_c<p<2$, as observed in~\eqref{good.range}. In this range solution to the Cauchy problem associated to~\eqref{ckn.fde.radial} conserve weighted mass ($\int_{\RN}\Phi(t,x)\,|x|^{-\aa}\dx=  \int_{\RN}\Phi(0,x)\,|x|^{-\aa}\dx$) once the initial datum $\Phi(0,x)\in \LL^1(\RN, |x|^{-\aa}\,\dx)$. Similarly, in this range~\eqref{ckn.fde.radial} also admits a family of self-similar solutions given by ~\eqref{bareblatt.profile.ckn}, that is commonly called Barenblatt solution, see~\cite{Bonforte2017a,Bonforte2017b,Bonforte_2019,Bonforte2020}. For parameters from~\eqref{good.range}  Barenblatt solutions are defined as follows
\begin{equation}\label{bareblatt.profile.ckn}
\mathfrak{B}_M(t,x)= t^\frac{1}{1-m}\left[\left(a_1\,t^{2\theta}\,M^{2\theta(m-1)}+a_2|x|^2\right)\right]^\frac{1}{1-m}\,,\quad\mbox{where}\quad \theta=\frac{1}{2-n(1-m)}\,.
\end{equation}
The constants $a_1$ and $a_2$ given by
\begin{equation}\label{a1-a2}
\int_{\RR^N}(a_1+ a_2 |x|^2)^\frac{1}{1-m}\,|x|^{-\aa} \dx=1\,,\quad\mbox{while}\quad a_2=\frac{1-m}{2m}\,\theta\,,
\end{equation}
where $\theta$ is as in \eqref{bareblatt.profile.ckn}. We remark that the mass $M$ of the profile $\mathfrak{B}_M$ is computed with respect to the measure $|x|^{-\aa}\dx$, that is
\[
M=\int_{\RR^N} \mathfrak{B}_M(t,x)\,|x|^{-\aa}\dx\,.
\]
Lastly, as the reader may suspect, the Barenblatt profile $\B_M$ of equation~\eqref{CPLE} and the one $\mathfrak{B}_{\overline{M}}$ of equation~\eqref{ckn.fde.radial} are related by formula~\eqref{radial.transformation}, namely
\begin{equation}\label{B:Bbar}
- \partial_r \B_M(t,r) = {\mathcal{D}} \, \varrho^{\frac{2}{m+1}} \mathfrak{B}_{\mathfrak{C}M}(t, \varrho)\,,
\end{equation}
where $r$, $\vr$, and $\mathcal{D}$ are as in Theorem~\ref{radial.transformation.thm}. We remark that the mass of $\mathfrak{B}$ is corrected by a multiplicative factor $\mathfrak{C}=\mathfrak{C}(N, p)>0 $ given by
\begin{equation*}
 \mathfrak{C}=\frac{pN}{2(p-1)\mathcal{D}}\,.
\end{equation*}

 When $m\leq \frac{n-2}{n}$, the Barenblatt solutions do not exist anymore as solutions generated by a {$\delta_0$} as initial datum. Nevertheless, a pseudo-Barenblatt profile is still available, for any $T>0$ and $D>0$ let us define
 \begin{equation}\label{bareblatt.profile.ckn.mc}
 \mathfrak{B}_{D, T}(t,x)=\mathfrak{R}_T(t)^{n}\,\mathfrak{U}_D(x\,\mathfrak{R}_T(t))\quad\mbox{where}\quad \mathfrak{U}_D(x):=\left(D + \frac{1-m}{2\,m}|x|^2\right)^\frac{1}{1-m}\,,
 \end{equation}
 and where
 \begin{equation}\label{time.rescale.mc}
   \mathfrak{R}_T(t):=\left(\frac{T-t}{|\theta|}\right)_{+}^\theta\quad\mbox{if}\,\,0<m<\frac{n-2}{n}\quad\mbox{and}\quad \mathfrak{R}_T(t):=\exp\{\mathfrak{l}(t+T)\}\quad\mbox{if}\,\, m=\frac{n-2}{n}\,,
 \end{equation}
 where $\theta$ is as in~\eqref{bareblatt.profile.ckn} (that is negative in the case $0<m<\frac{n-2}{n}$) and $\mathfrak{l}>0$ is a free parameter.
 We conclude this section by noticing that identity~\eqref{B:Bbar} also holds for the range $p_Y\leq p \leq p_c$, and we have that
 \begin{equation}\label{barenblat.derivative.pc}
    - \partial_r \B_{D,T}(t,r) = {\mathcal{
    D}} \, \varrho^{\frac{2}{m+1}} \mathfrak{B}_{\overline{\mathfrak{C}}D, T}(t, \varrho)\,,
  \end{equation}
where   $\overline{\mathfrak{C}}=\overline{\mathfrak{C}}(N, p)>0$ is given by
\begin{equation}\label{kappa.below.pc}
  \overline{\mathfrak{C}}=\left(\frac{1-m}{2m}\right)^{2\theta+1}\,.
\end{equation}

Let us explain how solutions to~\eqref{CPLE} are related to solutions to the Cauchy problem of~\eqref{ckn.fde.radial}. While at the level of solutions this is given directly by the transformation~\eqref{radial.transformation}, it is not clear what happens to the initial data and in what sense identity~\eqref{radial.transformation} should be understood.  Before giving a complete answer, let us fix the notation which will be used in what follows. Let us denote by $u:(0,\infty)\times \rn\rightarrow[0, \infty)$ the solution to~\eqref{CPLE} with a radial initial datum $\uo(x)$. Since the solution $u(t,x)$ is radial, with an abuse of notation, we shall denote $u(t,x)$ by $u(t,r)$ and $\uo(x)$ by $\uo(r)$. The function $\Phi:(0,\infty)\times \rn\rightarrow[0, \infty)$ will be a solution to~\eqref{ckn.fde.radial} with initial datum $\Phi_0(x)$. Again, in the case of radial initial datum $\Phi_0(x)=\Phi_0(r)$, we shall denote the solution $\Phi(t,x)$ by $\Phi(t, \vr)$. The following proposition answers the main questions of this section.

\begin{prop}\label{from.Cauchy.to.Cauchy}
Suppose that $N\ge2$\,, $1\leq p<2$, ${\mathcal{D}}$ is as in~\eqref{radial.transformation}, and $n$, $m$, $\aa$ are as in~\eqref{identity.parameters}. Let $u$ be a radial solution to~\eqref{CPLE} with a radial initial datum $\uo\in C^2(\R^N)$ satisfying  one of the following conditions:
\begin{itemize}
  \item[(i)] if $p_c < p <2$ we assume that~\eqref{assumptions.thm.derivative}  holds for some $A>0$ and $R_0>0$;
  \item[(ii)] if $N=2$ and $1=p_Y<p\leq p_c$, or  $3\leq N \leq 6$ and  $ p_Y\leq p \leq  p_c$, or $N>6$ and $p_2<p<p_c$, we assume that there exist $D_1$, $D_2>0$ and $T>0$ such that~\eqref{hp.derivative} holds;
  \item[(iii)] if $N>6$ and $p_Y\leq p \leq p_2$,  we assume that there exists $D_1$, $D_2>0$ and $T>0$ such that~\eqref{hp.derivative} holds and that there exist $\wt{D}>0$ and $f\in \LL^1((0,\infty), r^{n-1}\dr)$ with such that~\eqref{D} holds.
\end{itemize}
Then $\Phi_0$ given by
\begin{equation*}
\Phi_0(\vr):=-\tfrac{1}{\mathcal{D}}\vr^{-\frac{2}{1+m}}\,\left(\partial_r \uo\right)(\vr^\frac{2m}{1+m})\qquad\forall \vr>0\,,
\end{equation*}
 satisfies $0\le \Phi_0\in\LL^1_{\rm loc}(\rn, |x|^{-\aa}\dx)$  and the following Cauchy problem
\begin{equation}\label{fde.cauchy.problem}
\begin{cases}
  \partial_t \Phi= |x|^{\aa}\,\dv\left(|x|^{-\aa}\,\nabla \Phi^m\right)\,\quad\text{for }\ (t,x)\in\left(0, \infty\right)\times \RR^N\\
  \Phi(0,x)=\Phi_0(|x|)\, \quad\text{for }\ x\in\RR^N\,,
\end{cases}
\end{equation}
is solvable. Moreover, its solution $\Phi(t, \cdot)$ belongs to $\LL^\infty_{{\rm loc}}(\RR^N)$ for any $t>0$ and it is related to $u$ by transformation~\eqref{radial.transformation}.
\end{prop}
\begin{proof}
Let us start with justifying that for $p_Y\leq p<2$, the initial datum $0\le \Phi_0\in\LL^q_{{\rm loc}}(\rn, |x|^{-\aa}\dx)$ for any $0<q<\frac{n(1+m)}{2(1-m)}$. Simple computations show that $\frac{n(1+m)}{2(1-m)}>1$ under the current assumptions. To motivate the abovementioned integrability of $\Phi_0$, we recall that $\uo\in C^2(\rn)$, and $\uo$ is radial and decreasing, we have that $\partial_r u_r(0)=0$ and $|\partial_r u_r(\vr)|\leq C \vr$ close to the origin. Therefore, we find that
\[
\vr^{n-1}\, |\Phi_0(\vr)|^q\leq C\, \vr^{\fb}\quad\mbox{where}\quad \fb=q\,\tfrac{2(p-2)}{p}+2\,N\,\tfrac{(p-1)}{p}+1=n-1+q\tfrac{2(m-1)}{1+m}\,.
\]
Note that $\fb>-1$ as long as $q<\frac{n(1+m)}{2(1-m)}$. We recall that $\int_{|x|\leq 1} |\Phi_0(x)|^q |x|^{-\aa}\dx=\omega_N\,\int_0^1 |\Phi_0(\vr)|^q\, \vr^{N-1-\aa}\,{\rm d} \vr=\omega_N\,\int_0^1 |\Phi_0(\vr)|^q\, \vr^{n-1}\,{\rm d} \vr$, where $\omega_N$ is the area of the $N$-dimensional sphere. \newline

Further we proceed case by case.\newline

\noindent\textit{Case {\it (i)}: $p_c<p<2$.} Let us consider the integrability of $\Phi_0$. From the last inequality of~\eqref{assumptions.thm.derivative} we deduce that for any $\vr>R_0^\frac{1+m}{2m}$ we have
\begin{equation}\label{Phi_0-dec-1}
  \Phi_0(\vr)\leq C\,\vr^{-\frac{2}{1+m}-\frac{2}{1-m}\frac{2m}{1+m}}=C\, \vr^{-\frac{2}{1-m}}\,.
\end{equation}
Therefore, for any $\vr>R_0^\frac{1+m}{2m}$ it holds that
\begin{equation*}
\vr^{n-1}\,\Phi_0(\vr)\leq C\, \vr^{n-1-\frac{2}{1-m}}\,.
\end{equation*}
We notice that in this case ($p_c<p<2$ equivalently to $\frac{n-2}{n}<m<1$), we have that $\frac{2}{1-m}-n>0$. Therefore, the quantity $|\vr^{n-1}\,\Phi_0(\vr)|$ is integrable with respect to the Lebesgue measure and the initial datum $\Phi_0\in\LL^1(\rn, |x|^{-\aa}\,\dx)$.

We are in the position to pass to justification of solvability of~\eqref{fde.cauchy.problem}. Let us  notice that the equation in~\eqref{fde.cauchy.problem} is the same as in~\eqref{fde.radial}. As explained in~\eqref{good.range}, in the present range of parameters we always have $\aa<0$. Solutions for problem~\eqref{fde.cauchy.problem} have been constructed in~\cite{Bonforte2017b} in the same spirit as in \cite[Theorem 2.1]{Herrero_1985}. We stress that in~\cite[Proposition 7]{Bonforte2017b} the initial datum is assumed to be in $\LL^\infty(\RR^N)$. This assumption can be weakened to  merely $\Phi_0\in\LL^1(\RR^N, |x|^{-\aa}\dx)$ by a standard approximation procedure as it is done in the proof of~\cite[Theorem 2.1]{Herrero_1985}. Since $\Phi_0\geq 0$  and the comparison principle holds due to~\cite[Corollary 9]{Bonforte2017b}, we know that the solution is non-negative. In this range of parameters solutions are bounded since $\Phi_0\in \LL^1_{{\rm loc}}(\rn, |x|^{-\aa}\,\dx)$, see~\cite[Theorem 1.2]{Bonforte_2019}. It is also known that solutions are at least $C^\alpha$-regular close to the origin (see~\cite[Theorem 1.8]{Bonforte_2019}) and $C^\infty$-smooth outside of the origin. This has been already remarked in~\cite[Lemma 11]{Bonforte2017b}. See also~\cite[Section 21.5.3]{DiBenBook} where the authors affirm that local analyticity in space and, at least, Lipschitz continuity in time holds for solutions to a general equation of the form~\eqref{fde.cauchy.problem}. This considerations prove that the solution $\Phi$ to~\eqref{fde.cauchy.problem} exists and it has the wanted properties. It only remains to verify that a radial solution $u$ to~\eqref{CPLE} is related to the Cauchy problem~\eqref{fde.cauchy.problem} through the transformation~\eqref{radial.transformation}. Despite this seems obvious, it is not since we need to find a relevant relation between initial data.  Let us consider the auxiliary function
\begin{equation}\label{auxiliary.function}
  \overline{u}(t,r):= \,{\mathcal{D}}\, \int_{r}^\infty \Phi(t, s^\frac{p}{2(p-1)})\, s^\frac{1}{p-1}\, \ds\,,
\end{equation}
where ${\mathcal{D}}$ is as in~\eqref{radial.transformation}. By the global Harnack principle (see~\cite[Theorem 1.1]{Bonforte2020} and cf.~\eqref{ghp.inq}) for solutions to the radial problem~\eqref{fde.cauchy.problem} and by \eqref{Phi_0-dec-1} we infer that the following result holds: for any $\tau_0>0$ there exist $M_1(\tau_0), M_2(\tau_0)>0$ and $\tau_1(\tau_0), \tau_2(\tau_0)>0$, such that we have
\begin{equation}\label{ghp.ckn.radial}
\mathfrak{B}_{M_1}(t-\tau_1,\varrho)\leq \Phi(t,\varrho)\leq \mathfrak{B}_{M_2}(t+\tau_2,\varrho)\qquad\forall \varrho>0\quad \forall t\geq 2\,\tau_0\,.
\end{equation}
From the above estimates, we deduce that there exists a constant $C=C(t, \tau_0)>0$ such that for any $t \ge 2\tau_0$ we have
\begin{equation}\label{time.constant}
\Phi(t, s^\frac{p}{2(p-1)})\, s^\frac{1}{p-1} \leq C(t, \tau_0)\, \frac{s^\frac{1}{p-1}}{s^{\frac{2}{1-m}\frac{p}{2(p-1)}}} = C(t, \tau_0)\, s^{-\frac{2}{2-p}}\,,\quad\forall s>0
\,,\end{equation}
where  $C(t, \tau)\lesssim t^{-n\,\theta}$ for any $t>2\tau$ with $\theta$ as in~\eqref{bareblatt.profile.ckn}. The exponent $2/(2-p)>1$ since $2>p>1$, so we deduce from the above inequality that the function $\overline{u}$ is well defined. Furthermore, we have that $\overline{u}(t)\in\LL^\infty(0,\infty)$ for any $t>0$ (indeed, $\tau_0$ is chosen arbitrarily).  Let us now investigate the regularity of $\Phi$. The validity of the inequality~\eqref{ghp.ckn.radial} allows us to use the regularity information resulting from the proof of~\cite[Lemma 11]{Bonforte2017b}. By those results we have that $\Phi \in C^\infty(0,\infty)^2$, and, for any $\tau>0$, $\varepsilon>0$ and $k>0$ there exist $C_1=C_1(t, \tau, \varepsilon)>0$ and $C_2=C_2(t, \tau, \varepsilon, k)$ such that
\begin{equation}\label{regularity.decay}
  |\partial_t \Phi(t,\varrho)|\leq {C_1(t)}{\varrho^{-\frac{2}{1-m}}}\quad\mbox{and}\quad \left| \tfrac{\partial^k}{\partial \varrho^k}
  \big(\partial_t\Phi
  (t, \varrho)\big)\right| \leq {C_2(t)}{\varrho^{-\frac{2}{1-m}-k}}\quad \forall t\ge \tau\quad\forall \varrho\geq\varepsilon\,.
\end{equation}
The above estimates allow us to differentiate in $t$ and in $\varrho$ under the sign of the integral in~\eqref{auxiliary.function}. Consequently, $\overline{u}\in C^2(0,\infty)^2$ and $\overline{u}$ solves equation~\eqref{ple.radial} almost everywhere in $(t,r)\in \left(0,\infty\right)^2$. Furthermore, function $\overline{u}$ is a weak solution to the following Neumann problem
\begin{equation}\label{neumann.problem}
  \begin{cases}
    \partial_t \overline{u} = r^{1-N} \,  {\partial_r} \left(r^{N-1}\,|\partial_r \overline{u}|^{p-2}\,\partial_r \overline{u} \right)\quad \text{for }\ (t,r)\in\left(0,\infty\right)^2\,,\\
    \partial_r \overline{u}(t,0)  =0\,, \\
    \overline{u}(0,r)  =\overline{u}_0(r)\,\quad\text{for }\  r\in\left[0, \infty\right)\,.
  \end{cases}
\end{equation}
We briefly comment on the literature for the above problem in Remark~\ref{rem:neumann}. Let us continue with the rest of the proof. By using the result of~\cite[Theorem III.8.1]{DiBenedetto1990}, we know that the function $(t,x)\mapsto \nabla u(t,x) \in C^\alpha_{{\rm loc}}((0,\infty)\times \RR^N)$, which is enough to guarantee that $\partial_r {{u}}(t,0)={0}$ for all $t>0$. We conclude therefore that $u$ also solves problem~\eqref{neumann.problem}. Since $\uo=\overline{u}_0$ we would like to conclude that $u=\overline{u}$ by using the uniqueness result for~\eqref{neumann.problem}. This would be enough to conclude the proof, since by the construction we will have that $\Phi$ and $u$ satisfy the relation~\eqref{radial.transformation}. However, in order to apply the uniqueness result of~\cite[Theorem II.1]{DiBenedetto1990} we need to ensure that
\begin{equation}\label{time.derivative}
{\partial_t \overline{u}}(t,r)\leq \mathfrak{h}\, \overline{u}(t,r)\quad\mbox{a.e.}\,\,(t,r)\in(0,\infty)^2\,,
\end{equation}
where $\mathfrak{h}=\mathfrak{h}(N, p, t)$ is independent of $\overline{u}$. We notice that, as observed in~\cite[p.~45]{DiBenedetto1990}, solutions to problem~\eqref{CPLE} satisfy~\eqref{time.derivative} by the construction. In order to prove~\eqref{time.derivative} for solutions to~\eqref{neumann.problem} we shall use a modification of a trick due to B\'enilan and Crandall~\cite{Benilan1981} provided in~\cite[Lemma III.3.4]{DiBenedetto1990}. By the comparison principle proven in~\cite{Bonforte2017b}, the uniqueness for~\eqref{fde.cauchy.problem} is guaranteed. Let us consider $\Psi^{\lambda}$  being the unique solution to~\eqref{fde.cauchy.problem}  with initial datum
\[
\Psi_0^{\lambda}(\varrho):=\Psi^{\lambda}(0, \varrho)=\lambda^\frac{1}{m-1}\, \Phi_0(\varrho)\quad \text{for }\ \lambda>0 \,.
\]
Notice that if $\lambda\ge 1$, then $\Psi_0^{\lambda}(\varrho)\leq \Phi_0(\varrho)$ for all $\varrho\ge0$. The homogeneity of~\eqref{fde.cauchy.problem} implies that $\Psi$ can be written as
\[
\Psi^{\lambda}(t,\varrho) = \lambda^\frac{1}{m-1}\,\Phi(\lambda\, t, \varrho)\,.
\]
Therefore, again by the comparison principle, we have that $\Psi^{\lambda}(t,\varrho)\leq \Phi(t, \varrho)$ for all $(t,\varrho)\in(0, \infty)^2$. By setting $\lambda=1 + h/t$, for a small $h>0$, we obtain that, for any $(t, \varrho)\in (0, \infty)^2$ it holds
\begin{equation*}
\begin{split}
  \Phi(t+h, \varrho) - \Phi(t, \varrho) & = \Phi(\lambda\,t,\varrho) - \Phi(t, \varrho)   = \lambda^\frac{1}{1-m}\,\lambda^\frac{1}{m-1}\Phi(\lambda\,t,\varrho) - \Phi(t, \varrho)  = \lambda^\frac{1}{1-m}\,\Psi^{\lambda}(t,\varrho) -  \Phi(t, \varrho) \\& \leq \left(\lambda^\frac{1}{1-m}-1\right)\, \Phi(t, \varrho)\,.
\end{split}
\end{equation*}
By  the Mean Value Theorem  we infer that for some $\zeta\in(0, \frac{h}{t})$ it holds
\begin{equation*}
  \Phi(t+h, \varrho) - \Phi(t, \varrho) \leq \frac{h}{(1-m)\,t}\, \left(1+\zeta\right)^\frac{m}{1-m} \Phi(t, \varrho)\,.
\end{equation*}
By using the above inequality in $(t, s^\frac{p}{2(p-1)})$, multiplying by ${\mathcal{D}}$ (for ${\mathcal{D}}$ being as in~\eqref{radial.transformation}), and integrating with respect to the measure $s^\frac{1}{p-1}\ds$, we find that
\begin{equation}\label{time.estimate.bar}
\overline{u}(t+h, r) - \overline{u}(t,r) \leq \frac{{\mathcal{D}}\,h}{(1-m)\,t}\,\left(1+\zeta\right)^\frac{m}{1-m} \,\overline{u}(t,r)\,.
\end{equation}
We notice that, by the same computations with $\lambda\leq 1$, we shall establish inequality~\eqref{time.estimate.bar} with a reversed sign and $h<0$. Then we divide by $h$ both sides of~\eqref{time.estimate.bar} and take the limit for $h\rightarrow 0$. Let us point out that $\overline{u}\in C^1(0,\infty)^2$, so the left-hand side of~\eqref{time.estimate.bar} converges to $\partial_t \overline{u}(t,r)$. In turn, we find estimate~\eqref{time.derivative} with $t\, (1-m)\,\mathfrak{h}={\mathcal{D}}$. The proof in the case $p_c<p<2$ is complete.\newline

\noindent\textit{Case {\it (ii)}: $N=2$ and $1=p_Y<p\leq p_c$, or $3\leq N \leq 6$ and  $ p_Y\leq p \leq  p_c$, or $N>6$ and $p_2<p<p_c$}. We shall explain the main differences between this case and the above one. Let us consider first the case $p<p_c$ and then pass to $p=p_c$ in the end of the proof. First of all, we notice that, by assumption~\eqref{hp.derivative} and relation~\eqref{barenblat.derivative.pc}, there exists $\overline{D}_1, \overline{D}_2>0$ such that
\begin{equation}\label{initial.datum.mc}
  \mathfrak{B}_{\overline{D}_1, T}(0, \varrho) \leq \Phi_0(\varrho) \leq \mathfrak{B}_{\overline{D}_2, T}(0, \varrho)\quad\forall \varrho\ge0\,,
\end{equation}
the case $\varrho=0$ being obtained as a limit case. Notice also that $\overline{D}_i = \overline{\mathfrak{C}}\; D_i$, where $D_i$ is as in~\eqref{hp.derivative} and $\overline{\mathfrak{C}}$ as in~\eqref{kappa.below.pc}.  As in the previous case, the result~\cite[Propostion 7]{Bonforte2017b} is enough to establish the existence of a non-negative solution $\Phi\in \LL^\infty_{{\rm loc}}(\RN)$. The comparison principle has been established in~\cite[Corollary 9]{Bonforte2017b} for initial data which satisfies~\eqref{initial.datum.mc} and the following assumption: there exists $\overline{D}>0$ such that
\begin{equation}\label{initial.datum.mc.2}
  \Phi_0(\varrho)= \mathfrak{B}_{\overline{D},T}(0, \varrho) + f(\varrho) \quad\forall\varrho\ge0\,,
\end{equation}
for $f\in\LL^1((0, \infty), r^{n-1}\dr)$ (notice that parameter $n$ defined in Theorem~\ref{relative.error.radial.derivatives} is the same as in~\eqref{identity.parameters}). In the present range the assumption~\eqref{initial.datum.mc.2} easily follows from~\eqref{initial.datum.mc}, since the difference of two Barenblatts $\mathfrak{B}_{\overline{D}_2, T}-\mathfrak{B}_{\overline{D}_1, T}$ is always integrable if $m>m_\star$, where
\begin{equation*}
  m_\star=\frac{n-4}{n-2}= \frac{2N(p-1)-2p}{2N(p-1)}\,.
\end{equation*}
Two remarks are in order. Firstly, the fact that the difference of two Barenblatt is integrable can be proven by technniques similar to those in Lemma~\ref{decay.lemma}. We also refer to~\cite[Section 2.1]{Bonforte2017b} and~\cite[Introduction]{Blanchet2009} for a general discussion. Secondly, it is easy to see that $m>m_\star$ if and only if
\begin{equation*}
  2N(p-1)^2-2N(p-1)+2p\ge0\quad\mbox{and}\quad p>1\,.
\end{equation*}
A simple computation shows that the above condition always holds for $N< 6$. For $N\ge6$ it is satisfied for $p\in\left(1, p_1\right)\cup \left(p_2, 2\right)$, which explains the appearance of the exponent $p_2$.

We deduce that, by the comparison principle, inequality~\eqref{initial.datum.mc} continues to hold for $t>0$. More precisely
\begin{equation}\label{comparison.t.positive}
  \mathfrak{B}_{\overline{D}_1, T}(t, \varrho) \leq \Phi(t, \varrho) \leq \mathfrak{B}_{\overline{D}_2, T}(t, \varrho)\quad\forall \varrho\ge0\quad\mbox{and}\quad0< t\leq  T\,,
\end{equation}
which proves that $\Phi(t, \varrho)=0$ for all $t\ge T$ and $\varrho\ge0$. Inequality~\eqref{comparison.t.positive} plays the role of inequality~\eqref{ghp.ckn.radial} in this range of parameters. Indeed, from~\eqref{comparison.t.positive} one can deduce~\eqref{time.constant} for any $0<t<T$ which is enough to establish that $\overline{u}$ is well defined also in the present case. At the same time,  using again~\cite[Lemma 11]{Bonforte2017b}, we have that $\Phi \in C^\infty(0,\infty)^2$ and inequalities from~\eqref{regularity.decay} hold also in the present case. This is enough to show that $\overline{u}$ is a weak solution to~\eqref{neumann.problem}. Moreover, using the same argument as above, one can easily prove inequality~\eqref{time.derivative}, which is the missing condition to verify in order to use the uniqueness result of~\cite[Theorem II.1]{DiBenedetto1990}. We have explained in Section~\ref{ssec:reg} that in our setting for a solution $u$ to~\eqref{CPLE}, the function $(t,x)\mapsto \nabla u(t,x) \in C^\alpha((0, \infty)\times \RN)$. Hence, we are in the position to guarantee that $\partial_r u(t, 0)=0$ (which means that $u$ is also a solution to~\eqref{neumann.problem}) and therefore, by the uniqueness, we have that $u=\overline{u}$. This concludes the proof in the case $p<p_c$.

In the case $p=p_c$, the above proof is also valid. The only thing that changes is that~\eqref{comparison.t.positive} holds for any $t>0$. In this case, the solution lives for all $t>0$ as for $p_c<p<2$ and there is no extinction in finite time.\newline

\noindent\textit{Case (iii): $N \geq 6$ and  $ p_Y\leq p \leq  p_2$}.  The present case is very similar to \textit{Case {\it (ii)}}. Indeed, the main difference is that identity~\eqref{initial.datum.mc.2} is not a consequence of~\eqref{initial.datum.mc} but instead it needs to be assumed from the very beginning. Nevertheless, notice that~\eqref{initial.datum.mc.2} is exactly assumption~\eqref{D} rewritten after the use of transformation~\eqref{radial.transformation}. The rest of the proof follows exactly the lines of \textit{Case {\it (ii)}}. Therefore the proof is complete.

\end{proof}

\begin{rem}[On the Neumann problem]\rm \label{rem:neumann}
    To the best of our knowledge, the problem~\eqref{neumann.problem} has not been investigated yet. It seems, as well, that the Neumann problem for $p$-Laplacian type equation has been much less studied. For more information we refer to~\cite{Andreu_2011,Alikakos_1982} in the case of a Neumann problem in bounded domains, to~\cite{Alikakos_1981} for the Neumann problem for the Porous Medium Equation ($\partial_t u=\Delta u^m$, $m>1$), and to~\cite[Chapter~11]{JLVSmoothing} {for exposition of the background in detail}. In dimension $N=1$, the techniques used in the seminal paper~\cite{Esteban_1988} can be adapted (at least in the good range $p_c<p<2$) in order to prove the existence, uniqueness and comparison principle. We also stress that problem~\eqref{neumann.problem} is very similar to~\eqref{CPLE} and the techniques of~\cite{DiBenedetto1990} can be adapted in the whole generality for the entire range $1<p<2$.
\end{rem}

\subsection{Proof of the  convergence in the relative error of the radial derivatives}

It is convenient to rescale~\eqref{ckn.fde.radial} in the way we are able to consider at the same time the supercritical, critical and subcritical range. The following change of variables is very much in the same spirit of~\eqref{change.variables}. Consider $\Phi$ to be a~solution to~\eqref{ckn.fde.radial} and let us define $\Psi$ as
\begin{equation}\label{change.fkn}
  \Psi(\tau, y):= \mathfrak{R}_{T}(t)^n\, \Phi(t, x)\quad\mbox{where}\quad \tau=\log\frac{\mathfrak{R}_T(t)}{\mathfrak{R}_T(0)}\quad\mbox{and}\quad  y:=\frac{x}{\mathfrak{R}_T(t)}\,,
\end{equation}
where $\mathfrak{R}_T$ is as in~\eqref{time.rescale.mc}. We recall that the definition of $\mathfrak{R}_T$ differ when $p=p_c$ ($m=m_c$) and $p<p_c$ ($m<m_c$). However, in both cases, if $\Phi$ satisfies~\eqref{ckn.fde.radial} then the problem satisfied by $\Psi$ is the following
\begin{equation}\label{rescaled.ckn.fde.radial}
  \begin{cases}
  \partial_\tau \Psi= |y|^{\aa}\,\dv\big[|y|^{-\aa}\,\left(\nabla \Psi^m-y\,\Psi\right)\big]\,\quad\text{for }\ (\tau,y)\in\left(0, \infty\right)\times \RR^N\,,\\
  \Psi(0,y)=\Psi_0(|y|)\, \quad\text{for }\ y\in\RR^N\,,
\end{cases}
\end{equation}
where the initial datum $\Psi_0(y) = \mathfrak{R}_{T}(0)^n\, \Phi_0(x\,\mathfrak{R}_{T}(0)^n)$.

There are two main advantages which justify the introduction of the change of variables of \eqref{change.fkn}. The first reason is that, in the case $p<p_c$ ($m<m_c$), on the contrary to the solution to~\eqref{ckn.fde.radial} which extinguishes in finite time $T$ (as does the Barenblatt function), the rescaled solution $\Psi$ lives for any $0<\tau<\infty$. The second reason is that~\eqref{rescaled.ckn.fde.radial} admits the stationary solution
\begin{equation*}
  \mathfrak{U}_D(x):=\left(D + \frac{1-m}{2\,m}|x|^2\right)^\frac{1}{1-m}\,,\qquad D>0\,,
\end{equation*}
introduced in~\eqref{bareblatt.profile.ckn.mc}. When $m>\frac{n-2}{n}$, the parameter $D$ is related to the mass of $\mathfrak{U}_D$, i.e. $\int_{\rn} \mathfrak{U}_D |x|^{-\aa} \dx$.

We notice that conditions~\eqref{derivatives.supercritical}, \eqref{hp.derivative}, and~\eqref{initial.datum.mc} imply the existence of $\overline{D}_1$, $\overline{D}_2>0$ such that
\begin{equation}\label{H1.rescaled}
  \mathfrak{U}_{\overline{D}_1}(\varrho)\leq \Psi_0(\varrho)\leq \mathfrak{U}_{\overline{D}_2}(\varrho)\quad\forall\varrho\ge0\,,
\end{equation}
while, condition~\eqref{D} (or, equivalently,~\eqref{initial.datum.mc.2}) translates to the existence of $\overline{D}>0$ such that
\begin{equation}\label{H2.rescaled}
\Psi_0(\varrho) = \mathfrak{U}_{\overline{D}}(\varrho) + f(\varrho)\quad\forall\varrho\geq0\,,
\end{equation}
where $f\in\LL^1((0, \infty)\,,r^{n-1}\dr)$.

In what follows, we refer to the result~\cite[Theorem 5]{Bonforte2017b}  for $\gamma\leq0$, which can be stated in our language as follows. Under assumptions~\eqref{H1.rescaled} and~\eqref{H2.rescaled} there exist $D>0$,  $\tau_\bullet>0$, $C_\bullet>0$, and $\Lambda>0$ such that
\begin{equation}\label{convergence.derivatives}
  \left\|\frac{\Psi(\tau)}{\mathfrak{U}_{\overline{D}}}-1\right\|_{\LL^\infty(\RN)}\leq C_\bullet\, e^{-2\,\frac{(1-m)^2}{2-m}\Lambda (\tau - \tau_\bullet)}\quad\forall \tau\ge \tau_\bullet
\end{equation}
The convergence rate $\Lambda=\Lambda(n,m)$ is the optimal constant in a relevant Hardy--Poincar\'e inequality related to~\eqref{ckn.fde.radial}. We refer to~\cite[Proposition 3]{Bonforte2017b} for more information, see also~\cite{Chlebicka2022}. For the sake of completeness, we state here the different values of $\Lambda$ for various parameters. If $m\leq \frac{n}{n+2}$ (notice that $m=\frac{n}{n+2}$ means $p=p_M$), then
\begin{equation*}
  \Lambda=\Lambda_{\rm ess}:=\frac{\big((n-2)(1-m)-2\big)^2}{4(1-m)^2}\,.
\end{equation*}
On the other hand, when $m>\frac{n}{n+2} $ ($p>p_M$), then
\begin{equation*}
  \Lambda=\min\left\{ \Lambda_{\rm ess}\,, \frac{2\,\eta}{1-m}\,, \frac{2(2-n(1-m))}{1-m}\right\}\quad\mbox{where}\quad\eta:=\sqrt{N-1+\left(\frac{N-2-\aa}{2}\right)^2} - \frac{N-2-\aa}{2}\,.
\end{equation*}
A detailed inspection of the proof reveals that the time $\tau_\bullet$ cannot be quantified a priori and depends on the initial datum $\Phi_0$, see in particular~\cite[Proposition 14 and Section 3.2]{Bonforte2017b}. However, the constant $C_\bullet$ can be explicitly quantified and it depends on the initial datum $\Psi_0$ (through its entropy), the parameters $\overline{D}_1$, $\overline{D}_2$ and, of course, $n$ and $\aa$, see in particular~\cite[Proof of Theorem 4]{Bonforte2017b}.

\begin{proof}[Proof of Theorem~\ref{Thm:gradient.convergence}] Since $p>p_c$ and we assume~\eqref{assumptions.thm.derivative}, we can make use of Proposition~\ref{from.Cauchy.to.Cauchy}~{\it (i)}. Namely, the solution $u $ is related to the solution $\Phi $ of~\eqref{fde.cauchy.problem} by the transformation~\eqref{radial.transformation}.  We notice that assumption~\eqref{derivatives.supercritical} is more restrictive than~\eqref{assumptions.thm.derivative}, so also in this case Proposition~\ref{from.Cauchy.to.Cauchy}~{\it (i)} applies. Then, by~\cite[Theorem~1.4]{Bonforte2020b}, we know that
\begin{equation}\label{strong.convergence.derivatives}
  \left\|\frac{\partial_r u(t,\cdot)}{\partial_r \mathcal{B}_M(t,\cdot)}-1\right\|_{\LL^\infty(0,\infty)}=\left\|\frac{\Phi(t,\cdot)}{\mathfrak{B}_{\mathfrak{C}M}(t,\cdot)}-1\right\|_{\LL^\infty(0,\infty)}\longrightarrow0\quad\mbox{as}\quad t\longrightarrow\infty.
\end{equation}
Since the mass of $\Phi$ is conserved in time due to~\cite[Proposition 10]{Bonforte2017b}, the above display shows that
\[
\mathfrak{C}M=\int_{\RR^N}\,\Phi(t,x)\,|x|^{-\aa}\,\dx = \omega_N\,\int_0^\infty \Phi(t, \varrho)\,\varrho^{n-1}\,\textrm{d}\varrho\,.
\]
To conclude our proof, we only need to obtain a convergence rate towards zero for the uniform relative error which appears in the middle of~\eqref{strong.convergence.derivatives}. In the case $\frac{n}{n+2}\leq m<1$ (recall that $m=\frac{n}{n+2}$ means $p=p_M$),  this  can be inferred from~\cite{Bonforte_2023}. More precisely,  under a condition that we shall discuss below, from~\cite[Theorem 7]{Bonforte_2023} it follows that there exist explicit constants $\varepsilon_\star=\varepsilon_\star(m, N, \aa)>0$, $C_\star=C_\star(m, N, \aa,\Phi_0)>0$, and  $\lambda=\lambda(m, N, \aa)$ such that for any $0<\varepsilon<\varepsilon_\star$ we have
\begin{equation}\label{rate.convergence.derivatives}
\left\|\frac{\Phi(t,\cdot)}{\mathfrak{B}_{\mathfrak{C}M}(t,\cdot)}-1\right\|_{\LL^\infty(0,\infty)}\leq \varepsilon\qquad\forall t\ge C_\star \varepsilon^\frac{1}{\lambda}\,.
\end{equation}
While the result in~\cite[Theorem 7]{Bonforte_2023} is stated only for $\frac{n-1}{n}\leq m<1$, the method can be easily extended up to $m=\frac{n}{n+2}$ since it is based on a weaker form of assumption~\eqref{H1.rescaled}. Indeed, the proofs in~\cite{Bonforte_2023} are based on the Global Harnack Principle for equation~\eqref{ckn.fde.radial}, namely inequality~\eqref{ghp.ckn.radial}, and the fact that the second moment with respect to the measure $|x|^{-\aa}$ is finite, which holds exactly for $\frac{n}{n+2}<m<1$.

Once inequality~\eqref{rate.convergence.derivatives} is obtained, establishing the convergence rate of~\eqref{gradient.decay.convergence} requires only the inversion of the relation between $\varepsilon$ and~$t$. This has been done in detail in~\cite[Corollary 4.14]{Bonforte2020a}. In the case $\frac{n-2}{n}<m<\frac{n}{n+2}$, to infer the rate of convergence in~\eqref{strong.convergence.derivatives}, we invoke~\cite[Theorem 5]{Bonforte2017b}, indeed, that result guarantee an explicit convergence rate for~\eqref{strong.convergence.derivatives} (see inequality~\eqref{convergence.derivatives}) when the initial datum $\Psi_0$ satisfies both~\eqref{H1.rescaled} and~\eqref{H2.rescaled}. We remark that assumption~\eqref{H1.rescaled} is nothing than~\eqref{derivatives.supercritical}. Lastly, that in the current regime $\frac{n-2}{n}<m<\frac{n}{n+2}$, assumption~\eqref{H2.rescaled} can be easily obtained from~\eqref{H1.rescaled}, since in this regime the difference of two Bareblatts profile is always integrable.

To conclude the proof, let us briefly comment on the last restriction of~\cite[Theorem 7]{Bonforte_2023}:  the initial datum should satisfy
\begin{equation}\label{X.condition.phi}
\|\Phi_0\|_{\mathcal{Y}_m}:= \sup_{R>0} R^{\frac{2}{1-m}-n}\,\int_{|x|>R} \Phi_0(x)\,|x|^{-\aa}\,\dx\, <\infty\,.
\end{equation}
We stress that the condition $\|\Phi_0\|_{\mathcal{Y}_m}<\infty$  plays the same role for \eqref{fde.cauchy.problem} as~\eqref{UCRE/Xp} for \eqref{CPLE}.  This has already been pointed out in~\cite{Bonforte2020}.
In our setting $\Phi_0$ verifies \eqref{X.condition.phi} due to Proposition~\ref{from.Cauchy.to.Cauchy}. We also remark that such a condition is satisfied uniformly in $\frac{n-2}{n}<m<1$.
\end{proof}

\begin{proof}[Proof of Theorems~\ref{relative.error.radial.derivatives}  and~\ref{radial.derivatives.pc}]We consider first the case $p<p_c$. By reversing the change of variables~\eqref{change.fkn} and using the convergence rate~\eqref{convergence.derivatives} one gets inequality
\begin{equation}\label{wanted.inequality}
\Big\|\frac{\Phi(t,\cdot)}{\mathfrak{B}_{\overline{D},T}(t, \cdot)}-1\Big\|_{L^\infty(\RN)}\leq C_\diamond\, (T-t)^{-\lambda }\,\quad \forall\,\, t_\diamond < t <T\,,
\end{equation}
where $t_\diamond$ is such that
\begin{equation*}
  \tau_\bullet =\log\frac{\mathfrak{R}_T(t_\diamond)}{\mathfrak{R}_T(0)}\,,
\end{equation*}
\begin{equation*}
  C_\diamond=\frac{C_\bullet\, e^{2\,\frac{(1-m)^2}{2-m}\,\Lambda\,\tau_\bullet}}{|\theta|^\frac{2\,\theta\,(1-m)^2\,\Lambda}{2-m}}\,\mathfrak{R}_T(0)^\frac{2\,\theta\,(1-m)^2\,\Lambda}{2-m}\,,\qquad\mbox{and}\qquad \lambda=\frac{2\,\theta\,(1-m)^2}{2-m}\,\Lambda\,.
\end{equation*}
In the case $p=p_c$, the only change in inequality~\eqref{wanted.inequality} is the fact that the right-hand side is of the form $t^{-\lambda}$ and the inequality holds for any $t\ge t_\diamond$.
Once~\eqref{wanted.inequality} is obtained one can easily obtain the first inequality of~\eqref{convergence.rata.below.pc} (respectively, the first inequality of~\eqref{convergence.rate.pc}) by using the relations~\eqref{barenblat.derivative.pc} and~\eqref{radial.transformation}.

It only remains to prove the second inequality of~\eqref{convergence.rata.below.pc} (respectively, of~\eqref{convergence.rate.pc}).  We notice that, inequality~\eqref{wanted.inequality} can be rewritten in the following form. For any $T>t>t_\diamond$ and $\varrho\ge0$ it holds
\begin{equation*}
  -\varepsilon(t)\, \mathfrak{B}_{\overline{D},T}(t, \varrho) \leq \Phi(t,\varrho) - \mathfrak{B}_{\overline{D},T}(t, \varrho)\leq \varepsilon(t)\,\mathfrak{B}_{\overline{D},T}(t, \varrho)\quad\mbox{where}\quad \varepsilon(t)=C_\diamond\, (T-t)^{-\lambda }\,.
\end{equation*}
By integrating the above inequality as in~\eqref{auxiliary.function}, using the relation between $u$ and $\Phi$ (explained in Proposition~\ref{from.Cauchy.to.Cauchy}) and the relation among the different Barenblatt solutions (exposed in~\eqref{bareblatt.profile.ckn.mc}), one finds the following link between $u$ and $\B_{D,T}$:
\begin{equation*}
  -\varepsilon(t)\, \B_{D,T}(t, r) \leq  u(t, r) - \B_{D,T}(t, r)\leq \varepsilon(t)\,\B_{D,T}(t, r)\quad \forall\,r\ge0\quad\forall\,\, t_\diamond < t < T\,,
\end{equation*}
which is equivalent to the second inequality of~\eqref{convergence.rata.below.pc} (respectively, of~\eqref{convergence.rate.pc} upon choosing $\varepsilon(t)=C_\diamond\, t^{-\lambda }$). The proof is then concluded.
\end{proof}

\section{Proof of Theorem~\texorpdfstring{\ref{theo:Entropy-decay}}{2} and related results}\label{ssec:proof_of_theorem_2}

The proof of Theorem~\ref{theo:Entropy-decay} follows the lines of the proof of Propostions~\ref{prop:convergence.lp} and~\ref{convergence.lp.below.pm}. However, in order to achieve \emph{almost optimality} and the \emph{optimality} for a certain class of radial solutions, we need to be much more careful on constants of inequalities~\eqref{cal-I-est}, ~\eqref{inq.linearised.quantities}, and~\eqref{inq:linE-linI}. In order to do so we collect here several results of general ineterest. Let us begin with the counterpart of~\eqref{cal-I-est}.

\begin{prop}\label{prop:cal-I-est:improved}
Let $N\ge1$, $1<p<2$, $0\leq v\in C^{1, \alpha}(\mathbb{R}^N)$ for some $0<\alpha<1$, and $D>0$ such that $v$ satisfies~\ref{A1} for some $\ve\in\left(0,1\right)$ and
\begin{equation}\label{hp}
1-\varepsilon \leq \frac{\partial_r v}{\partial_r V_D} \leq 1+\varepsilon\,.
\end{equation}
Then we have
\begin{equation}\label{cal-I-est:improved}
\frac{(1-\varepsilon)^{\mathsf{a}}}{(1+\varepsilon)^{2-p}}\, (p-1)\, \linI^{(0)}_\gamma[v] \leq \cI[v|\VD] \leq \frac{(1+\varepsilon)^{\mathsf{a}}}{(1-\varepsilon)^{2-p}}\, (p-1)\, \linI^{(0)}_\gamma[v]\,,
\end{equation}
where $\mathsf{a}=1+(1-\gamma)(2-p)$.
\end{prop}
\begin{rem}
We notice that in Proposition~\ref{prop:cal-I-est:improved} the function $v$ is not necessarily a solution to any equation. Instead is just a smooth function that satisfy assumptions~\ref{A1} and~\eqref{hp}. We also observe that $\linI^{(0)}_\gamma[v]$ is nothing else than $\linI^{(\eta)}_\gamma[v]$ defined in~\eqref{def-I-ve-gamma} with the choice $\eta=0$.
\end{rem}
\begin{proof}
Let us call $\mathbf{G}(\phi,\psi)=(\nabla \phi-\nabla \psi)\cdot(\bb[\nabla \phi]-\bb[\nabla \psi])$ for functions $\phi,\psi:\mathbb{R}^N\to\mathbb{R}$, where $\bb[\nabla \phi]:=|\nabla \phi|^{p-2}\nabla \phi$. We have the identity $|\gamma-1|^p\,\cI[v|\VD]=\int_{\rn} v(y)\,\mathbf{G}(v^{\gamma-1}, \ \V_D^{\gamma-1})\dy$. As a first step, thanks to inequality~\ref{A1}, we deduce that
\begin{equation*}
  \frac{(1-\varepsilon)}{|\gamma-1|^p}\int_{\rn}\V_D(y)\,\mathbf{G}(v^{\gamma-1},  \V_D^{\gamma-1}) \dy \leq \cI[v|\VD] \\ \leq \frac{(1+\varepsilon)}{|\gamma-1|^p}\int_{\rn}\V_D(y)\,\mathbf{G}( v^{\gamma-1},  \V_D^{\gamma-1}) \dy\,.
\end{equation*}
We also notice that, since $\partial_r \V_D\leq 0$ and $\gamma\leq 1$, we conclude that $\partial_r \V_D^{\gamma-1}=(\gamma-1)\V_D^{\gamma-2}\partial_r \V_D\ge0$. By a very similar computation and by using~\eqref{hp} we find that $\partial_r v^{\gamma-1}\ge0$. In a similar way, we deduce that
\begin{equation}\label{needed.inequality.fisher}
  \frac{1-\varepsilon}{(1+\varepsilon)^{1-\gamma}}\partial_r \V_D^{\gamma-1}\leq \partial_r v^{\gamma-1} \leq \frac{1+\varepsilon}{\left(1-\varepsilon\right)^{1-\gamma}} \partial_r \V_D^{\gamma-1}\,.
\end{equation}
 Recall that, for a $C^1$ radial function $f(x)=f(r)$, we have that $\nabla f(x)= \partial_r f(r) \frac{x}{|x|}$ for any $x\in\mathbb{R}^N\setminus\{0\}$ and $\nabla f(0)=0$. Let us call that $\partial_r v^{\gamma-1}=\xi\ge0$ and $\partial_r \V_D^{\gamma-1}=\eta\ge0$, so that, for any $y\neq0$
\begin{equation*}\begin{split}
  \mathbf{G}(v^{\gamma-1},\V_D^{\gamma-1}) &= \big(\xi\frac{y}{|y|}-\eta\frac{y}{|y|}\big)\cdot\Big(|\xi|^{p-2}\xi\frac{y}{|y|}-|\eta|^{p-2}\eta\frac{y}{|y|}\Big) \\
  & = (\xi-\eta)(|\xi|^{p-2}\xi-|\eta|^{p-2}\eta)= |\xi|^{p}-|\eta|^{p-2}\eta\xi - |\xi|^{p-2}\xi\eta+|\eta|^p\,.
\end{split}\end{equation*}
Therefore, from inequality~\eqref{numerical.inequality.fisher} of Lemma~\ref{numerical_lemma_fisher} and inequality~\eqref{needed.inequality.fisher} we deduce that
\begin{multline*}
  \frac{(1-\varepsilon)^{(1-\gamma)(2-p)}}{(1+\varepsilon)^{2-p}}\,|\nabla \V_D^{\gamma-1}(y)|^{p-2}\,|\nabla v^{\gamma-1} - \nabla \V_D^{\gamma-1}|^2 \\ \leq \frac{\mathbf{G}(v^{\gamma-1}, \V_D^{\gamma-1})}{p-1} \leq \frac{(1+\varepsilon)^{(1-\gamma)(2-p)}}{(1-\varepsilon)^{2-p}}|\nabla \V_D^{\gamma-1}(y)|^{p-2}\,|\nabla v^{\gamma-1} - \nabla \V_D^{\gamma-1}|^2\,.
\end{multline*}
By integrating the above inequality one easily obtains~\eqref{cal-I-est:improved}.
\end{proof}

Let us now discuss the counterpart of inequality~\eqref{inq.linearised.quantities}. Since in Proposition~\ref{prop:cal-I-est:improved} we obtain the inequality which links $\cI[v(\tau), \VD]$ with $\linI^{(0)}_\gamma[v]$, therefore it is natural to consider inequality~\eqref{inq.linearised.quantities} with $\eta=0$. We notice that, under the same assumptions of Lemma~\ref{lem:ABC}, one can easily obtain the wanted inequality
\begin{equation}\label{inq.linearised.quantities.improved}
   \linI^{(0)}[v] \le \kappa_1(\ve)\,\linI^{(0)}_\gamma[v] + \kappa_2(\ve)\,\linE[v]\,,
 \end{equation}
 by taking the limit $\eta\to0$. Notice that in~\eqref{inq.linearised.quantities.improved} the constants $\kappa_1(\ve)$ and $\kappa_2(\ve)$ are as in~\eqref{kappas} and the constant $\mathcal{C}_{p, N}$, which enters into the definition of $\kappa_2$, can be easily computed from the proof of Lemma~\ref{lem:ABC} and its value is
 \[
   \mathcal{C}_{p,N}=\frac{N+p-\gamma-1}{(1-\gamma)^{2-p}}\,.
 \]
Lastly, we notice that the limit $\eta\to0$ could be easily justified by using the Monotone Convergence Theorem, since both $\linI^{(\eta)}[v]$ and $  \linI^{(\eta)}_\gamma[v]$ are monotone in $\eta$.

Let us now comment on the counterpart of inequality~\eqref{inq:linE-linI}, namely
\begin{equation}\label{inq:linE-linI.improved}
  \Lambda\,\linE[v]\leq \linI^{(0)}[v]\,.
\end{equation}
As it is clear from Lemma~\ref{lem:linE-leq-linI-eps} and Proposition~\ref{prop:hp}, inequality~\eqref{inq:linE-linI.improved} can be written as Hardy--Poincaré-type inequality of the form:
\begin{equation}\label{rewritten-HP}
 \frac{\Lambda}{2}\int_\rn |\phi|^2 \V_D^{2-\gamma} \dy \leq \frac{1}{|\gamma-1|^p}\int_\rn  \left|\nabla \phi \right|^2\,\VD\,|\nabla \V_D^{\gamma-1}|^{p-2} \dy\,,
\end{equation}
where $\phi=\V_D^{\gamma-2}(v- \V_D)$, under the addition hypothesis
\[
  \int_\rn \phi \V_D^{2-\gamma} \dy = 0 = \int_\rn (v-\VD)\dy=0\,.
\]
When $\phi$ is radial, and by using the change of variables $\frac{2-p}{p\,D}r^\frac{p}{p-1}=s^2$, inequality~\eqref{rewritten-HP} is equivalent
\begin{equation}\label{optimal-HP-radial}
   \Lambda_{\mathrm{opt}}\, \int_0^\infty g^2(s)\, \frac{s^\frac{2N(p-1)}{p}}{(1+s^2)^\frac{1}{2-p}}\,\frac{\ds}{s} \leq \int_0^\infty |g'(s)|^2\,\frac{s^{\frac{2N(p-1)}{p}}}{(1+s^2)^\frac{p-1}{2-p}}\,\frac{\ds}{s}\,,
\end{equation}
already introduced in~\eqref{optimal-HP.intro}, which holds under the additional assumption of
\[
  \int_0^\infty g(s)\frac{s^\frac{2N(p-1)}{p}}{(1+s^2)^\frac{1}{2-p}}\,\frac{\ds}{s}=0\,.
\]
Inequality~\eqref{optimal-HP-radial} has been investigated in~\cite{Bonforte2017a,Bonforte2017b,Bonforte_2023}, not only in the radial case, and the optimal value of $\Lambda_{\mathrm{opt}}$ is known. From the information contained in those papers by a simple, however a little tedious, computation, we learn that the optimal constant $\Lambda$ in~\eqref{inq:linE-linI.improved} is given by
\begin{equation}\label{value:Lambda:improved}
\Lambda= \begin{cases}
\frac{(p-1)}{(2-p)^2\beta}\quad&\mbox{when}\quad p_M<p<2\,,\\
\frac{\left[p-N(2-p)(p-1)\right]^2}{4p(2-p)^3}&\mbox{when}\quad p_c<p\leq p_M\,.
\end{cases}
\end{equation}
\begin{rem}
We notice that for $N\ge6$ it holds $\left[p-N(2-p)(p-1)\right]^2=\left(p-p_1\right)^2\left(p-p_2\right)^2$, where $p_1$ and $p_2$ are as in~\eqref{p1p2}.
\end{rem}
We are now in the position of proving Theorem~\ref{theo:Entropy-decay}.
\begin{proof}[Proof of Theorem~\texorpdfstring{~\ref{theo:Entropy-decay}}{2}]
In the general case $p_c<p<2$ without radiality assumptions, inequality~\eqref{entropy-decay} follows from the proofs of Proposition~\ref{prop:convergence.lp} and~\ref{convergence.lp.below.pm}. Let us consider the \emph{optimality} result in the radially decreasing case. Our first task is to prove the following claim. \\

\noindent\textbf{Claim~1}.
Assume that $0\le \vo\in C^2(\R^N) $ is radially symmetric and decreasing ($\partial_r \vo\leq 0$), satisfies~\eqref{stronger.GHP} and $\partial_r v_0$ satisfies~\eqref{hp}. Then inequality~\eqref{entropy-decay} holds true for any $\lambda\in\left(0, (p-1)\Lambda\right)$ where $\Lambda$ is as in~\eqref{value:Lambda:improved}. \\

Let us proceed with the proof of Claim~1. Under the current assumptions, by Theorem~\ref{theo:RECR} we have that $v(\tau)$ converges to $\VD$ uniformly in relative error. This implies that, for any $\ve\in\left(0,1\right)$ there exists $\tau_\ve'>0$ such that $v(\tau)$ satisfies~\ref{A1} (with the chosen $\ve$) for any $\tau\ge\tau_\ve'$.
Thanks to~\cite[Theorem 1.4]{Bonforte2020b} (see also our Theorem~\ref{Thm:gradient.convergence}), assumption~\eqref{hp.derivative} guarantees that $\partial_r v(\tau)$ converges to $\partial_r \VD$ uniformly in relative error. This implies that, for any $\ve\in\left(0,1\right)$ there exists $\tau_\varepsilon''>0$ such that $\partial_r v(\tau)$ satisfies inequality~\eqref{hp} for any $\tau\ge\tau_\varepsilon''$. Let us define $\tau_\ve=\max\{\tau_\ve', \tau_\ve''\}$. For any $\tau\ge\tau_\ve$, we can use both inequalities~\eqref{cal-I-est:improved} and~\eqref{inq.linearised.quantities.improved} to get
\begin{equation*}
  \linI^{(0)}[v(\tau)] \le \frac{1}{(p-1)}\frac{(1+\varepsilon)^{2-p}\kappa_1(\ve)}{(1-\varepsilon)^{\mathsf{a}}}\,\cI[v(\tau)|\VD] + \kappa_2(\ve)\,\linE[v(\tau)]\,,\quad\mbox{for any}\quad\tau\ge\tau_\ve\,.
\end{equation*}
where $\mathsf{a}=1+(1-\gamma)(2-p)$. From this inequality, by taking into account Lemma~\ref{lem:ent-lin-ent-2} and~\eqref{inq:linE-linI.improved}, we deduce that
\begin{equation}\label{inq:almost-optimal}
  \frac{\left(\Lambda-\kappa_2(\varepsilon)\right)}{(1-\varepsilon)^{\gamma-2}}\,\cE[v(\tau|\VD)]\leq \frac{1}{(p-1)}\frac{(1+\varepsilon)^{2-p}\kappa_1(\ve)}{(1-\varepsilon)^{\mathsf{a}}}\,\cI[v(\tau)|\VD]\,,\quad\mbox{for any}\quad\tau\ge\tau_\ve\,.
\end{equation}
Since $\kappa_1(\ve)\to(1-\gamma)^{-2}=\frac{(p-1)^2}{(2-p)^2}$ and $\kappa_2(\ve)\to0$ as $\ve\to0$ (see~\eqref{kappas}), we have that, for any $\lambda\in\left(0, \frac{(2-p)^2\Lambda}{(p-1)}\right)$, there exists $\tau_\lambda>0$ such that
\begin{equation*}
  \lambda\,\cE[v(\tau)|\VD]\leq \cI[v(\tau)|\VD]\,,\quad\mbox{for any}\quad\tau\ge\tau_\lambda\,.
\end{equation*}
By taking into account Lemmata~\ref{lem:ent/ent-prod} and~\ref{lem:ent/ent-prod-2}, and by using the Gromwall Lemma (as in Proposition~\ref{prop:convergence.lp} of~\ref{convergence.lp.below.pm}), we deduce that $\cE[v(\tau)|\VD]\leq \,C\, e^{-\lambda\,\tau}$, for any $\tau\ge\tau_\lambda$. This concludes the proof of the  claim. \\

We are now in the position to obtain the optimal rate $\lambda=\frac{(2-p)^2}{(p-1)}\Lambda$. The main idea is to take into account the dependence in time of $\varepsilon$ in inequality~\eqref{inq:almost-optimal}. Let us recall that $\varepsilon$ can be defined as
\[
  \varepsilon(\tau):=\max\left\{ \left|\frac{v(\tau)-\VD}{\VD}\right|\,, \left|\frac{\partial_r v(\tau) - \partial_r \VD}{\partial_r \VD}\right|\right\}\,,
\]
and that, by Theorems~\ref{theo:RECR} and~\ref{Thm:gradient.convergence}, there exists $\tilde{\sigma}>0$ and $\tau_{\tilde{\sigma}}>0$, such that $\varepsilon(\tau)\leq  e^{-\tilde{\sigma}\tau}$ for any $\tau\geq \tau_{\tilde{\sigma}}$. Inequality~\eqref{inq:almost-optimal} can be written as
\begin{equation*}
  \frac{(2-p)^2\Lambda}{(p-1)} \left(R(\varepsilon(\tau)) + 1\right)\,\cE[v(\tau)|\VD]   \leq \cI[v(\tau)|\VD]\,,\quad\mbox{for any}\quad\tau\ge\tau_{\tilde{\sigma}}\,,
\end{equation*}
where
\[
R(\varepsilon(\tau))=\frac{\Lambda -\kappa_2(\varepsilon(\tau))}{\Lambda}\frac{(p-1)^2}{(2-p)^2\kappa_1(\varepsilon(\tau))}\,\frac{(1-\varepsilon(\tau))^{\mathsf{a}+2-\gamma}}{(1+\varepsilon(\tau))^{2-p}}-1\,.
\]
By the Taylor expansion, from~\eqref{kappas} we deduce the existence of a constant $\mathsf{C}>0$, such that $|R(\varepsilon(\tau))|\leq\mathsf{C}\varepsilon(\tau)$ for any $\tau\geq\tau_{\tilde{\sigma}}$. Let us now consider the auxiliary function $Z(\tau):=\log(\cE[v(\tau)|\VD])$. It satisfies the following differential inequality
\[
  Z'(\tau)=-\frac{\cI[v(\tau)|\VD]}{\cE[v(\tau)|\VD]}\leq - \frac{(2-p)^2\Lambda}{(p-1)} + \frac{(2-p)^2\Lambda}{(p-1)} R(\varepsilon(\tau))\,.
\]
By applying Lemma~\ref{lem:improved-rate} to the function $\tau\mapsto Z(\tau)$, we retrieve the wanted inequality. The proof is concluded.
\end{proof}

\section{Proofs of Theorem~\texorpdfstring{\ref{theo:meta1}}{4} and~\texorpdfstring{\ref{theo:meta}}{5} capturing in particular the subcritical case~\texorpdfstring{$p\leq p_c$}{p<=pc}}\label{sec:comm}

\begin{proof}[Proof of Theorem~\ref{theo:meta1}]Through the proof we shall assume condition~\eqref{integrability.condition} and shall not distinguish the cases $p\leq p_c$ or $p_c<p<2$. One reasoning works in both cases. We shall follow several steps from~\cite{Agueh2009}, which relies on the ideas of~\cite{Blanchet2009}.
\newline

\noindent\textit{Step 1): identification of the limit when $\tau\rightarrow \infty$.}
For this step, we follow mainly \cite[Lemma 2.5]{Agueh2008}. We will first prove that $v(\tau, \cdot)$ converges to $V_D$  pointwise and in $\LL^p$-norms.  As in~\cite[Lemma 2.5]{Agueh2008}, let us define $v^h(\tau, y):=v(\tau+h, y)$, for any given $h>0$ and $\tau\in\left[0, 1\right]$.  By the comparison principle and thanks to assumption {\it (i)}, $\{v^h\}$ is uniformly bounded. Furthermore, it is uniformly continuous in $[0, 1]\times B_R$ thanks to assumption {\it (iv)}. By the Ascoli--Arzel\'a Theorem, for any sequence $h_n\rightarrow\infty$ (as $n\rightarrow\infty$) the sequence of functions $\{v^{h_n}\}$ converges uniformly (up to a subsequence) to a function $v^{\infty}$ on compact subsets of $[0, 1]\times \rn$. Moreover, we infer  that for any $R>0$ it holds that $\|v\|_{C^{1,\alpha}([0, 1]\times B_R)}<\infty$ and for $\tau\in\left[0, 1\right]$ function $v^\infty(\tau, \cdot)$ satisfies {\it (i)}. Since $N<\frac{p}{(2-p)(p-1)}$, we know that $\V_{D}^{\gamma-2}(v-\V_{D})^2\in \LL^1(\rn)$, cf.~\eqref{eq1} and the end of the proof of Lemma~\ref{lem:I-est}. Therefore, by using the arguments of Lemma~\ref{lem:ent-lin-ent-2}, we get that $\linE[v], \cE[v |\VD] <\infty $. Thanks to assumption {\it (iii)} the entropy functional $\cE[v |\VD]$ is non-negative and $\tau\mapsto\cE[v(\tau) |\VD]$ is decreasing in time. By the time monotonicity,  $\cE[v^{h_n}(\tau)|\VD]$ and $\cE[v^{h_n+1}(\tau)|\VD]$ have the same limit for $h_n\rightarrow\infty$. Therefore we infer that
\[
\int_0^1 \cI[v^{h_n}(\tau)|\VD] \dtau  = \int_{h_n}^{h_n+1}   \cI[v(\tau)|\VD] \dtau = \cE[v(h_n)|\VD] - \cE[v(h_n+1)|\VD] \xrightarrow[n\rightarrow\infty]{} 0\,.
\]
 By the positivity of $\cI[v^\infty(\tau)|\VD]$ and Fatou's Lemma, we infer therefore that $0\geq \int_0^1 \cI[v^\infty(\tau)|\VD] \dtau =0 $. Consequently, $\nabla (v^{\infty})^{\gamma-1}=\nabla V_{D^\star}^{\gamma-1}$ for some $D^\star>0$ and so $v^\infty=V_{D^\star}$.  Up to now, we have proven that $v^{h_n}$ converges pointwise towards $V_{D^\star}$ as $n\rightarrow\infty$. We only need to ensure that $D^\star=D$. By Lemma~\ref{decay.lemma}, we are in the position of using the Dominated Convergence Theorem to infer that $\{(v^{h_n}-V_{D^\star})\}$ converges in $\LL^1(\rn)$, which implies that
\[
\lim_{n\rightarrow \infty} \int_{\rn} \left(v^{h_n}(y) - V_{D^\star}(y) \right) \dy =0\,.
\]
At the same time, the above identity implies that necessarily $D^\star=D$. Indeed, otherwise  one would find that $\int_{\rn}(V_{D^\star}(y)-V_{D}(y))\dy =0$, which leads to a contradiction in the case $D\neq D^\star$.  Lastly, we observe that the limit does not depend on the sequence $\{v^{h_n}\}$ since the above reasoning is true for any possible convergent subsequence. Therefore, we conclude that $(v-V_D)$ converges to zero in the $\LL^1$-topology as $\tau\to\infty$.\newline

\noindent\textit{Step 2): from convergence in $\LL^1(\rn)$ to convergence in $\LL^\infty$.} In this step we follow the ideas of~\cite[Lemma~2.6]{Agueh2008} with a~few differences to be stressed. We notice that, by assumption {\it (i)}, the $\LL^\infty$-norm of the  function $y\mapsto |v(\tau, y)-V_D(\tau, y)|$ is bounded uniformly in $\tau$, see again Lemma~\ref{decay.lemma}. Therefore, by interpolation, one obtain that $v(\tau)$ converges to $V_D$ in the $\LL^q(\rn)$-topology, for any $1\leq q < \infty$. The convergence in $\LL^\infty(\rn)$ is more subtle. We shall first prove this convergence on balls. Let $R>0$ and $\mathfrak{d}\in(0,1)$. For any function $f\in C^{\mathfrak{d}}(B_{2R})\cap \LL^1(B_{2R})$ we have the following interpolation inequality whose proof can be found in~\cite{Bonforte2020a}:
\begin{equation}\label{interpolation.domains}
  \|f\|_{\LL^\infty(B_R)} \leq C_{N, \mathfrak{d}}\left(\|f\|_{C^\mathfrak{d}(B_{2R})}^\frac{N}{N+\mathfrak{d}}\,\|f\|_{\LL^1(B_{2R})}^\frac{\mathfrak{d}}{N+\mathfrak{d}}+ R^{-N}\|f\|_{\LL^1(B_{2R})}\right)\,.
\end{equation}
Let us fix $\varepsilon>0$. By using~\eqref{interpolation.domains}, assumption {\it (iv)}, and the already proven $\LL^1$-convergence, we infer that there exists $\tilde{\tau}=\tilde{\tau}(\varepsilon, \vo)>0$ such that
\begin{equation}\label{convergence.compact}
  \|v(\tau)-V_D\|_{\LL^\infty(B_R)}<\varepsilon\quad\,\forall\,\tau>\tilde{\tau}\,.
\end{equation}
Assume further that $R>C^{-1}\,\varepsilon^\frac{-(p-1)(2-p)}{p}$ where $C$ is as in Lemma~\ref{decay.lemma}. Thanks to assumption {\it (i)} and by Lemma~\ref{decay.lemma}, we infer that
\begin{equation}\label{tail.control}
  \left|v(\tau, y) - V_D(y) \right|\leq \varepsilon \quad \forall\,\, |y|\geq R\,.
\end{equation}
Inequalities~\eqref{convergence.compact} and~\eqref{tail.control} imply that $v$ converges to $V_D$ in $\LL^\infty(\rn)$ as $\tau\to\infty$.\newline

\noindent\textit{Step 3): from convergence in $\LL^\infty(\rn)$ to convergence in the uniform relative error. } We will prove first the convergence of the relative error in the $\LL^\infty$-norm and then obtain the general result by an interpolation argument. By using again Lemma~\ref{decay.lemma}, one can infer that the relative error decays as follows
\begin{equation}\label{tail.control.re}
  \left|\frac{v(\tau, y)-V_D(y)}{V_D(y)} \right|\leq \kappa_1 \frac{|y|^\frac{p}{2-p}}{|y|^\frac{p}{(p-1)(2-p)}}=\frac{\kappa_1}{|x|^\frac{p}{p-1}}\qquad\forall\,\, |y|\ge 1\,.
\end{equation}
On a ball of radius $R>0$, from the $\LL^\infty(\rn)$-convergence, one can infer the following:
\begin{equation}\label{compact.re}
  \sup_{|y|\leq R} \left|\frac{v(\tau, y)-V_D(y)}{V_D(y)} \right|\leq \kappa_2\, \|v(\tau, y)-V_D(y)\|_{\LL^\infty(\rn)}\, R^\frac{p}{2-p}\,.
\end{equation}
For fixed $\varepsilon>0$ there exist $R_\varepsilon>0$ and $\tau_\star=\tau_\star(\varepsilon, \vo)>0$ such that \[\text{$\kappa_1\, R_\varepsilon^\frac{-p}{p-1}<\varepsilon\quad$ and $\quad\kappa_2\, \|v(\tau, y)-V_D(y)\|_{\LL^\infty(\rn)}\, R_\varepsilon^\frac{p}{2-p}<\varepsilon$}\,.\]  Combining together inequalities~\eqref{tail.control.re} and~\eqref{compact.re} we get that for any $\varepsilon>0$ it holds
\begin{equation}\label{convergence.re.infty}
  \left\| \frac{v(\tau, y)-V_D(y)}{V_D(y)}\right\|_{\LL^\infty(\rn)}\leq \varepsilon\quad\forall\,\tau\ge \tau_\star\,,
\end{equation}
which justifies the uniform convergence in the relative error. It only remains to prove the convergence of the relative error in $\LL^q(\rn)$. Notice that the relative error $\frac{v(\tau,\cdot) -V_D(\cdot)}{V_D(\cdot)}$ is uniformly bounded in space and, thanks to inequality~\eqref{tail.control.re}, it is integrable for any $q> N\,\frac{(p-1)}{p}$. Indeed, for $\delta>0$ such that $2\,\delta< q - N\,\frac{(p-1)}{p}$ we have the following inequality
\begin{equation*}
  \left\|\frac{v(\tau, y)-V_D(y)}{V_D(y)}\right\|_{\LL^q(\rn)}\leq \left\|\frac{v(\tau, y)-V_D(y)}{V_D(y)}\right\|_{\LL^\infty(\rn)}^{q-N\frac{p-1}{p}-\frac{\delta}{2}}\, \int_{\rn}\left|\frac{v(\tau, y)-V_D(y)}{V_D(y)} \right|^{N\frac{(p-1)}{p}+ \frac{\delta}{2}} \dy \xrightarrow[\tau\to\infty]{}0 \,,
\end{equation*}
where we used the fact that $\int_{\rn}\left|\frac{v(\tau, y)-V_D(y)}{V_D(y)} \right|^{N\frac{(p-1)}{p}+ \frac{\delta}{2}} \dy $ is uniformly bounded in time. The proof is  complete.
\end{proof}

\begin{proof}[Proof of Theorem~\ref{theo:meta}] The strategy of the proof is to obtain first a convergence rate of the convergence in the $\LL^1$-topology and then improve it to the final result~\eqref{fin}. We remind that we stay under condition~\eqref{integrability.condition} and that the assumptions of Theorem~\ref{theo:meta1} are satisfied. From the convergence result~\eqref{fin1}, and assumption {\it (i)}, we deduce that there exists $\tau_\bullet=\tau_\bullet(\vo,D_1, D_2)>0$ such that \ref{A0}, \ref{A1} (for some $\varepsilon>0$), and \ref{A2} hold for $v(\tau,y)$ with every $\tau>\tau_\bullet$ and $y\in\rn$.  Since we assume the decay condition \eqref{gradient.assumption}, we can make use of the lines of the proof of Lemma~\ref{lem:ent-lin-ent} to justify that $\cE[v(\tau)|\VD]<\infty$.  Since $N<\frac{p}{(2-p)(p-1)}$ we know that $\V_{D}^{\gamma-2}(v-\V_{D})^2\in \LL^1(\rn)$. Analogously, using the arguments of Lemma~\ref{lem:ent-lin-ent-2}, we get that $\linE[v], \cE[v |\V_{D}] <\infty $, and $\left(1+\varepsilon\right)^{\gamma-2} \linE[v] \le \cE[v |\V_{D}] \le \left(1-\varepsilon\right)^{\gamma-2} \linE[v]\,$.

Let us now clarify the relation between the entropy and the Fisher information, both nonlinear and linearised versions. We notice that under assumptions \ref{A1}, \ref{A2} and {\it (ii)} we are able to repeat the proof of Lemma~\ref{lem:linE-leq-linI-eps} that implies that $\mathcal{C}(p,D,\ve)\linE[v]\leq \linI^{(\ve)}[v]$. At the same time, we get that $\cI[v(\tau)|\V_{D}]\geq C_\varepsilon\, \linI_\gamma^{(\varepsilon)} [v(\tau)]$ via Lemma~\ref{lem:ABC} and the reasoning of Lemma~\ref{lem:I-est}, where we make use of~\eqref{gradient.assumption} in the place of~\eqref{space.gradient.decay}. Collecting the above we infer that $\frac{d}{d\tau}\cE[v(\tau)|\VD]\leq -c \cE[v(\tau)|\VD],$ which via the Gronwall Lemma allow to state that for all $\tau>0$ it holds $\cE[v(\tau)|\VD]\leq e^{-\vt\, \tau} \cE[\vo|\VD]$ for some $\mathfrak{K}=\mathfrak{K}(p,N,\ve)$. On the other hand, due to {\it (i)} and {\it (ii)} and the Csisz\'ar--Kullback inequality provided in Lemma~\ref{CK.inq}, we know that $\|v-\V_{D}\|^2_{\LL^1(\rn)}
\leq c(\V_{D},p)\cE [v(\tau)|\V_{D}]$, where the right-hand side is finite. Therefore, we get that there exists $\wt{T}>0$ and $\wt{K}>0$, such that  we have that\begin{equation*}
\|v(\tau,\cdot)-\V_{D}(\cdot)\|_{\LL^{1}(\rn)}\le \wt{K}\,e^{-\mathfrak{K}\tau/2}\qquad\forall\, \tau\ge \wt{T}\,.
\end{equation*}

We can now get a convergence rate in the uniform relative error. Let $\varepsilon>0$ and $R=\left(\frac{\kappa_1}{\ve}\right)^\frac{(p-1)}{p}$ where $\kappa_1$ is as in~\eqref{tail.control.re}. Then we obtain from inequality~\eqref{tail.control.re} that
\[
 \left\|\frac{v(\tau, y)-V_D(y)}{V_D(y)}\right\|_{\LL^q(\rn)}\leq \ve\quad\forall \,|y|>R\,\,\forall\,\,\tau>0\,.
\]
At the same time, by using an interpolation inequality between the $\LL^\infty$, $C^1$, and $\LL^1$ norms on $\rn$ (directly resulting from~\eqref{interpolation.domains} by taking the limit $R\rightarrow\infty$) one finds that for any  $\tau\geq \tau_0$ (where $\tau_0$ is as in {\it (iv)}) we have
\[
\|v(\tau)-V_D\|_{\LL^\infty(\rn)} \leq C_{N} \|v(\tau, \cdot)-V_D(\cdot)\|_{C^1(\rn)}^\frac{N}{N+1} \|v(\tau, \cdot)-V_D(\cdot)\|_{\LL^1}^\frac{1}{N+1} \leq C_{N}\,C(\alpha, \tau_0, \vo)\, \wt{K}\,e^{-\mathfrak{K}\tau/2}\,,
\]
where $C(\alpha, \tau_0, \vo)$ is as in {\it (iv)}.
Combining the above estimate with~\eqref{compact.re} one get~\eqref{convergence.re.infty} with $\tau_\star=\max\{\tau_0, - \frac{2}{\mathfrak{K}}\log(H\, \varepsilon^\frac{3-2p}{2-p}) \}$, for a constant $H=H(\alpha, \tau_0, \vo)>0$. We notice that for $\ve$ small enough we have that $ - \frac{2}{\mathfrak{K}}\log(H\, \varepsilon^\frac{3-2p}{2-p})>\max\{0,\tau_0\}$. Once inequality~\eqref{convergence.re.infty} is obtained with an explicit functional relation between $\ve$ and $\tau_\star(\ve)$, one can compute the rate of convergence by inverting this relation, as it was done, for instance, in~\cite[Corollary 4.14]{Bonforte2020a}. This is enough to obtain the convergence result~\eqref{fin} for the $\LL^\infty$-norm. The result in the $\LL^q$-norm is obtained by interpolation as it is done in the proof of Theorem~\ref{theo:meta1}. The proof is complete.
\end{proof}

\section{Justification of the assumption for Theorems~\ref{theo:meta1} and~\ref{theo:meta} in the radial decreasing case}\label{sec:jus_radial}

In this section we shall prove that, for an initial datum $\vo$ that is radially decreasing (and satisfy an additional hypothesis), the assumptions of Theorems~\ref{theo:meta1} and~\ref{theo:meta} hold. Let us begin with the following claim.

\smallskip
\noindent\textbf{Claim 2.} Let $v$ be solution to~\eqref{eq:v} with an initial datum $\vo$ that satisfies assumption \textit{i)} of Theorem~\ref{theo:meta1}, i.e. inequality~\eqref{stronger.GHP.fast.range}.  If that intial datum $\vo$ satisfies
\begin{equation}\label{additional.assumption}
  \partial_r \V_{D_2}(r)\leq \partial_r \vo(r) \leq \partial_r \V_{D_1}(r)\quad\forall r\ge0\,,
\end{equation}
then
\begin{equation}\label{additional.consequence}
  \partial_r \V_{D_2}(r)\leq \partial_r v(\tau, r) \leq \partial_r \V_{D_1}(r)\quad\quad \forall \tau>0\quad\mbox{and}\quad\forall r\ge0\,.
\end{equation}

We notice that, when $p_Y\leq p<2$, the above claim is a consequence of the theory developed in Section~\ref{sec:conv-radial}. Since $p_Y\leq p_2$ when $N\ge6$, it only remains to prove the claim when $1<p<p_1$. We shall see below that this will be a consequence of the standard theory for weighted parabolic equation of type~\eqref{ckn.fde.radial}. For the moment, we shall focus on providing a proof for the other assumptions of Theorem~\ref{theo:meta1}.\newline

\noindent\textit{Proof that condition $(ii)$ of Theorem~\ref{theo:meta1} holds under the additional hypothesis~\eqref{additional.assumption}}. Let $D>0$, and consider $\VD$ (for the moment not necessarily such that $\int_{\mathbb{R}^N}(\vo-\VD)\dx=0$). By using the definition of weak solution (straightforwardly adapted to equation~\eqref{eq:v}) we find that, for any $s>t\geq0$ and any radial, compactly supported, smooth function $\phi:\mathbb{R}^N\to\left[0, \infty\right)$, we have
\begin{multline}\label{first.formula}
\int_{\RN}\left(v(s,y)-\VD(y)\right)\phi(y)\dy - \int_{\RN}\left(v(t,y)-\VD(y)\right)\phi(y)\dy  \\
= \int_s^t\int_{\RN}\big( |\nabla v(\tau,y)|^{p-2}\nabla v(\tau,y) - |\nabla \VD(y)|^{p-2}\nabla \VD(y)\big) \cdot \nabla \phi(y) \dx\dtau\,.
\end{multline}
Since both $v$ and $\VD$ are radial, we can rewrite the absolute value of right-hand-side term of~\eqref{first.formula} as
\[
 \left| \big( |\nabla v(\tau,y)|^{p-2}\nabla v(\tau,y) - |\nabla \VD(y)|^{p-2}\nabla \VD(y)\big) \cdot \nabla \phi(y) \right|  = \left||\partial_r v|^{p-1}(r)- |\partial_r \VD|^{p-1}(r)\right| |\partial_r \phi(r)|\,,
\]
where $r=|y|$. Let us consider $\psi\mathbb{R}^N\to\left[0, \infty\right)$ a radial cut-off function $\psi(y)=\psi(r)$ which is equal to $1$ when $0\leq r \leq 1$ and equal to $0$ when $r\ge2$. Let us take $\phi(y)=\psi(y/R)$ in~\eqref{first.formula}.  Since $|\nabla \psi|\leq c\, R^{-1}$, where $c>0$ is a dimensionless constant. Since $\partial_r \psi(r)\neq 0$ only when $R<r<2R$ we find, by applying Lemma~\ref{decay.lemma}, that
\[
  \left||\partial_r v|^{p-1}- |\partial_r \VD|^{p-1}(r)\right| |\partial_r \phi(r)| \leq C r^{-\frac{p}{(p-1)(2-p)}}\,,\quad\forall\,\, r\ge0\,.
\]
Since the power $r\mapsto r^{-\frac{p}{(p-1)(2-p)}}$ is integrable to infinity, we find that
\[
  \lim_{R\rightarrow\infty}\int_{\RN}\big( |\nabla v(\tau,y)|^{p-2}\nabla v(\tau,y) - |\nabla \VD(y)|^{p-2}\nabla \VD(y)\big) \cdot \nabla \phi(y) \dy=0\,,
\]
from which we deduce that, for any $s>t\geq0$
\[
  \int_{\RN}\left(v(s,y)-\VD(y)\right)\dy = \int_{\RN}\left(v(t,y)-\VD(y)\right)\dy\,.
\]
Since $\vo$ satisfies assumption \textit{i)} of Theorem~\ref{theo:meta1}, and in this range the difference of two profiles $\V_{D_2}-\V_{D_1}$ is integrable, there must exists $D>0$ such that $\int_{\mathbb{R}^N}(\vo-\VD)\dy=0$. The proof is concluded.\newline

\noindent\textit{Proof that condition $(iii)$ of Theorem~\ref{theo:meta1} holds under the additional hypothesis~\eqref{additional.assumption}. }
We shall prove that the Fisher information $\tau\to\cI[v(\tau)|\VD]$ is a $\LL^\infty_{\textrm{loc}}(0, \infty)$-function. In order to do so, we notice that, under the assumptions of Claim~2, the assumptions of Proposition~\ref{prop:cal-I-est:improved} are satisfied. Therefore, it is only needed to prove that $\tau\to\linI^{(0)}_\gamma[v(\tau)]$ is a locally (in time) finite function. By the proof of Proposition\ref{prop:cal-I-est:improved}, we have that
\[
  \linI^{(0)}_\gamma[v(\tau)] \leq C\, \int_{\RN} \VD|\nabla \V_D^{\gamma-1}(y)|^{p-2}\,|\nabla v^{\gamma-1}(\tau, y) - \nabla \V_D^{\gamma-1}(y)|^2 \dy
\]
for some constant $C>0$. We shall prove that the function $y\to \VD(y)|\nabla \V_D^{\gamma-1}(y)|^{p-2}\,|\nabla v^{\gamma-1}(\tau, y) - \nabla \V_D^{\gamma-1}(y)|^2$ is an integrable function. First of all, we notice that, since  $\nabla \V_{D}^{\gamma-1} = \partial_r \V_{D}^{\gamma-1}=(\gamma-1)\V_{D}^{\gamma-2}\partial_r \VD$ and $\partial_r \VD= r^\frac{1}{p-1}\V_{D}^\frac{1}{p-1}$, we have that
\begin{equation}\label{first.information}
\VD |\nabla \V_D^{\gamma-1}|^{p-2}=|1-\gamma|^{p-2}\,\V_{D}^{1+(\gamma-2)(p-2)+\frac{p-2}{p-1}} r^\frac{p-2}{p-1}= |1-\gamma|^{p-2}r^\frac{p-2}{p-1}\,\VD\,,
\end{equation}
since $1+(\gamma-2)(p-2)+\frac{p-2}{p-1}=1$. Under the current assumptions, in a neighbourhood of the origin, we can bound the term $|\nabla v^{\gamma-1}(\tau,y) - \nabla \V_D^{\gamma-1}(y)|^2$ as
\begin{equation}\label{second.information}
  |\nabla v^{\gamma-1}(\tau,y) - \nabla \V_D^{\gamma-1}(y)|^2 \leq C\, r^\frac{2}{p-1}\,,
\end{equation}
for a constant $C>0$. By taking into account~\eqref{first.information} with~\eqref{second.information}, we find the bound
\[
\VD |\nabla \V_D^{\gamma-1}|^{p-2}|\nabla v^{\gamma-1}(\tau,y) - \nabla \V_D^{\gamma-1}(y)|^2 \leq C\,r^\frac{p}{p-1}\,,
\]
which proves that the function $y\to \VD(y)|\nabla \V_D^{\gamma-1}(y)|^{p-2}\,|\nabla v^{\gamma-1}(\tau, y) - \nabla \V_D^{\gamma-1}(y)|^2$ is integrable close to the origin.
Let us tackle the issue of integrability at infinity. By elementary computation, we find that
\begin{equation}\label{third.information}
  |\nabla v^{\gamma-1} - \nabla \V_D^{\gamma-1}|^2\leq 2\,(\gamma-1)^2\, \left(|\partial_r v(\tau, r)|^2|v(\tau, r)^{\gamma-2}-\V_D^{\gamma-2}(r)|^2 + \V_D^{2(\gamma-2)}|\partial_r v(\tau, r) - \partial_r \VD|^2\right) \,.
\end{equation}
Let us consider the first term in the right-hand-side of~\eqref{third.information}. By inequality~\eqref{additional.consequence}, we deduce that there exists a constant $c_1>0$ such that $|\partial_r v(\tau, r)|^2\leq c_1 r^{-\frac{4}{2-p}}$. By using the Mean Value Theorem applied to the function $\xi\mapsto \xi^{\gamma-2}$ and by using the fact that $v(\tau)\leq \V_{D_2}$ and inequality~\eqref{better-decay} of Lemma~\ref{decay.lemma}, we deduce that there exists a constant $c_2>0$ such that
\begin{equation}\label{fourth.information}
  |v(\tau, r)^{\gamma-2}-\V_D^{\gamma-2}(r)|^2\leq (\gamma-2)^2\, |\V_{D_2}|^{2(\gamma-3)}\,|v(\tau, r)-\V_D(r)|^2\leq c_2\,r^\frac{2p}{2-p}\,.
\end{equation}
By considering identity~\eqref{first.information}, combined with inequality $|\partial_r v(\tau, r)|^2\leq c_1 r^{-\frac{4}{2-p}}$ and inequality~\eqref{fourth.information} we obtain
\[
  \VD |\nabla \V_D^{\gamma-1}|^{p-2}\,|\partial_r v(\tau, r)|^2\,|v(\tau, r)^{\gamma-2}-\V_D^{\gamma-2}(r)|^2 \leq c_3\, r^{\frac{p-2}{p-1}- \frac{p}{2-p}-\frac{4}{2-p}+\frac{2p}{2-p}}=c_3\,r^{-\frac{p}{(p-1)(2-p)}}\,,
\]
for a constant $c_3>0$ and $r$ large enough. We recall that, under the current assumptions the function $r\mapsto r^\frac{-p}{(p-1)(2-p)}$ is integrable at infinity. Lastly, we notice that, by taking into account point \textit{ii} of~\ref{decay.lemma}, we can estimate $\V_D^{2(\gamma-2)}|\partial_r v(\tau, r) - \partial_r \VD|^2$ in a similar manner and find that
\[
  \VD |\nabla \V_D^{\gamma-1}|^{p-2}\,\V_D^{2(\gamma-2)}|\partial_r v(\tau, r) - \partial_r \VD|^2\,\leq c_4 \,r^{-\frac{p}{(p-1)(2-p)}}\,,
\]
for a constant $c_4>0$ and $r$ large enough. The last two estimates prove that the function  $\tau\mapsto\linI^{(0)}_\gamma[v(\tau)]$ is $\LL^\infty_{\textrm{loc}}(0, \infty)$; by Proposition~\ref{prop:cal-I-est:improved} the same applies to $\cI[v(\tau)|\VD]$, therefore the proof of this point is complete.

\smallskip

\smallskip
\noindent\textit{Using radial equivalence between PLE and FDE. } We begin by proving Claim~2, since it will be used in what follows. We first recall the relation between radial derivatives of solutions $u$ to~\eqref{CPLE} (as function of $(r=|x|, t)$) and radial solutions $\Phi$ to a (weighted version of) the FDE, established in~\cite{Iagar2008}, and discussed above, see Section \ref{sec:conv-radial}.  Recall that, a radial solution $u$ to~\eqref{CPLE} is also a solution to the following equation
\begin{equation*}
\partial_t u = r^{1-N} \, {\partial_r} \left(r^{N-1}\,|{\partial_r u}|^{p-2}\,\partial_r u \right)\,.
\end{equation*}
By idenity~\eqref{radial.transformation}, i.e. $-\partial_r u(t,r)= {\mathcal{D}}\,\vr^\frac{2}{m+1}\,\Phi(t, \vr)$, we have that $\Phi: \mathbb{R}^N\times(0, \infty)\rightarrow \mathbb{R}$) is a nonnegative radial function of the variables $(\vr=|x|, t)$ which solves, by Theorem~\ref{radial.transformation.thm}, the following weighted equation
\begin{equation}\label{Phi-local}
\partial_t \Phi= \vr^{1-n} \, \partial_{\vr}\left(\vr^{n-1}\, \partial_{\vr} \Phi^m\right)\,,\qquad m=p-1\,.
\end{equation}

We notice that Claim 2 follows essentially from the \emph{comparison principle} (i.e. property~\eqref{comparison-principle}) for derivatives of equation~\eqref{eq:v}. As explained in Section~\ref{sec:conv-radial}, radial derivatives of a solution to~\eqref{eq:v} (up to a power of $\rho$) satisfy equation~\eqref{rescaled.ckn.fde.radial}, which is a rescaled version of~\eqref{Phi-local}. We concluded therefore that, in order to prove Claim 2, it suffices to prove a comparison principle for equation~\eqref{Phi-local}. There are two cases to be considered: the simpler one is when $n=2+2N\frac{(p-1)}{p}$ (which is always bigger than 2) is an integer, and in this case $\Phi$ can be seen as the radial solution to a FDE with $m=p-1$ with the classical Laplacian in dimension $n$, and the claim follows by \cite[Theorem 4]{Blanchet2009}. However, this does not happen in general, hence we introduce an operator with suitable CKN-type weights: consider $\Phi$ as a radial function of $R^N$, then we have, according to \eqref{parameter.transformation}, that
\[
n= 2+2\tfrac{N}{p'}=N-\aa>2 \qquad\mbox{where}\qquad \aa = N-2-2\tfrac{N}{p'}=(N+2)\frac{p_Y-p}{p}\,.
\]
Notice that we always have $n>2$, while $\aa$ has no sign when $p\in (1,2)$ and it is zero when $p=p_Y=\frac{2N}{N+2}$, as remarked in Section \ref{ssec:information_on_related_radial_classical_and_weighted_fast_diffusion_equations}. Of course, the proof of this case (which we perform below) covers also the case when $n$ is integer.

\smallskip
\noindent\textit{Proof of Claim 2: }We shall deal with a slightly more general case: consider the operator
\[
\mathcal{L}_\aa\Phi(x)= |x|^\aa\,\dv\left(|x|^{-\aa}\,\nabla \Phi \right)= \Delta\Phi(x)  + |x|^\aa\,\nabla\left(|x|^{-\aa}\right)\cdot \nabla \Phi\,,
\]
with  $\aa= N-2-2\tfrac{N}{p'}$,  
 which is well defined for  $\Phi\in C^\infty_c(\RR^N\setminus\{0\})$ in which case the second equality above holds true.

The comparison principle for solutions to the Cauchy problem for the associated FDE
\begin{equation}\label{FDE-CKN-aa}
\partial_t \Phi= |x|^\aa\,\dv\left(|x|^{-\aa}\,\nabla \Phi^m\right)\,.
\end{equation}
essentially follows by a Kato type inequality for the operator $\mathcal{L}_\aa$. In the case of the standard laplacian $\Delta$,  the Kato inequality states the following: whenever $\Delta f \in \LL^1_{\rm loc}$, the following inequality holds true in the sense of distributions:
\begin{equation}\label{Kato.Lapl}
\Delta (f)_+ - H(f)\Delta f \ge 0\,.
\end{equation}
where $(f)_+$ denotes the positive part of $f$ and $H$ is the Heaviside function. This implies the Kato inequality for $\mathcal{L}_\aa$\,:
\[\begin{split}
\mathcal{L}_\aa (f)_+ - H(f)\mathcal{L}_\aa f
&= \Delta (f)_+ + |x|^\aa\,\nabla\left(|x|^{-\aa}\right)\cdot \nabla (f)_+
- H(f)\Delta f - H(f) |x|^\aa\,\nabla\left(|x|^{-\aa}\right)\cdot \nabla f
=\Delta (f)_+ - H(f)\Delta f \ge 0
\end{split}
\]
where in the second equality we have used that $\nabla (f)_+=H(f)\nabla f$ a.e. and in the last inequality we have just used \eqref{Kato.Lapl}. All the above inequalities are intended in distributional sense.

We are now going to show how Kato inequality implies comparison, more precisely that given any two solutions $\Phi_1(t,\cdot), \Phi_2(t,\cdot)\in L^1(|x|^{-\aa}\dx)$ to \eqref{FDE-CKN-aa}, we have
\begin{equation}\label{comp.Phi.01}
\int (\Phi_1(t,x)-\Phi_2(t,x))_+ \psi(|x|) |x|^{-\aa}\dx
\le \int (\Phi_1(0,x)-\Phi_2(0, x))_+ \psi(|x|) |x|^{-\aa}\dx\qquad\mbox{for all }t>0\,.
\end{equation}
for a suitable positive and radially decreasing $\psi$, such that $\mathcal{L}_\aa \psi\le 0$. We shall see an example of it below. The above inequality easily implies that if $\Phi_1(0,\cdot) \le \Phi_2(0,\cdot)$ then $\Phi_1(t,\cdot) \le \Phi_2(t,\cdot)$ for all $t>0$. Claim 2 follows from this last step, as explained above.

Let us show an example of $\psi$ that satisfies the above conditions: the radial expression of $\mathcal{L}_\aa$ is:
\[
\mathcal{L}_\aa\psi(r)=\psi''(r)+\frac{n-1}{r}\psi'(r)\qquad\mbox{recalling that }\qquad n= 2+2\tfrac{N}{p'}=N-\aa>2\,,
\]
hence
\[
\psi(r)=(1+r^2)^{-\frac{n-2}{2}}\qquad\mbox{is such that}\qquad \mathcal{L}_\aa\psi =-n(n-2)\psi^{\frac{n+2}{n-2}}\le 0
\]
We present now a proof of~\eqref{comp.Phi.01} that holds for strong solutions\footnote{
 For mild solutions, i.e. nonlinear gradient flows on the Banach space $L^1(|x|^{-\aa}\dx)$, this property is often called $T$-Contraction or ``well ordering'' of the nonlinear semigroup. This always holds for mild solutions, see for instance \cite{Vazquez2007} and references therein. Strong solutions are particular cases of mild solutions, hence the result holds. We have decided to sketch the proof, since it gives an idea of how it works for weak solution: this is done by a careful choice of test function in the weak formulation: take a smooth, compactly supported $\varphi_{i,j,k}(\tilde t,x)=\eta_i(\tilde t)\psi_k(|x|)\zeta_j(\tilde t,x)$, where $\eta_i(\tilde t)\xrightarrow[]{i\to\infty}\chi_[0,t](\tilde t)$ and $\psi_k\xrightarrow[]{k\to\infty}\psi$ and with some additional properties, cf. Appendix of \cite{Bonforte_2019}. We also take a smooth approximation $\zeta_i(\tilde t,x)\xrightarrow[]{j\to\infty}H(\Phi_1(\tilde t,x)-\Phi_2(\tilde t,x))$.}, i.e. for those solutions such that the equation holds almost everywhere, i.e. for which $\partial_t \Phi = \mathcal{L}_\aa \Phi^m \in \LL^1_{\rm loc}((0,\infty)\times\RR^N)$. Notice that for strong solutions we always have that
\begin{equation*}
\begin{split}
\int (\Phi_1(t,x)-\Phi_2(t,x))_+ \psi(|x|) |x|^{-\aa}\dx &- \int (\Phi_1(0,x)-\Phi_2(0, x))_+ \psi(|x|)|x|^{-\aa}\dx\\
&= \int_{0}^{t}\frac{\rd}{\rd \tilde t}\int (\Phi_1(\tilde t,x)-\Phi_2(\tilde t,x))_+ \psi(|x|)|x|^{-\aa}\dx \rd \tilde t\,,
\end{split}\end{equation*}
hence it only remains to prove that the part inside the time integral and showing that is non-positive:
\begin{align*}
\frac{\rd}{\dt}\int (\Phi_1(t,x)-\Phi_2(t,x))_+ \psi(|x|) |x|^{-\aa}\dx
&= \int   H(\Phi_1(t,x)-\Phi_2(t,x)) \mathcal{L}_\aa (\Phi_1^m(t,x)-\Phi_2^m(t,x)) \psi(|x|)|x|^{-\aa}\dx\\
&= \int   H(\Phi_1^m(t,x)-\Phi_2^m(t,x)) \mathcal{L}_\aa (\Phi_1^m(t,x)-\Phi_2^m(t,x)) \psi(|x|)|x|^{-\aa}\dx\notag \\
&\le \int    \mathcal{L}_\aa (\Phi_1^m(t,x)-\Phi_2^m(t,x))_+ \psi(|x|)|x|^{-\aa}\dx\notag\\
&=\int    (\Phi_1^m(t,x)-\Phi_2^m(t,x))_+ \left(\mathcal{L}_\aa \psi(|x|)\right)|x|^{-\aa}\dx\le 0\notag
\end{align*}
The proof for the general case, i.e., for nonnegative weak solutions, follows by approximation\footnote{
 The rigorous proof can be quite long and technical, but it is standard, as indicated in the previous footnote, hence we have decided to omit it: it relies on careful choices of admissible test functions, see for instance the Appendix of \cite{Bonforte_2019}.}, or by noticing that weak regular bounded solutions are indeed strong, following the ideas in \cite{Lad,Vazquez2007}. This concludes the proof of the comparison for solutions to the FDE \eqref{FDE-CKN-aa}, hence the proof of Claim~2. \qed

\smallskip
It only remains  to prove the validity of condition $(iv)$ and inequality~\eqref{gradient.assumption} of Theorem~\ref{theo:meta} under the additional hypothesis~\eqref{additional.assumption}.

\smallskip

\noindent\textit{Proof that condition $(iv)$ of Theorem~\ref{theo:meta1} holds under the additional hypothesis~\eqref{additional.assumption}. }The proof of the uniform regularity estimates for radial functions follows by extending regularity estimates for the FDE with CKN-weights to the present (easier since radial) case: here for some values of $\aa$ we fall out the ``classical'' CKN-setting of \cite{Bonforte_2019,Bonforte2020}, since a priori we do not have the corresponding CKN inequalities that allow to perform the Nash-Moser iteration. However, since we have opposite powers inside and outside the divergence, and we are in the radial setting, we luckily have ``the right'' weighted Sobolev type inequality valid for all $n>2$ that reads:
\[
\left(\int_0^\infty |f(r)|^{\frac{2n}{n-2}} r^{n-1}\rd r\right)^{\frac{n}{n-2}}\le \mathsf c_n\int_0^\infty |f'(r)|^2 r^{n-1}\rd r\,.
\]
See \cite[Section 1.2.1.3]{Bonforte2020a} for a proof, together with an explicit expression of $\mathsf c_n>0$ and a proof that $\psi(r)$ of the previous step satisfies equality. This inequality allows the local methods of \cite{Bonforte_2019} to work and provide uniform $C^\alpha$ estimates over a unit space-time cylinder. Also we notice that a GHP holds true, as a consequence of our assumptions: when $m\in (m_c,1)$ they allow to use the results of \cite{Bonforte2020}, while when $m\in(0,m_c]$, we they hold simply by comparison, i.e. Claim~2. These two ingredients can be combined as in \cite[Lemma 11]{Bonforte2017b} and allow to prove the desired uniform $C^\alpha$ estimates. Undoing the change of function and going back to solution of the p-Laplacian evolution, we deduce the desired $C^{1,\alpha}$ estimates. \qed

\smallskip
\noindent\textit{Proof of inequality~\eqref{gradient.assumption} of Theorem~\ref{theo:meta} under the additional hypothesis~\eqref{additional.assumption}. }
We notice that it is a straightforward application of Claim~2 in the beginning of the present section. \qed

\section{Comments, Extensions, and Open Problems}\label{sec:open}
In this paper we presented several results on the long-time behaviour for solutions to~\eqref{CPLE}. Let us summarize our results and open problems in the view of Questions~\eqref{Q1}, \eqref{Q2}, and \eqref{Q3} from Introduction. They open several directions in which we our study might be extended.

\begin{enumerate}[{\it (i)}]
\item   Theorem~\ref{theo:RECR} and Proposition~\ref{convergence.L1norm} give convergence rate towards the Barenblatt solution. It is clear, from the examples in Introduction, that our rates are not always sharp. It is, therefore, an interesting open problem to obtain the optimal rates.

\item In the range $p_c<p<p_M$, the entropy method requires an additional assumption (i.e. {\it (ii)} of Theorem~\ref{theo:RECR}) in order to give convergence rates (both for the $\LL^1$-norm and the uniform relative error). It is known, however, that when $u_0\in\X$, solutions still converge to the Barenblatt in the uniform relative error. We pose a question: is it possible to obtain the convergence rate without the additional assumption {\it (ii)} of Theorem~\ref{theo:RECR}? This does not seem an easy task. It was done in the case of the fast diffusion equation by exploiting the (very good) regularity properties of  solutions in that case, see~\cite{Denzler2015} and the shortest version~\cite{Denzler2016}. We notice, however, that solutions to~\eqref{CPLE} do not enjoy the same regularity properties.

\item Let us have a closer look on the convergence results of Theorem~\ref{Thm:gradient.convergence}. It is unclear how to extend the convergence result for radial derivatives to the non-radial case. We shall expect the following. For an intial datum $u_0\in\X$, and for $|x|/t$ large enough, the gradient of a solution to~\eqref{CPLE} behaves as $|\nabla u(t,x)| \sim t^\frac{1}{2-p} |x|^{-\frac{2}{2-p}}$. Therefore we propose the following question: prove or disprove that, when $p_c<p<2$, for an initial datum $u_0\in\X$ with mass $M=\int_{\rn}u_0\dx$, we have that
\[
\lim_{t\rightarrow\infty}\,\, t^{-\frac{1}{2-p}}\,\left\|\left(\nabla u(t)- \nabla \B_M(t)\right)\left(1+|x|^\frac{2}{2-p}\right) \right\|_{\LL^\infty(\rn)}=0\,.
\]
Of course, the same question should be asked in $1<p<p_c$ for solutions expected to converge to the pseudo-Barenblatt profile.
\item As a partial answer to~\eqref{Q3} Theorem~\ref{theo:meta} provides sufficient conditions that has to be satisfied by solutions along time so that the entropy methods work. What is the full description of the basin of attraction of the Barenblatt solutions for $p$ satisfying $1<N<\frac{p}{(2-p)(p-1)}$?  The most interesting information would be giving explicit convergence rates in the relative error under conditions imposed on the initial data only.
\item Despite $p_D$ used to be treated as an important threshold in the analysis of $p$-Laplace Cauchy problem (see Section~\ref{ssec:conv-above-p_D}), we have shown that is only  a technical one restricting the use of the optimal transportation approach, not the dynamics itself. Are the special values we apply: $p_c$, $p_M$, and $p_Y$ essential or technical thresholds?
\end{enumerate}
Lastly, let us comment on two very natural directions that may arise after the present work: the doubly nonlinear equation and anisotropic $p$-Laplace evolution equation. By the doubly nonlinear diffusion equation we mean $\partial_t u= \Delta_p(u^m)$. The fast diffusion regime is when $p(m-1)<1$. It is known, at least in the corresponding \emph{good diffusion range}, that (non-negative and integrable) solutions to the Cauchy problem behave for large times as the corresponding Barenblatt profiles, see for instance~\cite{Agueh2003,Agueh2008, Agueh2009}. Of course, the very natural question is how much of what has been proven in this work also applies to doubly nonlinear case. We believe that the available regularity theory, see for instance~\cite{Bogelein:2023aa,BoDuLi,BoDuSch,LiSch}, allows to try to address questions~\eqref{Q1}, \eqref{Q2}, and~\eqref{Q3}.

The second direction that we believe it is natural to explore are equations of the form $\partial_t u=\sum_{i=1}^N \partial_i\left( |\partial_i u|^{p_i-2} \nabla u\right)$ for possibly different values of $p_i\in\left(1, 2\right)$, $i=1,\dots,N$. In these models, there are several difficulties, starting from the regularity theory to the existence and uniqueness of a fundamental solution. Nevertheless, these models seem to have attracted more and more attention, cf.~\cite {Feo_2021, Feo_2023,Vaz2023singular}. In our analysis, the main difficulty would be understanding the right behaviour for large enough $|x|$. It is unclear whether a class as $\X$ can be found. At the same time, an interesting challenge is adapting the entropy method to those models as the fundamental solution, when it exists, is not explicit.

\section*{Acknowledgements} M.B. acknowledges the financial support from the Spanish Ministry of Science and Innovation, through the Projects PID2020-113596GB-I00 and PID2023-150166NB-I00 and through the ``Severo Ochoa Programme for Centres of Excellence in R\&D'' (CEX2019-000904-S and CEX2023-001347-S).

This work has been (partially) supported by the Project Conviviality ANR-23-CE40-0003 of the French National Research Agency. The authors are glad to acknowledge  funding from Excellence Initiative -- Research University (IDUB) at the University of Warsaw.

The authors are deeply grateful to Professor Jean Dolbeault for bringing to their attention the reference~\cite{Del_Pino_2003}, which had not been included in earlier versions of this manuscript.  
\\

\smallskip
\noindent{\scriptsize\copyright\,2025\ by the authors. This paper may be reproduced, in its entirety, for non-commercial purposes. \href{https://creativecommons.org/licenses/by/4.0/legalcode}{CC-BY 4.0}}

\section*{Appendix}
\renewcommand\thesection{\Alph{section}}
\setcounter{section}{1}
\setcounter{coro}{0}

\subsection{Lemmata}
We present here some  general facts that do not rely strongly on our setting.\newline

Let us present a  Csisz\'ar--Kullback-type inequality yielding that the relative entropy $\cE[v|\VD]$ with respect to the Barenblatt profile of the same mass as $v$ controls the $\LL^1$-distance to the Barenblatt profile. The proof we shall give is inspired from~\cite[Lemma 2.12]{Bonforte2020a}, see also~\cite{Carrillo-Jungel-Markowich-Toscani-Unterreiter} for a previous contribution and more information about this inequality.
\begin{lem}[Csisz\'ar--Kullback inequality]\label{CK.inq} Let $1<p<2$ and $ v:\rn\rightarrow\left[0,\infty\right)$ be a measurable function. Suppose that there exists $D>0$ such that
\[
v-\VD\in\LL^1(\rn)\,,\quad\int_\rn (v-\VD)\dy=0\,,\quad\mbox{and}\quad \cE[v|\VD]<\infty\,.
\]
Then the following inequality holds true
    \[\|v-\VD\|_{\LL^1(\rn)}^2\leq 8\|\V_D^{2-\gamma}\|_{\LL^1(\rn)}\, \cE[v|\VD]\,.\]
\end{lem}\begin{proof}
By the Mean Value Theorem, we know that for $0\leq t\leq s$ it holds
\begin{equation}\label{MVT}
t^{\gamma}-s^{\gamma}-\gamma{\tfrac{1}{\gamma-1}}s^{\gamma-1}(t-s)=\tfrac{\gamma(\gamma-1)}{2}\xi^{\gamma-2}(t-s)^2,\quad\text{with some }\ \xi\in[t,s]\,.
\end{equation}
Since $\xi\leq s$ and $\gamma-2\leq 0$, we infer that
\begin{equation}
    \label{Langrange}
s-t\leq \sqrt{2s^{2-\gamma}}\sqrt{\tfrac{1}{\gamma(\gamma-1)}t^{\gamma}-\tfrac{1}{\gamma(\gamma-1)}s^{\gamma}-\tfrac{1}{\gamma-1}s^{\gamma-1}(t-s)}\,.
\end{equation}
From the assumption $\int_\rn (v-\VD)\dy=0$ we deduce that $\int_{\{v\leq \VD\}} (\VD-v)\dy=\int_{\{\VD\leq v\}}(v-\VD)\dy$, and hence
\[
\frac{1}{2}\int_\rn |v-\VD|\dy=\frac12\left(\int_{\{v\leq \VD\}} (\VD-v)\dy+ \int_{\{\VD\leq v\}}(v-\VD)\dy\right)= \int_{\{v\leq \VD\}} (\VD-v)\dy\,.
\]
Therefore, recalling the very definition of $\cE$ and using inequality~\eqref{Langrange} with the Cauchy-Schwarz inequality we find
\begin{align*}
    \frac{1}{4}\left(\int_\rn |v(\tau)-\VD|\dy\right)^2&\leq {\| 2 \V_D^{2-\gamma}\|_{\LL^1(\rn)}}\ {\cE[v(\tau)|\VD]}\,.
\end{align*}
\end{proof}

Let us establish the decay of functions trapped between two Barenblatt profiles. Recall that $\VD$ is defined in~\eqref{vD}.
\begin{lem}\label{decay.lemma} Let  $p\in(1,2)$,  $D_1\geq D_2>0$ and let $v:\rn\rightarrow\rp$ be a measurable function.
\begin{itemize}
\item[(i)] If for any $|y|\geq 1$ it holds,
 \[
 \V_{D_1}(y)\leq v(y)\leq \V_{D_2}(y)\,,\]
then, for any $D\geq 0$, there exists a constant $C_1=C_1(D, D_1, D_2, p)>0$ such that
\begin{equation}\label{better-decay}
|v(y)-\VD(y)|\leq C_1 |y|^{-\frac{p}{(p-1)(2-p)}}\qquad\text{for }\ |y|\geq 1\,.\end{equation}
\item[(ii)] If $v(y)=v(r)$, where $r=|y|$, and for any $r\geq 1$ it holds,
 \[
 \partial_r\V_{D_2}(r)\leq \partial_r v(r)\leq \partial_r\V_{D_1}(r)\,,\]
 then, for any $D\geq 0$, there exists a constant $C_2=C_2(D, D_1, D_2, p)>0$ and $C_3=C_3(D, D_1, D_2, p)$ such that
 \[\left|\partial_r v(r)-\partial_r\VD(r)\right|\leq C_2 r^{-\frac{p(3-p)}{(p-1)(2-p)}}\qquad\text{for }\ r\geq 1\,,\]
 and
 \[
   \left|\left|\partial_r v(r)\right|^{p-1}-\left|\partial_r \VD(r)\right|^{p-1}\right|\leq C_3\, r^{-\frac{p}{(p-1)(2-p)}+1}
 \]
\end{itemize}
\end{lem}
\begin{proof} Let us begin with~\textit{(i)}. Notice that
\begin{align*}
    \frac{\partial}{\partial D}\VD(y)&=\frac{\partial}{\partial D}\left[\left(D+\frac{2-p}{p}\left|y\right|^\frac{p}{p-1} \right)^\frac{p-1}{p-2}\right]=\frac{-\frac{p-1}{2-p}\left(D+\frac{2-p}{p}\left|y\right|^\frac{p}{p-1} \right)^{\frac{p-1}{p-2}-1}}{\left(D+\frac{2-p}{p}\left|y\right|^\frac{p}{p-1} \right)^{2\frac{p-1}{p-2}}}
=-\frac{p-1}{2-p} \V_D^\frac{1}{p-1}(y)\,.
\end{align*}
Since $D_1>D_2>0$ we can write, for any $y\in\rn$
\[\begin{split}
0\leq \V_{D_2}(y)-\V_{D_1}(y)&=\int_{D_1}^{D_2} \frac{\partial}{\partial D}\V_{D}(y)\dD=-\frac{p-1}{2-p} \int_{D_1}^{D_2} \V_{D}^\frac{1}{p-1}(y)\dD\leq \frac{p-1}{2-p}|D_2-D_1|\, \V_{D_1}^\frac{1}{p-1}(y)\,,
\end{split}\]
which is integrable for the prescribed range of $p$. The lower bound can be shown in the same way. Therefore
\[\tfrac{p-1}{2-p} |D_2-D_1|\, \V_{D_2}^\frac{1}{p-1}(y)\leq \V_{D_2}(y)-\V_{D_1}(y)\leq \tfrac{p-1}{2-p}|D_2-D_1|\, \V_{D_1}^\frac{1}{p-1}(y)\,.\]
In turn, for any $D$, we can estimate
\begin{equation*}\begin{split}|v(y)- \VD(y)|&\leq |\V_{D_2}(y)- v(y)|+|\VD(y)- \V_{D_2}(y)|\leq  (\V_{D_2}(y)- \V_{D_1}(y))+|\VD(y)- \V_{D_2}(y)|\\
&\leq \tfrac{p-1}{2-p} (D_2-D_1+|D-D_2|) \V_{{D_0}}^\frac{1}{p-1}(y)\,,
\end{split}\end{equation*}
where ${D_0}=\min\{D,D_1,D_2\}$. By taking into account that for $|y|\geq 1$ it holds $\V_{{D_0}}(y)\leq C|y|^{-\frac{p}{2-p}}$ we get the claim.
In the case of~\textit{(ii)} we notice that
\[
  \partial_r \VD(r)=-r^\frac{1}{p-1}\,\V_{D}^\frac{1}{p-1}(r)\,.
\]
By using this observation in a similar way as above we get the claim for~\textit{ii)}.
\end{proof}
\begin{lem}[Lemma~3.1,~\cite{Borowski2022}]\label{lem:p-male} Suppose $1<p< 2$. Then there exist $c_1,c_2>0$ such that for all $\xi,\eta\in\rn$ such that $\xi\neq 0$ we have \begin{equation*}
\left\langle |\xi|^{p-2}\xi-|\eta|^{p-2}\eta,\xi-\eta\right\rangle\geq c_1\frac{|\xi-\eta|^2}{|\xi|^{2-p}+|\eta|^{2-p}}\,,
\end{equation*}
where the optimal constant is achieved when $\langle\xi,\eta\rangle =|\xi|\,|\eta|$ and is given by $c_1 = \min\{1, 2(p - 1)\}$, and
\begin{equation*}
    \left\langle |\xi|^{p-2}\xi-|\eta|^{p-2}\eta,\xi-\eta\right\rangle\geq c_2 \frac{|\xi-\eta|^2}{(|\xi|+|\eta|)^{2-p}}\,.
\end{equation*}
where $c_2=c_1/2$.
\end{lem}

\begin{lem}\label{numerical_lemma_fisher}
  Let $0\leq \xi, \eta\in\mathbb{R}$ and $1<p<2$, then
  \begin{equation}\label{numerical.inequality.fisher}
    \max\{\xi, \eta\}^{p-2}\,|\xi-\eta|^2 \leq \frac{\xi^p-\xi^{p-1}\eta - \eta^{p-1}\xi-\eta^p}{p-1} \leq \min\{\xi, \eta\}^{p-2}\,|\xi-\eta|^2
  \end{equation}
\end{lem}
\begin{proof}
  Let us call $F(\xi, \eta)=\xi^p-\xi^{p-1}\eta - \eta^{p-1}\xi-\eta^p$, we notice that
  \[
    F(\xi, \eta)=(\xi^{p-1}-\eta^{p-1})(\xi-\eta)\ge0\quad\forall\,\,\xi,\eta\ge0.
  \]
Since $F(\xi, \eta)=F(\eta, \xi)\ge0$, we can assume, without loss of generality, that $\xi\ge\eta$. If $\eta=0$, there is nothing to prove in~\eqref{numerical.inequality.fisher}, so let us suppose that $\eta>0$. From the concavity of the function $t\mapsto t^{p-1}$ we deduce that, for all $t\ge1$, we have $t^{p-1}\leq 1+(p-1)(t-1)$. By applying this last inequality to $t=\xi/\eta\ge1$, we obtain that
\[
  \xi^{p-1}-\eta^{p-1}=\eta^{p-1}\left(t^{p-1}-1\right)\leq (p-1)\eta^{p-1}(t-1)=(p-1)\eta^{p-2}(\xi-\eta)\,,
\]
from which the right-hand-side of~\eqref{numerical.inequality.fisher} could be easily deduced. The left-hand-side in~\eqref{numerical.inequality.fisher} is deduced analogously from the convexity of the function $s\mapsto -s^{p-1}$. The proof is concluded.
\end{proof}

We give here a modified version of the Gronwall-type lemma. We are sure it is known, but since we were not able to find a relevant reference, we present it with a proof.
\begin{lem}\label{weak.gronwall}
Let $u:[0, \infty)\rightarrow [0, \infty)$ be bounded, decreasing and satisfying the inequality
\[
s\,\int_{t}^\infty u(\tau)\dtau \le u(t)\quad\forall\, t>0\,,
\]
where $s>0$. Then, there exists $C=C(\uo, s)>0$ such that  $u(t)\le C\,e^{-s\,t}$ for all $t>0$. 
From the proof is clear that $C(\uo,s)=\frac{e^s}{s}\,u(0)$.
\end{lem}
\begin{proof}
Let us define
\[
v(t):=\int_t^\infty u(\tau)\dtau\,.
\]
By the properties of $u$, we infer that $v\in W^{1,\infty}(0,\infty)$. We may apply the classical version of the Gronwall lemma to $v$, since it satisfies the inequality
\[
v'(t)=-u(t)\leq -s \int_t^\infty u(\tau)\dtau = -s\, v(\tau)\,,
\]
and thus $v(t)\le v(0) e^{-s\, t}$. Since, by hypothesis, we have $s\,v(0)\le u(0)$, we shall find $v(t)\le v(0) e^{-s\, t}\le \frac{\uo(0)}{s}\,e^{-s\, t}\,.$ Notice that, by definition of $v(t)$, we have
\[
v(t-1)-v(t)=\int_{t-1}^t u(\tau)\dtau\ge0\,.
\]
At the same time, since $u$ is nonincreasing, i.e., $u(t)\le u(s)$ for any $s\in[t-1,t]$, we have that
\[
u(t)=\int_{t-1}^t u(t)\dtau \le \int_{t-1}^t u(\tau)\dtau=v(t-1)-v(t)\le v(t-1)\,.
\]
Combining all the above estimates we deduce that $u(t)\le \frac{u(0)}{s}\,e^{-s(t-1)}\,.$
\end{proof}
\begin{lem}\label{lem:improved-rate}
Let $\lambda\in\left(0,\infty\right)$  and $u, R:\left[0, \infty\right)\rightarrow \mathbb{R}$ be differentiable functions and that satifies the inequality
\[
  u'(t)\leq -\lambda + R(t)\quad\forall t\in\left(0, \infty\right)\,,
\]
such that $I=\int_0^\infty R(\tau) \dtau<\infty$. Then, for any $t_2>t_1>0$ we have that
\[
  u(t_2)\leq u(t_1) -\lambda (t_2-t_1) + I\,.
\]
\end{lem}
\begin{proof}
The proof is a simple application of the Gronwall lemma.
\end{proof}

\subsection{Parameters}

\noindent {\bf Thresholds for $p$.} Their role is described in more details in Introduction.

\begin{tabular}{l|lll}
symbol & introduced & info\\\hline
$p_c=\frac{2N}{N+1}\in(\frac{3}{2},2)$ &\eqref{p_c}& for $p>p_c$ solutions to \eqref{CPLE} conserve mass\\
$p_M\in(p_c,2)$ &\eqref{p_M}& for $p>p_M$ solutions to \eqref{CPLE} have finite weighted $|x|^{p'}$-moments\\
$p_1,p_2$ &\eqref{p1p2}& integrability threshold for $|\B_{M_1}-\B_{M_2}|$ defined if $N\geq 6$;\\ &&$N<\frac{p}{(2-p)(p-1)}\iff p\in(1,p_1)\cup (p_2,2)$\\
$p_Y=\frac{2N}{N+2}$&\eqref{pY}& Yamabe exponent and gradient regularity threshold, cf. Section~\ref{ssec:reg}\\
$p_D=\frac{2N+1}{N+1}\in (p_M,2)$ &\eqref{pD}& for $p>p_D$ the entropy functional is displacement convex, cf. Introduction\\
\end{tabular}

\bigskip

\noindent {\bf Main characters}

\begin{tabular}{l|lll}
symbol & introduced & info\\\hline
$u$ & \eqref{CPLE} & a solution to $p$-Laplace Cauchy problem with $\uo$ as initial datum; proven to converge\\
&& to $\B_M$ \\
$v$ & \eqref{eq:v} & a solution to Nonlinear Fokker--Planck problem with $\vo$ as initial datum\\
$\Phi$ &  \eqref{fde.cauchy.problem} & a solution to a radial FDE problem
\end{tabular}

\medskip

\noindent{\bf Other symbols}

\begin{tabular}{l|lll}
symbol & introduced & info\\\hline
$\VD$ &\eqref{vD} & stationary solution to the Fokker--Planck equation~\eqref{eq:v}\\
$\beta$ &\eqref{beta} & parameter for definition of Barenblatt profile; $\beta(p-p_c)\geq 0$ \\
$b_1,b_2$&\eqref{b1}, \eqref{b2} & parameter for definition of Barenblatt profile\\
$\ell$ &\eqref{time.rescaling-pc} & free parameter for definition of Barenblatt profile when $p=p_c$ \\
 $R_T(t)$ &\eqref{R(t)} or \eqref{time.rescaling-pc} & time rescaling for definition of Barenblatt profile \\
&&
$R_T(t)=\begin{cases}
    \{(T-t)_+/|\beta|\}^\beta\,, & 1<p<p_c\,,\\
    \exp\{\ell (T+t)\}\,,& p=p_c\,,\\
    \{(T+t)/|\beta|\}^\beta\,, & p_c<p<2\,.\\
\end{cases}$
\\
$\B_M(t+\beta,x)=R^{-N}_\beta(t)\VD(y)$ & \eqref{Barenblatt-via-vD} & the Bareblatt solution to~\eqref{CPLE} with $D$ as in \eqref{Barenblatt-via-vD}; for $p>p_c$ \\
$\B_{{D},T}(t+\beta,x)=R_T^{-N}(t)\VD(y)$ & \eqref{Barenblatt-via-vD} & the Bareblatt solution to~\eqref{CPLE}; for $p\leq p_c$ \\
$M_\star$ &\eqref{def:Mstar}& $M_\star=\|V_1\|_{\LL^1}$\\
$\gamma=\frac{2p-3}{p-1}$ &\eqref{cal-E} & parameter of the entropy\\
$\cE$ &\eqref{cal-E} & relative entropy functional\\
$\cI$ &\eqref{cal-I-def} & relative Fisher information\\
$m=p-1\in (0,1)$ &\eqref{parameter.transformation} & exponent for the radial FDE\\
$n=2+2\tfrac{N}{p'}> 2$ & \eqref{parameter.transformation} & artificial dimension for the radial FDE\\
$\aa=N-2-2\frac{N}{p'}$ &\eqref{parameter.transformation} & $\aa=N-n$ \\
$a_1,a_2>0$&\eqref{a1-a2} & depend on $m,n,\aa$ or $p,N$\\
$\theta>0$&\eqref{bareblatt.profile.ckn} &  FDE-Barenblatt parameter\\
$\mathfrak{B}_M$&\eqref{bareblatt.profile.ckn} &  the Bareblatt solution to the radial FDE~\eqref{ckn.fde.radial}
\end{tabular}

\bibliographystyle{abbrv}
\bibliography{reference}

\end{document}